\documentclass[11pt]{article}
\usepackage{amssymb,amsmath,mathrsfs,txfonts,graphicx,tikz,color}
\usepackage[title]{appendix}
\graphicspath{{figs/}}

\allowdisplaybreaks

\let\oldsection\section
\renewcommand\section{\setcounter{equation}{0}\oldsection}

\newtheorem{theorem}{\indent Theorem}[section]
\newtheorem{lemma}{\indent Lemma}[section]
\newtheorem{proposition}{\indent Proposition}[section]

\newtheorem{remark}{\indent Remark}[section]
\newtheorem{corollary}{\indent Corollary}[section]

\begin{document}

\title{\LARGE
Optimal decay rates of the compressible Euler equations
with time-dependent damping in $\mathbb R^n$: (II) over-damping case}
\author{
Shanming Ji$^{a,c}$,
Ming Mei$^{b,c,}$\thanks{Corresponding author, email:ming.mei@mcgill.ca}
\\
\\
{ \small \it $^a$School of Mathematics, South China University of Technology}
\\
{ \small \it Guangzhou, Guangdong, 510641, P.~R.~China}
\\
{ \small \it $^b$Department of Mathematics, Champlain College Saint-Lambert}
\\
{ \small \it Quebec, J4P 3P2, Canada, and}
\\
{ \small \it $^c$Department of Mathematics and Statistics, McGill University}
\\
{ \small \it Montreal, Quebec, H3A 2K6, Canada}
}
\date{}

\maketitle

\begin{abstract}
This paper is concerned with the multi-dimensional compressible Euler equations
with time-dependent over-damping of the form $-\frac{\mu}{(1+t)^\lambda}\rho\boldsymbol u$
in $\mathbb R^n$, where $n\ge2$, $\mu>0$, and $\lambda\in[-1,0)$.
This continues our previous work dealing with the under-damping case for $\lambda\in[0,1)$.
We show the optimal decay estimates of the solutions such that
for $\lambda\in(-1,0)$ and $n\ge2$,
$\|\rho-1\|_{L^2(\mathbb R^n)}\approx(1+t)^{-\frac{1+\lambda}{4}n}$
and $\|\boldsymbol u\|_{L^2(\mathbb R^n)}\approx
(1+t)^{-\frac{1+\lambda}{4}n-\frac{1-\lambda}{2}}$,
which indicates that a stronger damping gives rise to solutions decaying optimally slower.
For the critical case of $\lambda=-1$, we prove
the optimal logarithmical decay of the perturbation of density
for the damped Euler equations such that
$\|\rho-1\|_{L^2(\mathbb R^n)}\approx |\ln(e+t)|^{-\frac{n}{4}}$
and $\|\boldsymbol u\|_{L^2(\mathbb R^n)}\approx
(1+t)^{-1}\cdot|\ln(e+t)|^{-\frac{n}{4}-\frac{1}{2}}$
for $n\ge7$. The over-damping effect reduces the decay rates of the solutions to be slow, which causes us
some technical difficulty in obtaining the optimal decay rates by the Fourier analysis method and the Green function method.
Here, we propose a new idea to overcome such a difficulty by artfully combining the Green function method and the time-weighted energy method.
\end{abstract}


\tableofcontents

\baselineskip=15pt

\section{Introduction}
\subsection{Modeling equations and research background}
We consider the multi-dimensional compressible Euler equations with time-dependent damping
\begin{equation} \label{eq-Euler}
\begin{cases}
\displaystyle
\partial_t \rho+\nabla\cdot(\rho \boldsymbol u)=0, \\
\displaystyle
\partial_t (\rho \boldsymbol u)+\nabla\cdot(\rho \boldsymbol u\otimes \boldsymbol u)
+\nabla p(\rho)=-\frac{\mu}{(1+t)^\lambda}\rho \boldsymbol u, \\
\displaystyle
\rho|_{t=0}=\rho_0(x):=1+\tilde\rho_0(x), \quad \boldsymbol u|_{t=0}=\boldsymbol u_0(x),
\end{cases}
\end{equation}
where $x\in \mathbb R^n$, $n\ge2$, $\mu>0$, $\lambda\in[-1,0)$.
Here, the unknown functions $\rho(t,x)$ and $\boldsymbol u(t,x)$
represent the density and velocity of the fluid,
and the pressure $p(\rho)=\frac{1}{\gamma}\rho^\gamma$ with $\gamma>1$.
The initial data satisfy
\begin{equation}\label{newnew-2020}
\rho_0(x) \to 1, \ \ \mbox{i.e.,} \ \
\tilde\rho_0(x) \to 0, \ \mbox{ and } \
\boldsymbol u_0(x) \to \boldsymbol 0 , \ \ \mbox{ as } |x|\to \infty.
\end{equation}
The under-damping case of $\lambda\in[0,1)$ is considered in the first part \cite{Ji-Mei-1} of our series of study,
where we shown that weaker damping leads to faster decays.
Here in this paper, we focus on the over-damping case of $\lambda\in[-1,0)$
and we prove that stronger damping gives rise to optimally slower decays.

As we mentioned in the first part \cite{Ji-Mei-1} of this series of study, the damping effect plays a key role in the structure of solutions to the compressible Euler equations.
Without damping effect, the solutions of Euler equations  usually possess singularity like shock waves and exhibit blow-up for their gradients \cite{C-D-S-W,Courant-F,Dafermos,Lax,Smoller}.
When the Euler system of equations are with damping effect, the structure of the solutions becomes more complicated and various according to the size of the damping effect, and of course, the study is more challenging.
When $\lambda=0$ and $\mu>0$, the regular case of damping effect in the form of
$-\mu\rho \boldsymbol u$, once the initial data and their gradients are small enough,
the damping effect can prevent the formation of shocks for the damped Euler equations \cite{Sideris-Thomas-Wang}, and makes the solutions to behave time-asymptotically as the so-called diffusion waves for the corresponding nonlinear diffusion (porous media) equations \cite{Hsiao-Liu,Marcati-Milani,Mei,Nishihara,Nishihara-Wang-Yang}; while, once the gradients of the initial data are bigger, the blow-up phenomena for the solutions of Euler equations with regular damping still occur \cite{Hailiang-Li,Wang-Chen}.
When $\lambda>0$ and $\mu>0$, the damping effect $-\frac{\mu}{(1+t)^\lambda}\rho \boldsymbol u$ becomes weaker as $\lambda$ increases, the so-called under-damping case.
Here, for $\lambda\in (0,1)$ and $\mu>0$, once the initial data and their gradients are  small enough, the weak damping effect can still guarantee the global existence of the  solutions for the Euler equations with under-damping \cite{Chen-Li-Li-Mei-Zhang,Cui-Yin-Zhang-Zhu,Hou-Witt-Yin,Hou-Yin,Li-Li-Mei-Zhang,Pa1,Sugiyama1}; while the solutions will blow up at finite time when the gradients of the initial data are big \cite{Chen-Li-Li-Mei-Zhang}.
However, when $\lambda>1$ with $\mu>0$, the damping effect is too weak, and the Euler system with such a weak damping essentially behaves like the pure Euler system so that the singularity of shocks cannot be avoided, no matter how smooth and small the initial data are \cite{Chen-Li-Li-Mei-Zhang,Hou-Witt-Yin,Hou-Yin,Pa2,Sugiyama2}.
Such blow-up phenomena in this super under-damping case of $\lambda>1$ are determined by the mechanism of the dynamic system itself, rather than the selection of the initial data \cite{Chen-Li-Li-Mei-Zhang}.
When $\lambda=1$, this is the critical case, where the solutions globally exist for $\mu>3-n$ as shown in \cite{Hou-Witt-Yin,Hou-Yin}
(see also \cite{Chen-Li-Li-Mei-Zhang,Geng-Lin-Mei,Pa2,Sugiyama2} for $1$-D case)
and occur blow-up for $\mu\le 3-n$  as studied in \cite{Hou-Witt-Yin,Hou-Yin}.

For the global solutions of the dynamic system of partial differential equations, one of the fundamental problems from both mathematical and physic points of view is to investigate the asymptotic behavior at large-time.
For the time-dependent damped Euler equations \eqref{eq-Euler}, when $\lambda=0$, the optimal decay rates were technically obtained by Sideris-Thomas-Wang \cite{Sideris-Thomas-Wang} when the initial data are in certain Sobolev space
and by Tan {\it et al.} \cite{TanZ-JDE13,TanZ-JDE12} in some Besov spaces.
For $\lambda\in (0,1)$, the methods for deriving the decay estimates of the solutions adopted in the previous studies for $\lambda=0$ case in \cite{Sideris-Thomas-Wang,TanZ-JDE13,TanZ-JDE12} cannot be directly applied, due to the complexity of the damping effect involving the time $t$.
In our study \cite{Ji-Mei-1}, we apply the technical Fourier analysis to derive the optimal decay estimates for the linearized system which can be formally expressed by the implicit Green functions, then use the weighted-energy method with some
new  developments to obtain the optimal decay rates of the solutions for the nonlinear Euler equations with time-dependent under-damping:
\[
\|\partial_x^\alpha (\rho-1)\|_{L^2(\mathbb R^n)}\approx(1+t)^{-\frac{1+\lambda}{2}(\frac{n}{2}+|\alpha|)},
\quad \|\partial_x^\alpha \boldsymbol u\|_{L^2(\mathbb R^n)}\approx
(1+t)^{-\frac{1+\lambda}{2}(\frac{n}{2}+|\alpha|)-\frac{1-\lambda}{2}},
\qquad \lambda\in[0,1).
\]
The new point observed in \cite{Ji-Mei-1} is that, for $\lambda\in [0,1)$, the weaker under-damping effect makes the faster decay of the solutions,
namely, the decay of the solutions at $\lambda=0$ is weakest, while the decays of the solutions around $\lambda = 1^-$ are much faster.

However, for $\lambda<0$, the so-called over-damping case, the relevant study for the damped Euler equations is almost nothing,  to the best of our knowledge.
This will be the main concern of the present paper.
We consider the case for $\lambda\in [-1,0)$ and $\mu>0$.
First of all, we focus ourselves on the case of $\lambda \in (-1,0)$, and show the optimal decay of the implicit Green functions by using  Fourier analysis to the high frequency part and the low frequency part respectively, and further obtain the optimal decay estimates for the solutions to the nonlinear Euler equations with time-dependent over-damping \eqref{eq-Euler} by the Green function method with some restriction on $\lambda$.
That is,
\[
\|\partial_x^\alpha (\rho-1)\|_{L^2(\mathbb R^n)}\approx(1+t)^{-\frac{1+\lambda}{2}(\frac{n}{2}+|\alpha|)},
\quad \|\partial_x^\alpha \boldsymbol u\|_{L^2(\mathbb R^n)}\approx
(1+t)^{-\frac{1+\lambda}{2}(\frac{n}{2}+|\alpha|)-\frac{1-\lambda}{2}},
\qquad \lambda\in(-{\textstyle\frac{n}{n+2}},0),
\]
but we have to restrict $\lambda\in (-\frac{n}{n+2},0)$ due to the bad effect of the over-damping.
In fact, from the above decay estimates, we realize that the over-damping effect for $\lambda\in (-1,0)$ makes the decay of the solutions to get slower and slower, as $\lambda\to -1^+$.
Namely, the strongest over-damping at $\lambda=-1^+$ reduces the solutions decay slowest.
Just because of this, we cannot close the high-order decay estimates for all $\lambda\in (-1,0)$ by the Green function method, and have to leave the case of $\lambda \in (-1, -\frac{n}{n+2}]$ open.
In order to delete such a gap for $\lambda \in (-1, -\frac{n}{n+2}]$,
we propose a new technique, which is an artful combination of the Green function method
and the time-weighted energy method.
The Green function method cannot perfectly treat the high-order decay estimates
for $\lambda$ near $-1$,
and the time-weighted energy method is also short in deriving the optimal decay estimates,
but it is very efficient to treat the high-order estimates.
Hence we try to combine these two methods together to get the optimal decay estimates for all $\lambda\in (-1,0)$.
In fact, the procedure to cleverly combine both existing methods is still technical as we know. Thus, we can finally prove the optimal decay estimates for all
$\lambda\in (-1,0)$ as follows
\[
\|\rho(t,x)-1\|_{L^2(\mathbb R^n)}\approx (1+t)^{-\frac{1+\lambda}{4}n}, \qquad
\|\boldsymbol u(t,x)\|_{L^2(\mathbb R^n)}\approx
(1+t)^{-\frac{1+\lambda}{4}n-\frac{1-\lambda}{2}},
\qquad \lambda\in(-1,0).
\]
Secondly, we consider the critical case of $\lambda=-1$, the most interesting but also the most difficult case. We further show the optimal decay rates as follows
\[
\|\rho(t,x)-1\|_{L^2(\mathbb R^n)}\approx |\ln(e+t)|^{-\frac{n}{4}}, \qquad
\|\boldsymbol u(t,x)\|_{L^2(\mathbb R^n)}\approx
(1+t)^{-1}\cdot|\ln(e+t)|^{-\frac{n}{4}-\frac{1}{2}},
\qquad \lambda=-1.
\]
But we have to restrict the space dimension $n\ge 7$ for some technical reason.

For the other topics with vacuum for the damped Euler equations, we refer to the significant works \cite{Geng-Huang,Huang-Marcati-Pan,Huang-Pan,Huang-Pan-Wang,Luo-Zeng}. For the linear wave equations with time-dependent damping, we refer to the pioneering studies by Wirth in \cite{Wirth-JDE06,Wirth-JDE07,Wirth-MMAS04}. For the time-dependent damped Klein-Gordon equations, we refer to the interesting results by Burq-Raugel-Schlag in  \cite{Burq-Raugel-Schlag-2015,Burq-Raugel-Schlag-2018}.

\subsection{Transformation of equations and notations}

In order to study the system \eqref{eq-Euler}, we switch it to a symmetric system. Let $v=\frac{2}{\gamma-1}(\sqrt{p'(\rho)}-1)
=\frac{2}{\gamma-1}(\rho^\frac{\gamma-1}{2}-1)$ and $\varpi=\frac{\gamma-1}{2}$.
Then $(v,\boldsymbol u)$ satisfies the following symmetric system
\begin{equation} \label{eq-vbdu}
\begin{cases}
\displaystyle
\partial_t v+\nabla\cdot\boldsymbol u
=-\boldsymbol u\cdot\nabla v-\varpi v\nabla\cdot \boldsymbol u,\\
\displaystyle
\partial_t \boldsymbol u+\nabla v+\frac{\mu}{(1+t)^\lambda}\boldsymbol u
=-(\boldsymbol u\cdot\nabla) \boldsymbol u-\varpi v\nabla v,\\
\displaystyle
v|_{t=0}=v_0(x), \quad \boldsymbol u|_{t=0}=\boldsymbol u_0(x),
\end{cases}
\end{equation}
where $v_0(x)=\frac{2}{\gamma-1}((1+\tilde\rho_0(x))^\frac{\gamma-1}{2}-1)$,
which behaves like $\tilde\rho_0(x)$ if the initial perturbation is small.

{\bf Notations.}
We denote $D_t=-i\partial_t$, and $\hat v(\xi)=\mathscr{F}(v)$ the $n$-dimensional Fourier transform
of a function $v(x)$.
We use $H^s=H^s(\mathbb R^n)$, $s\in\mathbb R$, to denote Sobolev spaces, and
$L^p=L^p(\mathbb R^n)$, $1\le p\le \infty$, to denote the $L^p$ spaces.
The spatial derivatives $\partial_x^\alpha$ stands for
$\partial_{x_1}^{\alpha_1}\cdots\partial_{x_n}^{\alpha_n}$
with nonnegative multi-index $\alpha=(\alpha_1,\dots,\alpha_n)$
(the order of $\alpha$ is denoted by $|\alpha|=\sum_{j=1}^{j=n}\alpha_j$)
and $\partial_x^{|\alpha|}$ stands for all the
spatial partial derivatives of order $|\alpha|$.
The pseudo differential operator $\Lambda$ is defined by
$\Lambda^s v:=\mathscr{F}^{-1}(|\xi|^s\hat v(\xi))$ for $s\in\mathbb R$.
The norm $\|v\|_X^l$ stands for the $\|\cdot\|_X$ norm of the low frequency
part $v^l:=\mathscr{F}^{-1}(\chi(\xi)\hat v(\xi))$ of $v$,
while $\|v\|_X^h$ stands for the $\|\cdot\|_X$ norm of the high frequency
part $v^h:=\mathscr{F}^{-1}((1-\chi(\xi))\hat v(\xi))$ of $v$,
where $0\le\chi(\xi)\le1$ is a smooth cut-off function supported in $B_{2R}(0)$
and $\chi(\xi)\equiv1$ on $B_R(0)$ for a given $R>0$.

Throughout this paper, we
denote $b(t)=\frac{\mu}{(1+t)^\lambda}$ with $\mu>0$ and $\lambda\in[-1,0)$
and we let $C$ (or $C_j$ with $j=1,2,\dots$)
denote some positive universal constants
(may depend on $n$, $\lambda$, $\mu$, $\gamma$, and $\alpha$).
We use $f\lesssim g$ or $g\gtrsim f$ if $f\le Cg$,
and denote $f\approx g$ if $f\lesssim g$ and $g\gtrsim f$.
For simplicity, we define $\|(f,g)\|_X:=\|f\|_X+\|g\|_X$
and $\int f:=\int_{\mathbb R^n}f(x)dx$.
The norm $\|\cdot\|_{L^2}$ will be simplified as $\|\cdot\|$  without confusion.
For a matrix $A=(A_{j,k})$, the norm $\|A\|_{\max}:=\max_{j,k}|A_{j,k}|$
is the maximum absolute value of all its elements.
We define the following time decay function
\begin{equation} \label{eq-Gamma}
\Gamma(t,s):=
\begin{cases}
[1+(1+t)^{1+\lambda}-(1+s)^{1+\lambda}]^{-\frac{1}{2}}, \qquad & \lambda\in(-1,0),
\\[2mm]
\displaystyle
\Big[1+\ln\Big(\frac{1+t}{1+s}\Big)\Big]^{-\frac{1}{2}}, \qquad & \lambda=-1.
\end{cases}
\end{equation}

\subsection{Main results}

For the over-damping case with $\lambda\in[-1,0)$,
our main results for the global existence and uniqueness of the solutions as well as the optimal decay esitmates are stated  as follows.

\begin{theorem}[Optimal $L^2$ decay estimates of nonlinear Euler system] \label{th-nonlinear}
For the dimension $n\ge2$ and $\lambda\in(-\frac{n}{n+2},0)$,
there exists a constant $\varepsilon_0>0$, such that
the solution $(v,\boldsymbol u)$ of the nonlinear system \eqref{eq-vbdu}
corresponding to initial data $(v_0,\boldsymbol u_0)$
with small energy $\|(v_0,\boldsymbol u_0)\|_{L^1\cap H^{[\frac{n}{2}]+3}}\le\varepsilon_0$
exists time-globally and satisfies
\begin{equation} \label{eq-nonlinear}
\begin{cases}
\|\partial_x^\alpha v\|\lesssim (1+t)^{-\frac{1+\lambda}{4}n-\frac{1+\lambda}{2}|\alpha|},
\quad &0\le |\alpha|\le [\frac{n}{2}]+1,\\
\|\partial_x^\alpha \boldsymbol u\|\lesssim
(1+t)^{-\frac{1+\lambda}{4}n-\frac{1+\lambda}{2}(|\alpha|+1)+\lambda},
\quad &0\le |\alpha|\le [\frac{n}{2}],\\
\|(v,\boldsymbol u)\|_{H^{[\frac{n}{2}]+3}}\lesssim 1.
\end{cases}
\end{equation}
The first two decay estimates in \eqref{eq-nonlinear}
(i.e., the decay estimates on $\|\partial_x^\alpha v\|$
with $0\le|\alpha|\le [\frac{n}{2}]+1$ and
$\|\partial_x^\alpha \boldsymbol u\|$ with $0\le|\alpha|\le [\frac{n}{2}]$) are optimal
and consistent with the linearized hyperbolic system.
\end{theorem}

\begin{theorem}[Optimal $L^q$ decay estimates of nonlinear Euler system]
\label{th-nonlinear-Lp}
For $n\ge2$, $q\in[2,\infty]$, $k\ge 3+[\gamma_{2,q}]$
with $\gamma_{2,q}:=n(1/2-1/q)$, and $\lambda\in(-\frac{n}{n+2},0)$,
let $(v,\boldsymbol u)$ be the solution to the nonlinear system \eqref{eq-vbdu}
corresponding to initial data $(v_0,\boldsymbol u_0)$
with small energy such that
$\|(v_0,\boldsymbol u_0)\|_{L^1\cap H^{[\frac{n}{2}]+k}}\le\varepsilon_0$,
where $\varepsilon_0>0$ is a small constant only depending on $n,q,k$ and
the constants $\gamma,\mu,\lambda$ in the system.
Then $(v,\boldsymbol u)\in L^\infty(0,+\infty;H^{[\frac{n}{2}]+k})$ and satisfies
\begin{equation} \label{eq-nonlinear-Lp}
\begin{cases}
\|\partial_x^\alpha v\|_{L^q}\lesssim
(1+t)^{-\frac{1+\lambda}{2}\gamma_{1,q}-\frac{1+\lambda}{2}|\alpha|},
\quad &0\le |\alpha|\le 1,\\
\|\boldsymbol u\|_{L^q}\lesssim
(1+t)^{-\frac{1+\lambda}{2}\gamma_{1,q}-\frac{1-\lambda}{2}},
\end{cases}
\end{equation}
where $\gamma_{1,q}=n(1-1/q)$.
The decay estimates in \eqref{eq-nonlinear-Lp} are optimal.
\end{theorem}

\begin{remark}
The above optimal $L^2$ and $L^q$ decays are formulated by means of the
technical Fourier analysis and the Green function method.
The restriction of $\lambda\in(-\frac{n}{n+2},0)$ comes from the following two main
difficulties caused by the over-damping:

(i) The optimal decay of $\|\partial_x^\alpha v\|$ for the linearized hyperbolic system
of \eqref{eq-vbdu} is slow,
\begin{equation} \label{eq-decay-slow}
\begin{cases}
\|\partial_x^\alpha\mathcal{G}_{11}(t,0)v_0\|\approx
(1+t)^{-\frac{1+\lambda}{4}n-\frac{1+\lambda}{2}|\alpha|},
\quad & \lambda\in(-1,0),
\\
\|\partial_x^\alpha\mathcal{G}_{11}(t,0)v_0\|\approx
|\ln(e+t)|^{-\frac{n}{4}-\frac{|\alpha|}{2}},
\quad & \lambda=-1,
\end{cases}
\end{equation}
where $\mathcal{G}(t,s)$ is the Green matrix (see \eqref{eq-Duhamel}).
One should be careful in calculating the estimates of $\int_0^t \mathcal{G}(t,s)Q(s)ds$
involving general nonlinear terms $Q(t)$.

(ii) The over-damping $b(t)$ causes trouble in the estimates on
$b(t)\partial_x^k \boldsymbol u\cdot \partial_x^{k+1}v$,
which is crucial for the high-order energy estimates on $\|\partial_x^{k+1}v\|$
in the closure of the a priori assumption.
\end{remark}

\begin{remark}
The solutions to the linearized hyperbolic system of \eqref{eq-vbdu} decay optimally slower
for the over-damping case.
We may understand it as follows:
when the over-damping is stronger as $\lambda\in[-1,0)$,
the high frequencies decay faster as $e^{-C(1+t)^{1-\lambda}}$ (super-exponential),
while the low frequencies decay slower as
$$
\begin{cases}
e^{-C|\xi|^2(1+t)^{1+\lambda}}, \quad & \mathrm{for~} \lambda\in(-1,0),
\\
e^{-C|\xi|^2\ln(e+t)}, \quad & \mathrm{for~critical~} \lambda=-1,
\end{cases}
$$
and on the whole the solutions decay slower.
\end{remark}

In order to formulate the decay estimates for all $\lambda\in(-1,0)$ and
especially for the critical case of $\lambda=-1$,
we develop a time-weighted iteration scheme,
which is a combined time-weighted energy estimates and Green functions we build up in the above,
to close the decay estimates.

\begin{theorem}[Optimal decay estimates for $\lambda\in(-1,0)$] \label{th-com-0-in}
For $n\ge2$, $N\ge[\frac{n}{2}]+2$ and $\lambda\in(-1,0)$,
there exists a constant $\varepsilon_0>0$ such that
the solution $(v,\boldsymbol u)$ of the nonlinear system \eqref{eq-vbdu}
corresponding to small initial data
$\|(v_0,\boldsymbol u_0)\|_{L^1\cap H^N}\le \varepsilon_0$
exists globally and satisfies
\begin{equation} \label{eq-com-0-in}
\begin{cases}
\|v(t,\cdot)\|\lesssim (1+t)^{-\frac{1+\lambda}{4}n}, \\
\|\boldsymbol u(t,\cdot)\|\lesssim (1+t)^{-\frac{1+\lambda}{4}n-\frac{1-\lambda}{2}}.
\end{cases}
\end{equation}
The above decay estimates are optimal and consistent with the linearized hyperbolic system.
\end{theorem}

\begin{remark}
Theorem \ref{th-nonlinear} shows the optimal decay rates of all derivatives of solutions
$\|\partial_x^\alpha v\|$ with $0\le|\alpha|\le [\frac{n}{2}]+1$ and
$\|\partial_x^\alpha \boldsymbol u\|$ with $0\le|\alpha|\le [\frac{n}{2}]$, but
$\lambda$ is restricted in $(-\frac{n}{n+2},0)$.
Based on Theorem \ref{th-nonlinear}, applying the new developed
time-weighted energy method, we further improve the optimal decay rates of $\|(v,\boldsymbol u)\|$ in Theorem \ref{th-com-0-in} for all $\lambda\in(-1,0)$ and $n\ge2$.
But for the optimal decay rates to the derivatives of the solutions as
$\lambda \in (-1,-\frac{n}{n+2})$, they still remain open.
\end{remark}

\begin{theorem}[Optimal logarithmic decays for the critical case of $\lambda=-1$] \label{th-com-1-in}
For $n\ge7$, $\lambda=-1$ and $N\ge[\frac{n}{2}]+2$,
there exists a constant $\varepsilon_0>0$ such that
the solution $(v,\boldsymbol u)$ of the nonlinear system \eqref{eq-vbdu}
corresponding to small initial data
$\|(v_0,\boldsymbol u_0)\|_{L^1\cap H^N}\le \varepsilon_0$
exists globally and satisfies
\begin{equation} \label{eq-com-1-in}
\begin{cases}
\|v(t,\cdot)\|\lesssim |\ln(e+t)|^{-\frac{n}{4}}, \\
\|\boldsymbol u(t,\cdot)\|\lesssim (1+t)^{-1}\cdot|\ln(e+t)|^{-\frac{n}{4}-\frac{1}{2}}.
\end{cases}
\end{equation}
The above decay estimates are optimal and consistent with the linearized hyperbolic system.
\end{theorem}

\begin{remark}
For the critical $\lambda=-1$, the optimal decay of $\|v\|$ of the nonlinear Euler system
\eqref{eq-vbdu} is powers of the logarithmic function, i.e. $|\ln(e+t)|^{-\frac{n}{4}}$,
which differs from the classical algebraical decays.
To the best of our knowledge, this is the first result that
shows the optimal logarithmical decays of the damped Euler equations.
\end{remark}

All the above decay estimates are valid for the Euler equation \eqref{eq-Euler}.

\begin{corollary}
For $n\ge2$ and $\lambda\in(-\frac{n}{n+2},0)$,
there exists a constant $\varepsilon_0>0$, such that
the solution $(\rho,\boldsymbol u)$ of the nonlinear system \eqref{eq-Euler}
corresponding to initial data $(\rho_0,\boldsymbol u_0)$
with small energy $\|(\rho_0-1,\boldsymbol u_0)\|_{L^1\cap H^{[\frac{n}{2}]+3}}\le\varepsilon_0$
exists globally and satisfies
\begin{equation} \label{eq-nonlinear-Euler}
\begin{cases}
\|\partial_x^\alpha (\rho-1)\|\lesssim (1+t)^{-\frac{1+\lambda}{4}n-\frac{1+\lambda}{2}|\alpha|},
\quad &0\le |\alpha|\le [\frac{n}{2}]+1,\\
\|\partial_x^\alpha \boldsymbol u\|\lesssim
(1+t)^{-\frac{1+\lambda}{4}n-\frac{1+\lambda}{2}(|\alpha|+1)+\lambda},
\quad &0\le |\alpha|\le [\frac{n}{2}],\\
\|(v,\boldsymbol u)\|_{H^{[\frac{n}{2}]+3}}\lesssim 1.
\end{cases}
\end{equation}
The first two decay estimates in \eqref{eq-nonlinear-Euler}
(i.e., the decay estimates on $\|\partial_x^\alpha (\rho-1)\|$
with $0\le|\alpha|\le [\frac{n}{2}]+1$ and
$\|\partial_x^\alpha \boldsymbol u\|$ with $0\le|\alpha|\le [\frac{n}{2}]$) are optimal.

For $n\ge2$, $q\in[2,\infty]$, $k\ge 3+[\gamma_{2,q}]$
with $\gamma_{2,q}:=n(1/2-1/q)$, and $\lambda\in(-\frac{n}{n+2},0)$,
let $(\rho,\boldsymbol u)$ be the solution to the nonlinear system \eqref{eq-vbdu}
corresponding to initial data $(\rho_0,\boldsymbol u_0)$
with small energy such that
$\|(\rho_0-1,\boldsymbol u_0)\|_{L^1\cap H^{[\frac{n}{2}]+k}}\le\varepsilon_0$,
where $\varepsilon_0>0$ is a small constant only depending on $n,q,k$ and
the constants $\gamma,\mu,\lambda$ in the system.
Then $(\rho-1,\boldsymbol u)\in L^\infty(0,+\infty;H^{[\frac{n}{2}]+k})$ and satisfies
\begin{equation} \label{eq-nonlinear-Lp-Euler}
\begin{cases}
\|\partial_x^\alpha (\rho-1)\|_{L^q}\lesssim
(1+t)^{-\frac{1+\lambda}{2}\gamma_{1,q}-\frac{1+\lambda}{2}|\alpha|},
\quad &0\le |\alpha|\le 1,\\
\|\boldsymbol u\|_{L^q}\lesssim
(1+t)^{-\frac{1+\lambda}{2}\gamma_{1,q}-\frac{1-\lambda}{2}},
\end{cases}
\end{equation}
where $\gamma_{1,q}=n(1-1/q)$.
The decay estimates in \eqref{eq-nonlinear-Lp-Euler} are optimal.
\end{corollary}

\begin{corollary}
For $n\ge2$, $N\ge[\frac{n}{2}]+2$ and $\lambda\in(-1,0)$,
there exists a constant $\varepsilon_0>0$ such that
the solution $(\rho,\boldsymbol u)$ of the nonlinear system \eqref{eq-vbdu}
corresponding to small initial data
$\|(\rho_0-1,\boldsymbol u_0)\|_{L^1\cap H^N}\le \varepsilon_0$
exists globally and satisfies
\begin{equation*}
\begin{cases}
\|\rho(t,x)-1\|\lesssim (1+t)^{-\frac{1+\lambda}{4}n}, \\
\|\boldsymbol u(t,x)\|\lesssim (1+t)^{-\frac{1+\lambda}{4}n-\frac{1-\lambda}{2}}.
\end{cases}
\end{equation*}
The above decay estimates are optimal.

For $n\ge7$, $N\ge[\frac{n}{2}]+2$ and $\lambda=-1$,
there exists a constant $\varepsilon_0>0$ such that
the solution $(\rho,\boldsymbol u)$ of the nonlinear system \eqref{eq-vbdu}
corresponding to small initial data
$\|(\rho_0-1,\boldsymbol u_0)\|_{L^1\cap H^N}\le \varepsilon_0$
exists globally and satisfies
\begin{equation*}
\begin{cases}
\|\rho(t,x)-1\|\lesssim |\ln(e+t)|^{-\frac{n}{4}}, \\
\|\boldsymbol u(t,x)\|\lesssim (1+t)^{-1}\cdot|\ln(e+t)|^{-\frac{n}{4}-\frac{1}{2}}.
\end{cases}
\end{equation*}
The above decay estimates are optimal.
\end{corollary}

The paper is organized as follows.
We first leave the optimal decay estimates of the time-dependent damped
wave equations and the linearized system \eqref{eq-vu} into Appendix.
In Section 2 we formulate the optimal decay rates of the solutions with high-order derivatives
up to $[\frac{n}{2}]$-th order
for the nonlinear system \eqref{eq-vbdu} with  $\lambda\in(-\frac{n}{n+2},0)$.
In Section 3, by developing a new approach combined the Green function method with the time-weighted energy method,  we further improve the optimal decay rates of $\|(v,\boldsymbol u)\|$ for all $\lambda\in(-1,0)$.
Finally, the critical case of $\lambda=-1$
with optimal logarithmic decays is considered in Section 4.

\section{Green function method}

In this section we apply the technical Fourier analysis and the Green function method
to the study of the asymptotic behavior of nonlinear system \eqref{eq-vbdu}.
We rewrite \eqref{eq-vbdu} as
\begin{equation} \label{eq-non-sys}
\partial_t
\begin{pmatrix}
v \\
\boldsymbol u
\end{pmatrix}
=
\begin{pmatrix}
0 & -\nabla\cdot \\
-\nabla & -\frac{\mu}{(1+t)^\lambda}
\end{pmatrix}
\begin{pmatrix}
v \\
\boldsymbol u
\end{pmatrix}
+
\begin{pmatrix}
-\boldsymbol u\cdot\nabla v-\varpi v\nabla\cdot \boldsymbol u \\
-(\boldsymbol u\cdot\nabla) \boldsymbol u-\varpi v\nabla v
\end{pmatrix},
\end{equation}
and the solution can be expressed,  by the Duhamel principle, as follows
\begin{equation} \label{eq-Duhamel}
\begin{pmatrix}
v(t,x) \\
\boldsymbol u(t,x)
\end{pmatrix}
=\mathcal{G}(t,0)
\begin{pmatrix}
v(0,x) \\
\boldsymbol u(0,x)
\end{pmatrix}
+
\int_0^t\mathcal{G}(t,s)Q(s,x)ds,
\end{equation}
where
$$
Q(s,x)=
\begin{pmatrix}
Q_1(s,x) \\
Q_2(s,x)
\end{pmatrix}
: =
\begin{pmatrix}
-\boldsymbol u\cdot\nabla v-\varpi v\nabla\cdot \boldsymbol u \\
-(\boldsymbol u\cdot\nabla) \boldsymbol u-\varpi v\nabla v
\end{pmatrix},
\quad
\mathcal{G}(t,s)
=
\begin{pmatrix}
\mathcal{G}_{11}(t,s) & \mathcal{G}_{12}(t,s) \\
\mathcal{G}_{21}(t,s) & \mathcal{G}_{22}(t,s)
\end{pmatrix}.
$$
The Green matrix $\mathcal{G}(t,s)$ represents the evolution
of the linear system starting from time $s$ to $t$.
It should be noted that $\mathcal{G}(t,s)\ne \mathcal{G}(t-s,0)$
since the time-asymptotically growing damping $\frac{\mu}{(1+t)^\lambda}$ on $(s,t)$
is completely different from the damping on $(0,t-s)$.
Moreover, there is no explicit (matrix exponential type) expression
of the Green matrix $\mathcal{G}(t,s)$
due to the time-dependent coefficient $b(t)$.
In fact, the abstract expression of $\mathcal{G}(t,s)$ based on the Peano-Baker formula
(see Proposition A.3 in \cite{Wirth-JDE06} for example) is
$$
\mathcal{G}(t,s)
=I+\sum_{k=1}^{\infty}\int_s^t\mathcal{A}(t_1,\xi)\int_s^{t_1}
\mathcal{A}(t_2,\xi)\cdots\int_s^{t_{k-1}}
\mathcal{A}(t_k,\xi)dt_k\cdots dt_2dt_1,
$$
with the non-commutative
($\mathcal{A}(t,\xi)\mathcal{A}(s,\xi)\ne \mathcal{A}(s,\xi)\mathcal{A}(t,\xi)$
for general $s\ne t$) matrix
$$\mathcal{A}(t,\xi):=
\begin{pmatrix}
0 & -i\xi^\mathrm{T} \\
-i\xi & -\frac{\mu}{(1+t)^\lambda}I_{n\times n}
\end{pmatrix},$$
where $(\cdot)^\mathrm{T}$ is the transpose of a vector.
The exact time decay estimates of $\mathcal{G}(t,s)$
are shown in Theorem \ref{th-linear} in Appendix,
where we write the Green function of time and space $\mathcal{G}(t,s;x,\xi)$
as $\mathcal{G}(t,s)$ for the sake of simplicity. Here and hereafter, in order to emphasize the effect of time $t$ for a given function $v(t,x)$,
we often  simply write $v(t)$ instead of $v(t,x)$ if there is no confusion.

The linearized system of \eqref{eq-vbdu} (or \eqref{eq-non-sys}) is
\begin{equation} \label{eq-vu-linear}
\begin{cases}
\displaystyle
\partial_t v+\nabla\cdot\boldsymbol u=0,\\
\displaystyle
\partial_t \boldsymbol u+\nabla v+\frac{\mu}{(1+t)^\lambda}\boldsymbol u=0,\\
\displaystyle
v|_{t=0}=v_0(x), \quad \boldsymbol u|_{t=0}=\boldsymbol u_0(x).
\end{cases}
\end{equation}
Let $u:=\Lambda^{-1}\nabla\cdot\boldsymbol u$ and
$\boldsymbol w:=\Lambda^{-1}\mathrm{curl}\,\boldsymbol u$
(with $(\mathrm{curl}\,\boldsymbol u)_j^k:=\partial_{x_j}u^k-\partial_{x_k}u^j$
for $\boldsymbol u=(u^1,\dots,u^n)$), see \cite{TanZ-JDE12} for example,
where the pseudo differential operator $\Lambda$ is defined by
$\Lambda^s v:=\mathscr{F}^{-1}(|\xi|^s\hat v(\xi))$ for $s\in\mathbb R$.
Then the linearized system \eqref{eq-vu-linear} is equivalent to
\begin{equation} \label{eq-vu}
\begin{cases}
\displaystyle
\partial_t v+\Lambda u=0,\\
\displaystyle
\partial_t u-\Lambda v+\frac{\mu}{(1+t)^\lambda} u=0,\\
\displaystyle
\partial_t \boldsymbol w+\frac{\mu}{(1+t)^\lambda} \boldsymbol w=0,\\
\displaystyle
v|_{t=0}=v_0(x), \quad u|_{t=0}=u_0(x), \quad \boldsymbol w|_{t=0}=\boldsymbol w_0(x),
\end{cases}
\end{equation}
where $u_0(x)=\Lambda^{-1}\nabla\cdot\boldsymbol u_0(x)$
and $\boldsymbol w_0(x)=\Lambda^{-1}\mathrm{curl}\,\boldsymbol u_0(x)$.
We note that the estimates on $(v,\boldsymbol u)$ are equivalent to
the estimates on $(v,u,\boldsymbol w)$.
From the equation \eqref{eq-vu}$_3$, the vorticity $\boldsymbol w(t,x)$
of the linearized system decays to zero super-exponentially
(as $\boldsymbol w_0(x)e^{-\mu(1+t)^{1-\lambda}/(1-\lambda)}$
with $\lambda\in[-1,0)$),
which is faster than any algebraical decays.
So we only consider the first two equations of \eqref{eq-vu}.

In order to formulate the optimal decay rates of the linearized system \eqref{eq-vu},
we consider the following two kinds of wave equations with time-dependent damping
\begin{equation} \label{eq-Pv}
\begin{cases}
\displaystyle
\partial_t^2 v-\Delta v+\frac{\mu}{(1+t)^\lambda} \partial_t v=0,
\quad x\in\mathbb R^n,
\\
\displaystyle
v|_{t=0}=v_1(x), \quad \partial_tv|_{t=0}=v_2(x),
\end{cases}
\end{equation}
and
\begin{equation} \label{eq-Pu}
\begin{cases}
\displaystyle
\partial_t^2 u-\Delta u+\partial_t\Big(\frac{\mu}{(1+t)^\lambda}u\Big)=0,
\quad x\in\mathbb R^n,
\\
\displaystyle
u|_{t=0}=u_1(x), \quad \partial_tu|_{t=0}=u_2(x),
\end{cases}
\end{equation}
which are satisfied by the solutions $v(t,x)$ and $u(t,x)$ of \eqref{eq-vu} respectively.

We show that the optimal decay rate of $u(t,x)$ in the damped wave equation \eqref{eq-Pu}
is faster than the optimal decay rate of $v(t,x)$ in the wave equation \eqref{eq-Pv},
and further we prove that $u(t,x)$ in the damped linear system \eqref{eq-vu}
decays optimally faster than all the damped wave equations \eqref{eq-Pv} and \eqref{eq-Pu}.
Therefore, there are cancellations between the evolution of initial data
if we regard $u(t,x)$ in the linear system \eqref{eq-vu} as a solution of
the wave equation \eqref{eq-Pu}
with initial data $u_1(x)=u_0(x)$ and $u_2(x)=\Lambda v_0(x)-\mu u_0(x)$.

The optimal decay estimates of the time-dependent damped linearized system \eqref{eq-vu},
together with the optimal decays of the wave equations \eqref{eq-Pv} and \eqref{eq-Pu},
are proved in Appendix (Theorem \ref{th-wave} and Theorem \ref{th-linear})
by means of the technical Fourier analysis.

Compared with the under-damping case $\lambda\in[0,1)$ in \cite{Ji-Mei-1},
here the over-damping case $\lambda\in[-1,0)$ gives rise to two main difficulties
in the decay estimates of the nonlinear system:

(i) $\|\partial_x^\alpha v\|$ decays slowly
since $\|\mathcal{G}_{11}(t,0)v_0\|\approx (1+t)^{-\frac{1+\lambda}{4}n}$ for $\lambda\in(-1,0)$
and $\|\mathcal{G}_{11}(t,0)v_0\|\approx |\ln(e+t)|^{-\frac{n}{4}}$ for $\lambda=-1$.
One should be careful in calculating the estimates on $\int_0^t\mathcal{G}(t,s)Q(s)ds$.

(ii) The high-order energy estimates on $\|\partial_x^{[n/2]+3}(v,\boldsymbol u)\|$
are deduced through energy method,
but the estimate on $\|\partial_x^{[n/2]+3}v\|$
needs the estimate
$$\int b(t)\partial_x^{[n/2]+2}\boldsymbol u(t)\cdot \nabla\partial_x^{[n/2]+2}v(t),$$
where the over-damping coefficient
$b(t)=\frac{\mu}{(1+t)^\lambda}=\mu(1+t)^{|\lambda|}$ for $\lambda\in[-1,0)$
is growing and causes trouble for $\lambda$ near $-1$.

\subsection{High-order energy estimates with over-damping}

For the closure of the decay estimates of nonlinear system \eqref{eq-vbdu},
we need to formulate high-order energy estimates.
Note that the over-damping coefficient $b(t)=\frac{\mu}{(1+t)^\lambda}$
is growing for $\lambda\in[-1,0)$.

\begin{lemma} \label{le-energy}
Let $(v_0,\boldsymbol u_0)\in H^{[\frac{n}{2}]+k}$ with $k\ge2$, and
$(v, \boldsymbol u)(x,t)$ be the solutions of the nonlinear system \eqref{eq-vbdu} for $t\in [0,T]$ with a positive number $T$, and satisfy
\begin{equation} \label{eq-apriori-energy}
\|(v(t),\boldsymbol u(t))\|_{H^{[\frac{n}{2}]+2}}\le \delta_0\frac{1}{b(t)},
\end{equation}
where $\delta_0>0$ is a small number.
Then it holds
\begin{equation} \label{eq-energy-high}
\|(v,\boldsymbol u)\|_{H^{[\frac{n}{2}]+k}}^2
+\int_0^t \Big(\frac{1}{b(s)}\|\nabla v(s)\|_{H^{[\frac{n}{2}]+k-1}}^2
+b(s)\|\boldsymbol u(s)\|_{H^{[\frac{n}{2}]+k}}^2\Big)ds
\lesssim \|(v_0,\boldsymbol u_0)\|_{H^{[\frac{n}{2}]+k}}^2, \ \ \ t\in [0,T].
\end{equation}
\end{lemma}
{\it\bfseries Proof.}
The case of time-independent damping (i.e. $\lambda=0$) is proved in \cite{TanZ-JDE13},
and the under-damping case $\lambda\in(0,1)$ is proved in \cite{Ji-Mei-1}. But, different from the previous studies,
for the over-damping case with $\lambda \in[-1,0)$, here the main difficulty lies in the absence of uniform upper bound
of the over-damping coefficient. We divide the proof into four steps.

{\it Step 1}.   For $0\le j\le [\frac{n}{2}]+k-1$, we have
\begin{equation} \label{eq-zenergy1}
\frac{d}{dt}\|\partial_x^j(v,\boldsymbol u)\|^2
+b(t)\|\partial_x^j \boldsymbol u\|^2
\lesssim \|(v,\boldsymbol u)\|_{H^{[\frac{n}{2}]+2}}
\cdot (\|\partial_x^{j+1}v\|^2+\|\partial_x^j\boldsymbol u\|^2).
\end{equation}
This can be proved by applying $\partial_x^j$ to \eqref{eq-vbdu}
and then multiplying the resultant equations by $\partial_x^j(v,\boldsymbol u)$,
summing them up and integrating it with respect to $x$ over $\mathbb R^n$.
Here we omit the details.

{\it Step 2}.   By applying $\partial_x^{j+1}$ to \eqref{eq-vbdu} with $0\le j\le [\frac{n}{2}]+k-1$,
and multiplying the resultant equations by $\partial_x^{j+1}(v,\boldsymbol u)$, and
summing them up and integrating it over $\mathbb R^n$,  we have
\begin{equation} \label{eq-zenergy2}
\frac{d}{dt}\|\partial_x^{j+1}(v,\boldsymbol u)\|^2
+b(t)\|\partial_x^{j+1} \boldsymbol u\|^2
\lesssim \|(v,\boldsymbol u)\|_{H^{[\frac{n}{2}]+2}}
\cdot (\|\partial_x^{j+1}v\|^2+\|\partial_x^{j+1}\boldsymbol u\|^2).
\end{equation}

{\it Step 3}.   For $0\le j\le [\frac{n}{2}]+k-1$, we can obtain
\begin{equation} \label{eq-zenergy3}
\frac{d}{dt}\int\partial_x^j\boldsymbol u\cdot \nabla\partial_x^{j}v
+\|\partial_x^{j+1}v\|^2
\lesssim
b(t)\|\partial_x^{j}\boldsymbol u\|\cdot\|\partial_x^{j+1}v\|
+\|\partial_x^{j+1}\boldsymbol u\|^2
+\|(v,\boldsymbol u)\|_{H^{[\frac{n}{2}]+2}}
\cdot (\|\partial_x^{j+1}v\|^2+\|\partial_x^{j+1}\boldsymbol u\|^2).
\end{equation}
In fact, this can be proved by applying $\partial_x^j$ to \eqref{eq-vbdu}$_2$
and multiplying it by $\partial_x^{j+1}v$ (specifically, $\nabla \partial_x^j v$),
utilizing \eqref{eq-vbdu}$_1$ to dealing with the mixed space-time derivative term
$\int \partial_x^j\partial_t\boldsymbol u\cdot \partial_x^{j+1}v$,
that is,
\begin{equation*}
\|\partial_x^j\nabla v\|^2
+\int \partial_t(\partial_x^j \boldsymbol u)\cdot \nabla \partial_x^j v
+\int b(t)\partial_x^j \boldsymbol u\cdot \nabla \partial_x^j v
=\int\partial_x^j Q_2\cdot \nabla \partial_x^j v,
\end{equation*}
and
\begin{align*}
\int \partial_t(\partial_x^j \boldsymbol u)\cdot \nabla \partial_x^j v
=&\frac{d}{dt}\int \partial_x^j \boldsymbol u\cdot \nabla \partial_x^j v
+\int \partial_x^j (\nabla \cdot \boldsymbol u) \cdot \partial_t\partial_x^j v
\\
=&\frac{d}{dt}\int \partial_x^j \boldsymbol u\cdot \nabla \partial_x^j v
-\| \partial_x^j (\nabla \cdot \boldsymbol u)\|^2
+\int \partial_x^j (\nabla \cdot \boldsymbol u) \cdot \partial_x^j Q_1.
\end{align*}
Applying Cauchy's inequality to \eqref{eq-zenergy3}, we then arrive at
\begin{equation} \label{eq-zenergy3a}
\frac{d}{dt}\int\partial_x^j\boldsymbol u\cdot \nabla\partial_x^{j}v
+\|\partial_x^{j+1}v\|^2
\lesssim
b^2(t)\|\partial_x^{j}\boldsymbol u\|^2
+\|\partial_x^{j+1}\boldsymbol u\|^2
+\|(v,\boldsymbol u)\|_{H^{[\frac{n}{2}]+2}}
\cdot (\|\partial_x^{j+1}v\|^2+\|\partial_x^{j+1}\boldsymbol u\|^2).
\end{equation}
Next,  we multiply \eqref{eq-zenergy3a} by $\frac{1}{b(t)}$,
for $0\le j\le [\frac{n}{2}]+k-1$, to have
\begin{align*}
&\frac{d}{dt}\Big(\frac{1}{b(t)}\int\partial_x^j\boldsymbol u\cdot \nabla\partial_x^{j}v\Big)
+\frac{1}{b(t)}\|\partial_x^{j+1}v\|^2
\\
&\lesssim
\frac{|b'(t)|}{b^2(t)}\int\big|\partial_x^j\boldsymbol u\cdot \nabla\partial_x^{j}v\big|
+b(t)\|\partial_x^{j}\boldsymbol u\|^2
+\frac{1}{b(t)}\|\partial_x^{j+1}\boldsymbol u\|^2
+\frac{1}{b(t)}\|(v,\boldsymbol u)\|_{H^{[\frac{n}{2}]+2}}
\cdot (\|\partial_x^{j+1}v\|^2+\|\partial_x^{j+1}\boldsymbol u\|^2)
\\
&\lesssim
\varepsilon_1\frac{1}{b(t)}\|\partial_x^{j+1}v\|^2
+b(t)\|\partial_x^{j}\boldsymbol u\|^2
+\frac{1}{b(t)}\|\partial_x^{j+1}\boldsymbol u\|^2
+\frac{1}{b(t)}\|(v,\boldsymbol u)\|_{H^{[\frac{n}{2}]+2}}
\cdot (\|\partial_x^{j+1}v\|^2+\|\partial_x^{j+1}\boldsymbol u\|^2),
\end{align*}
where $\varepsilon_1>0$ is a small number.
Therefore, for $0\le j\le [\frac{n}{2}]+k-1$, we have
\begin{align} \nonumber
&\frac{d}{dt}\Big(\frac{1}{b(t)}\int\partial_x^j\boldsymbol u\cdot \nabla\partial_x^{j}v\Big)
+\frac{1}{b(t)}\|\partial_x^{j+1}v\|^2
\\ \label{eq-zenergy4}
\lesssim&
b(t)\|\partial_x^{j}\boldsymbol u\|^2
+\frac{1}{b(t)}\|\partial_x^{j+1}\boldsymbol u\|^2
+\frac{1}{b(t)}\|(v,\boldsymbol u)\|_{H^{[\frac{n}{2}]+2}}
\cdot (\|\partial_x^{j+1}v\|^2+\|\partial_x^{j+1}\boldsymbol u\|^2).
\end{align}

{\it Step 4}. Multiplying \eqref{eq-zenergy4} by a small number $\varepsilon_2>0$,
summing it up with \eqref{eq-zenergy1} and \eqref{eq-zenergy2},
we have
$$
\frac{d}{dt}\|(v,\boldsymbol u)\|_{H^{[\frac{n}{2}]+k}}^2
+\frac{d}{dt}\Big(\varepsilon_2\sum_{j=0}^{[{n}/{2}]+k-1}
\frac{1}{b(t)}\int\partial_x^j\boldsymbol u\cdot \nabla\partial_x^{j}v\Big)
+\frac{1}{b(t)}\|\nabla v\|_{H^{[\frac{n}{2}]+k-1}}^2
+b(t)\|\boldsymbol u\|_{H^{[\frac{n}{2}]+k}}^2\le0,
$$
provided with the a priori assumption \eqref{eq-apriori-energy}. Let us choose
 $\varepsilon_2>0$ to be small such that
$$
\Big|\varepsilon_2\sum_{j=0}^{[{n}/{2}]+k-1}
\frac{1}{b(t)}\int\partial_x^j\boldsymbol u\cdot \nabla\partial_x^{j}v\Big|
\le \frac{1}{2}\|(v,\boldsymbol u)\|_{H^{[\frac{n}{2}]+k}}^2,
$$
then we obtain \eqref{eq-energy-high}.
The proof is completed.
$\hfill\Box$

\vskip2mm
The most tricky part lies in the treatment of
$b(t)\|\partial_x^{j}\boldsymbol u\|\cdot\|\partial_x^{j+1}v\|$
in \eqref{eq-zenergy3}, where $\|\partial_x^{j+1}v\|^2$ is the only good term,
therefore $b^2(t)\|\partial_x^{j}\boldsymbol u\|^2$ arises
(if Cauchy's inequality is applied)
and grows faster than $b(t)\|\partial_x^{j}\boldsymbol u\|^2$
in \eqref{eq-zenergy1}.
This is the reason of the a priori assumption \eqref{eq-apriori-energy}.
We can prove that \eqref{eq-apriori-energy} is satisfied for $\lambda$ near zero.
However, the decay estimates required in \eqref{eq-apriori-energy} are
not true for $\lambda\in[-1,0)$ near $-1$, especially for the case $\lambda=-1$.
In fact,
$$
\|\partial_x^\alpha\mathcal{G}_{11}(t,0)v_0\|
\approx |\ln(e+t)|^{-\frac{n}{4}-\frac{|\alpha|}{2}}, \text{~for~} \lambda=-1,
$$
and the decay condition
$\|v(t)\|_{H^{[n/2]+2}}\le \delta_0(1+t)^{-1}$ in
\eqref{eq-apriori-energy} is not valid.

We can relax the decay condition of high-order estimates in \eqref{eq-apriori-energy}
to a wider range of $\lambda$.
The crucial point is to avoid the decay conditions of
$\|v(t)\|_{H^{[n/2]+2}}$.
For application, we prove the following inequality
which can be regarded as a generalized Gr\"onwall's inequality with relaxation.

\begin{lemma}[Gr\"onwall's inequality with relaxation] \label{le-Gronwall}
Assume that $\omega(t)$, $g(t)$, and $H(t)$ are nonnegative functions,
$C_2\ge C_1>0$, $\theta\in(0,1)$, $\eta>0$, all are constants,
and
$F(t)$ satisfies (note that $F(t)$ is not necessarily nonnegative)
\begin{equation}\label{new-1}
C_1H(t)-g(t)\le F(t) \le C_2H(t)+g(t),
\end{equation}
and the following differential inequality
\begin{equation} \label{eq-Gronwall}
\frac{d}{dt}F(t)+\eta F(t)\le \omega(t)H^\theta(t)+g(t), \quad \forall t>0,
\end{equation}
then
\begin{equation}\label{new-2}
F(t)\lesssim \max\{F(0),\sup_{s\in(0,t)}
((\omega(s)/\eta)^\frac{1}{1-\theta}+g(s)(1+1/\eta))\}.
\end{equation}
and
\begin{equation}\label{new-3}
H(t)\lesssim \max\{F(0),\sup_{s\in(0,t)}
((\omega(s)/\eta)^\frac{1}{1-\theta}+g(s)(1+1/\eta))\}.
\end{equation}
Furthermore, if  $\omega(t)$ and $g(t)$
are monotonically decreasing, then
\begin{equation} \label{new-4}
F(t)
\lesssim
F(0)e^{-\frac{\eta}{2}t}
+\Big(\frac{1}{\eta^{1/(1-\theta)}}\omega^\frac{1}{1-\theta}(0)
+\big(1+\frac{1}{\eta}\big)g(0)\Big)e^{-\frac{\eta}{8}t}
+\frac{1}{\eta^{1/(1-\theta)}}\omega^\frac{1}{1-\theta}\big(\frac{t}{2}\big)
+\big(1+\frac{1}{\eta}\big)g\big(\frac{t}{2}\big).
\end{equation}
and
\begin{equation} \label{eq-G-decay}
H(t)
\lesssim
F(0)e^{-\frac{\eta}{2}t}
+\Big(\frac{1}{\eta^{1/(1-\theta)}}\omega^\frac{1}{1-\theta}(0)
+\big(1+\frac{1}{\eta}\big)g(0)\Big)e^{-\frac{\eta}{8}t}
+\frac{1}{\eta^{1/(1-\theta)}}\omega^\frac{1}{1-\theta}\big(\frac{t}{2}\big)
+\big(1+\frac{1}{\eta}\big)g\big(\frac{t}{2}\big).
\end{equation}
\end{lemma}
{\it\bfseries Proof.}
We may assume that $C_1=\frac{1}{2}$, $C_2=2$, and $\theta=\frac{1}{2}$.
Other situation follows similarly.
For any $t>0$, if $F(t)>F(0)$, then two cases happen:
(i) $F(t)=\sup_{s\in(0,t)}F(s)$, such that $F'(t)\ge0$;
(ii) there exists a number $s\in(0,t)$, such that $F'(s)=0$ and $F(s)> F(t)$.
In both cases, we can find a number
$s\in(0,t]$, such that $F'(s)\ge0$ and $F(s)\ge F(t)$.
Therefore, according to the differential inequality \eqref{eq-Gronwall},
we have
\[
\eta F(s)\le \omega(s)H^\frac{1}{2}(s)+g(s)
\le \frac{1}{\eta}\omega^2(s)+\frac{\eta}{4}H(s)+g(s)
\le \frac{1}{\eta}\omega^2(s)+\frac{1}{2}\eta F(s)+(\frac{\eta}{2}+1)g(t),
\]
which implies
\[
F(t)\le F(s)\lesssim \omega^2(s)/\eta^2+g(s)(1+1/\eta).
\]
This immediately guarantees \eqref{new-2}. On the other hand, \eqref{new-1} implies
\[
H(t)\lesssim F(t)+g(t).
\]
This together with \eqref{new-2} proves \eqref{new-3}.

If $\omega(t)$ and $g(t)$ are monotonically decreasing, then
according to \eqref{eq-Gronwall} and Young's inequality
\begin{align*}
\frac{d}{dt}F(t)+\eta F(t)
&\le \omega(t)H^\theta(t)+g(t)
\\
&\le \frac{1}{2(C_1\eta)^{\theta/(1-\theta)}}\omega^{1/(1-\theta)}(t)
+\frac{C_1}{2}\eta H(t)+g(t)
\\
&\le \frac{1}{2(C_1\eta)^{\theta/(1-\theta)}}\omega^{1/(1-\theta)}(t)
+\frac{1}{2}\eta F(t)+(\frac{\eta}{2}+1)g(t),
\end{align*}
we have
\begin{align*}
\frac{d}{dt}(e^{\frac{\eta}{2}t}F(t))
&=e^{\frac{\eta}{2}t}\Big(\frac{d}{dt}F(t)+\frac{\eta}{2}F(t)\Big)
\\
&\le
e^{\frac{\eta}{2}t}\Big(\frac{1}{2(C_1\eta)^{\theta/(1-\theta)}}\omega^{1/(1-\theta)}(t)
+(\frac{\eta}{2}+1)g(t)\Big)
\\
&\lesssim e^{\frac{\eta}{2}t}(\omega^{1/(1-\theta)}(t)+g(t)),
\end{align*}
where we have slightly abused the notion ``$\lesssim$''
such that the inequality depends on $\eta$ and $\theta$
and the dependence is clear.
Integrating it with respect to $t$ over $(0,t)$ gives
\begin{align*}
F(t)&\lesssim F(0)e^{-\frac{\eta}{2}t}
+\int_0^t e^{-\frac{\eta}{2}(t-s)}\Big(\omega^{1/(1-\theta)}(s)+g(s)\Big)ds
\\
&\lesssim
F(0)e^{-\frac{\eta}{2}t}
+\int_0^\frac{t}{2} e^{-\frac{\eta}{4}t}\Big(\omega^{1/(1-\theta)}(0)+g(0)\Big)ds
+\int_\frac{t}{2}^t e^{-\frac{\eta}{2}(t-s)}\Big((\omega^{1/(1-\theta)}(\frac{t}{2})+g(\frac{t}{2})\Big)ds
\\
&\lesssim
F(0)e^{-\frac{\eta}{2}t}
+\Big(\omega^{1/(1-\theta)}(0)+g(0)\Big)e^{-\frac{\eta}{8}t}
+\omega^{1/(1-\theta)}(\frac{t}{2})+g(\frac{t}{2}),
\end{align*}
since $te^{-\frac{\eta}{4}t}\lesssim e^{-\frac{\eta}{8}t}$.
Thus, \eqref{new-4} and \eqref{eq-G-decay} are immediately obtained.
The proof is completed.
$\hfill\Box$

We modify the high-order estimates Lemma \ref{le-energy}
such that $\|v(t)\|_{H^{[\frac{n}{2}]+2}}$
does not necessarily decay as fast as $\frac{1}{b(t)}$.
The key ingredient is to avoid the estimate on $\|v(t)\|$
such that Step 1 in the proof of Lemma \ref{le-energy} is excluded.

\begin{lemma} \label{le-energy-G}
Let $(v_0,\boldsymbol u_0)\in H^{[\frac{n}{2}]+k}$ with $k\ge2$, and let
$(v, \boldsymbol u)(x,t)$ be the solutions of the nonlinear system \eqref{eq-vbdu} for $t\in [0,T]$ with a positive number $T$, and satisfy
\begin{equation} \label{eq-apriori-energy-G}
\|(v(t),\boldsymbol u(t))\|_{H^{[\frac{n}{2}]+2}}\le \delta_0,
\quad
\|\boldsymbol u(t)\|_{H^{[\frac{n}{2}]+k-1}}\le \delta_0\frac{\omega(t)}{b(t)}, \ \ t\in [0,T],
\end{equation}
where $\delta_0>0$ is a small number
and $\omega(t)$ is a nonnegative decreasing function.
Then it holds
\begin{equation} \label{eq-energy-high-G}
\|\nabla(v,\boldsymbol u)\|_{H^{[\frac{n}{2}]+k-1}}^2
\lesssim \|\nabla(v_0,\boldsymbol u_0)\|_{H^{[\frac{n}{2}]+k-1}}^2
+\delta_0^2\cdot\omega^2(t/2), \ \ \ t\in [0,T].
\end{equation}
\end{lemma}
{\it\bfseries Proof.}
According to the estimates \eqref{eq-zenergy2} and \eqref{eq-zenergy3}
in Step 2 and Step 3 of the proof of Lemma \ref{le-energy},
for $0\le j\le [\frac{n}{2}]+k-1$,
we have
\begin{equation} \label{eq-zenergy-G1}
\frac{d}{dt}\|\partial_x^{j+1}(v,\boldsymbol u)\|^2
+b(t)\|\partial_x^{j+1} \boldsymbol u\|^2
\lesssim \|(v,\boldsymbol u)\|_{H^{[\frac{n}{2}]+2}}
\cdot \|\partial_x^{j+1}v\|^2.
\end{equation}
and
\begin{equation} \label{eq-zenergy-G2}
\frac{d}{dt}\int\partial_x^j\boldsymbol u\cdot \nabla\partial_x^{j}v
+\|\partial_x^{j+1}v\|^2
\lesssim
b(t)\|\partial_x^{j}\boldsymbol u\|\cdot\|\partial_x^{j+1}v\|
+\|\partial_x^{j+1}\boldsymbol u\|^2,
\end{equation}
where we have used the a priori assumption \eqref{eq-apriori-energy-G}
such that $\|(v(t),\boldsymbol u(t))\|_{H^{[\frac{n}{2}]+2}}\le \delta_0$
with a small $\delta_0$.
Multiplying \eqref{eq-zenergy-G2} by a small number $\tilde\varepsilon_1$ (only depending on the dimension $n$) such that
$$
{\tilde\varepsilon_1}  \int\big|\partial_x^j\boldsymbol u\cdot \nabla\partial_x^{j}v\big|
\le \frac{1}{2}\|\partial_x^{j+1}v\|^2+\|\partial_x^{j}\boldsymbol u\|^2,
$$
and  making addition of  $\tilde\varepsilon_1\cdot$\eqref{eq-zenergy-G2}+\eqref{eq-zenergy-G1},
then we have
\begin{align} \nonumber
&\frac{d}{dt}\Big(\|\partial_x^{j+1}(v,\boldsymbol u)\|^2+
{\tilde\varepsilon_1} \int\partial_x^j\boldsymbol u\cdot \nabla\partial_x^{j}v\Big)
+b(t)\|\partial_x^{j+1} \boldsymbol u\|^2
+{\tilde\varepsilon_1} \|\partial_x^{j+1}v\|^2
+\frac{{\tilde\varepsilon_1} ^2}{2} \int\partial_x^j\boldsymbol u\cdot \nabla\partial_x^{j}v
\\ \nonumber
\lesssim &
{\tilde\varepsilon_1}  b(t)\|\partial_x^{j}\boldsymbol u\|\cdot\|\partial_x^{j+1}v\|
+\frac{{\tilde\varepsilon_1} ^2}{2} \int\partial_x^j\boldsymbol u\cdot \nabla\partial_x^{j}v
+\delta_0 \|\partial_x^{j+1}v\|^2
+{\tilde\varepsilon_1} \|\partial_x^{j+1}\boldsymbol u\|^2
\\ \label{eq-zenergy-G3}
\lesssim &
{\tilde\varepsilon_1}  b(t)\|\partial_x^{j}\boldsymbol u\|\cdot\|\partial_x^{j+1}v\|
+{\tilde\varepsilon_2}{\tilde\varepsilon_1} \|\partial_x^{j+1}v\|^2
+\frac{{\tilde\varepsilon_1} ^3}{{\tilde\varepsilon_2} }\|\partial_x^{j} \boldsymbol u\|^2
+\delta_0 \|\partial_x^{j+1}v\|^2
+{\tilde\varepsilon_1} \|\partial_x^{j+1}\boldsymbol u\|^2,
\end{align}
where $\tilde\varepsilon_2>0$ is another small number (only dependent on $n$)
such that ${\tilde\varepsilon_2} {\tilde\varepsilon_1} \|\partial_x^{j+1}v\|^2$
is dominated by $\frac{1}{4}{\tilde\varepsilon_1} \|\partial_x^{j+1}v\|^2$.
Noticing that $b(t)$ is growing, ${\tilde\varepsilon_1} $ can be chosen small enough,
and $\delta_0$ is small, too, we rewrite \eqref{eq-zenergy-G3} into
\begin{align} \nonumber
&\frac{d}{dt}\Big(\|\partial_x^{j+1}(v,\boldsymbol u)\|^2+
{\tilde\varepsilon_1}  \int\partial_x^j\boldsymbol u\cdot \nabla\partial_x^{j}v\Big)
+\frac{{\tilde\varepsilon_1} }{2}\|\partial_x^{j+1} \boldsymbol u\|^2
+\frac{{\tilde\varepsilon_1} }{2}\|\partial_x^{j+1}v\|^2
+\frac{{\tilde\varepsilon_1} ^2}{2} \int\partial_x^j\boldsymbol u\cdot \nabla\partial_x^{j}v
\\ \nonumber
&\lesssim
{\tilde\varepsilon_1}  b(t)\|\partial_x^{j}\boldsymbol u\|\cdot\|\partial_x^{j+1}v\|
+\frac{{\tilde\varepsilon_1} ^3}{\varepsilon_2}\|\partial_x^{j} \boldsymbol u\|^2
\\ \label{eq-zenergy-G4}
&\lesssim b(t)\|\partial_x^{j}\boldsymbol u\|\cdot\|\partial_x^{j+1}(v,\boldsymbol u)\|
+\|\partial_x^{j} \boldsymbol u\|^2.
\end{align}

Let
$$
F(t):=\|\partial_x^{j+1}(v,\boldsymbol u)\|^2+
{\tilde\varepsilon_1}  \int\partial_x^j\boldsymbol u\cdot \nabla\partial_x^{j}v,
\quad
H(t):=\|\partial_x^{j+1}(v,\boldsymbol u)\|^2,
\quad
g(t):=\|\partial_x^{j}\boldsymbol u\|^2,
$$
then
$$
\frac{1}{2}H(t)-g(t)\le F(t) \le 2H(t)+g(t),
$$
and
$$
\frac{d}{dt}F(t)+\frac{{\tilde\varepsilon_1} }{2}F(t)\lesssim
\delta_0\omega(t)H^\frac{1}{2}(t)+g(t),
$$
with $\delta_0\omega(t)=b(t)\|\partial_x^{j}\boldsymbol u\|$ decreasing,
provided the a priori assumption \eqref{eq-apriori-energy-G}.
Applying the generalized Gr\"onwall's inequality with relaxation in
Lemma \ref{le-Gronwall},
we have
$$
H(t)=\|\partial_x^{j+1}(v,\boldsymbol u)\|^2
\lesssim
F(0)e^{-\frac{\varepsilon_1}{4}t}
+\Big(\delta_0^2\omega^2(0)+g(0)\Big)e^{-\frac{\varepsilon_1}{16}t}
+\delta_0^2\omega^2(t/2)
+g(t/2).
$$
The proof is completed.
$\hfill\Box$

\subsection{Optimal $L^2$ decay estimates}

We start with the optimal $L^1$-$L^2$ decay estimates of the nonlinear
system \eqref{eq-vbdu} for the over-damping case of $\lambda\in[-1,0)$.

\begin{lemma} \label{le-decay-G}
For $\lambda\in[-1,0)$ and $t\ge s\ge T_0$
($T_0\ge0$ is a universal constant only depending on the constants $\lambda$ and $\mu$), then there hold
\begin{align} \nonumber
&\|\partial_x^\alpha \mathcal{G}_{11}(t,s)\phi(x)\|\lesssim
\Gamma^{\frac{n}{2}+|\alpha|}(t,s)\cdot
(\|\phi\|_{L^1}^l+\|\partial_x^{|\alpha|}\phi\|^h),
\\ \nonumber
&\|\partial_x^\alpha \mathcal{G}_{12}(t,s)\phi(x)\|\lesssim
(1+s)^\lambda\cdot
\Gamma^{\frac{n}{2}+|\alpha|+1}(t,s)\cdot
(\|\phi\|_{L^1}^l+\|\partial_x^{|\alpha|}\phi\|^h),
\\ \nonumber
&\|\partial_x^\alpha \mathcal{G}_{21}(t,s)\phi(x)\|\lesssim
(1+t)^\lambda\cdot
\Gamma^{\frac{n}{2}+|\alpha|+1}(t,s)\cdot
(\|\phi\|_{L^1}^l+\|\partial_x^{|\alpha|}\phi\|^h),
\\ \label{eq-decay-G}
&\|\partial_x^\alpha \mathcal{G}_{22}(t,s)\phi(x)\|\lesssim
\Big(\frac{1+t}{1+s}\Big)^\lambda\cdot
\Gamma^{\frac{n}{2}+|\alpha|}(t,s)\cdot
(\|\phi\|_{L^1}^l+\|\partial_x^{|\alpha|}\phi\|^h).
\end{align}
Furthermore,
\begin{align} \label{eq-decay-G-opt}
\|\partial_x^\alpha \mathcal{G}_{22}(t,s)\phi(x)\|\lesssim&
(1+t)^\lambda(1+s)^\lambda\cdot
\Gamma^{\frac{n}{2}+|\alpha|+2}(t,s)\cdot
(\|\phi\|_{L^1}^l+\|\partial_x^{|\alpha|+1}\phi\|^h),
\\ \nonumber
\|\partial_x^\alpha \mathcal{G}_{22}(t,s)\phi(x)\|\lesssim&
\Big(\frac{1+t}{1+s}\Big)^\lambda\cdot
\Gamma^{\frac{n}{2}+|\alpha|}(t,s)
\\ \label{eq-decay-G-can}
&\cdot
\Big(
(1+s)^{2\lambda}\cdot\Gamma^{2}(t,s)
+\frac{1}{(1+s)^{\lambda-1}}
+C_\kappa\Gamma^{\kappa}(t,s)
\Big)
\cdot
(\|\phi\|_{L^1}^l+\|\partial_x^{|\alpha|}\phi\|^h),
\end{align}
where $\kappa\ge2$ can be chosen arbitrarily large
and $C_\kappa>0$ is a constant depending on $\kappa$.
\end{lemma}
{\it\bfseries Proof.}
These estimates are simple conclusions of Theorem \ref{th-linear} in Appendix.
$\hfill\Box$

\begin{lemma} \label{le-min}
For $\beta>0$, $\gamma>0$, and $\lambda\in(-1,0)$, there holds
\begin{equation} \label{eq-min}
\int_0^t
\Gamma^{\beta}(t,s)\cdot
(1+s)^{-\gamma}ds
\lesssim
\begin{cases}
(1+t)^{-\min\{\frac{1+\lambda}{2}\beta,\gamma\}},
&\max\{\frac{1+\lambda}{2}\beta,\gamma\}>1,\\
(1+t)^{-\min\{\frac{1+\lambda}{2}\beta,\gamma\}}\cdot\ln(e+t),
&\max\{\frac{1+\lambda}{2}\beta,\gamma\}=1,\\
(1+t)^{-\gamma-\frac{1+\lambda}{2}\beta+1},
&\max\{\frac{1+\lambda}{2}\beta,\gamma\}<1.
\end{cases}
\end{equation}
\end{lemma}
{\it\bfseries Proof.}
This can be proved by dividing the interval of integration into $(0,\frac{t}{2})$
and $(\frac{t}{2},2)$. For details,  see the first part of our series of studies \cite{Ji-Mei-1}(Lemma 4.2)  for example.
$\hfill\Box$

\vskip2mm
We are now going to prove the optimal $L^1$-$L^2$ decay rates in  Theorem \ref{th-nonlinear}
for the nonlinear system \eqref{eq-vbdu}.

\vskip2mm
{\it\bfseries Proof of Theorem \ref{th-nonlinear}.}
The outline of proof is similar to that of the under-damping case $\lambda\in[0,1)$ in
Theorem 1.3 as we show in \cite{Ji-Mei-1}.
But the details are totally different.

Suppose that the local solution $(v,\boldsymbol u)$ exists for $t\in(0,T)$.
Since we are concerned with the large time behavior, we may assume that
the constant $T_0=0$ in Lemma \ref{le-decay-G}.
Denote the weighted energy function by
\begin{align*} \nonumber
E_n(\tilde t):=\sup_{t\in(0,\tilde t)}\Big\{&
\sum_{0\le|\alpha|\le [n/2]+1}
(1+t)^{\frac{1+\lambda}{4}n+\frac{1+\lambda}{2}|\alpha|}\|\partial_x^\alpha v\|,
\sum_{0\le|\alpha|\le [n/2]}
(1+t)^{\frac{1+\lambda}{4}n+\frac{1+\lambda}{2}(|\alpha|+1)-\lambda}
\|\partial_x^\alpha \boldsymbol u\|,
\\
&\sum_{|\alpha|=[n/2]+1}
(1+t)^{\frac{1+\lambda}{4}n+\omega_{|\alpha|}}
\|\partial_x^\alpha \boldsymbol u\|,
\sum_{|\alpha|=[n/2]+2}
(1+t)^{\frac{1+\lambda}{4}n+\theta_{|\alpha|}}
\|\partial_x^\alpha v\|,
\\
&\sum_{|\alpha|=[n/2]+2}
(1+t)^{\frac{1+\lambda}{4}n+\omega_{|\alpha|}}
\|\partial_x^\alpha \boldsymbol u\|,
\sum_{|\alpha|=[n/2]+3}
\|\partial_x^\alpha (v,\boldsymbol u)\|
\Big\},
\end{align*}
where $\omega_{[n/2]+1}$, $\omega_{[n/2]+2}$,
and $\theta_{[n/2]+2}$ are constants depending on $n$ and $\lambda$, and $\tilde t\in(0,T)$.
We claim that under the smallness of the initial data:
$\|(v_0,\boldsymbol u_0)\|_{L^1\cap H^{[\frac{n}{2}]+3}}\le\varepsilon_0$,
there holds
\begin{equation} \label{eq-apriori-n4}
E_n(\tilde t)\lesssim \delta_0, \quad \forall \tilde t\in(0,T),
\end{equation}
where $\varepsilon_0>0$ and $\delta_0>0$ are some small numbers  to be determined later.

The global existence and the a priori assumption \eqref{eq-apriori-n4}
are proved through the following three steps.
For the sake of simplicity, we take the case $n=3$ for example.
Other cases with $n\ge2$ follow similarly.

{\it Step 1}: Basic energy decay estimates.

According to the Duhamel principle \eqref{eq-Duhamel}
and the decay estimates of the Green matrix $\mathcal{G}(t,s)$
in Lemma \ref{le-decay-G}, we have
\begin{align*}
\|v(t)\|
\lesssim&
\|\mathcal{G}_{11}(t,0)v_0\|+\|\mathcal{G}_{12}(t,0)\boldsymbol u_0\|
+\int_0^t\|\mathcal{G}_{11}(t,s)Q_1(s)\|ds
+\int_0^t\|\mathcal{G}_{12}(t,s)Q_2(s)\|ds
\\
\lesssim&
\varepsilon_0(1+t)^{-\frac{1+\lambda}{4}n}
+\int_0^t
\Gamma^{\frac{n}{2}}(t,s)\cdot
(\|Q_1(s)\|_{L^1}^l+\|Q_1(s)\|^h)ds
\\
&+\int_0^t
(1+s)^\lambda\cdot
\Gamma^{\frac{n}{2}+1}(t,s)\cdot
(\|Q_2(s)\|_{L^1}^l+\|Q_2(s)\|^h)ds
\\
\lesssim&
\varepsilon_0(1+t)^{-\frac{1+\lambda}{4}n}
+E_n^2(t)\int_0^t
\Gamma^{\frac{n}{2}}(t,s)\cdot
(1+s)^{-\frac{1+\lambda}{2}n-1}ds
\\
&+E_n^2(t)\int_0^t
(1+s)^\lambda\cdot
\Gamma^{\frac{n}{2}+1}(t,s)\cdot
(1+s)^{-\frac{1+\lambda}{2}n-\frac{1+\lambda}{2}}ds
\\
\lesssim&
\varepsilon_0(1+t)^{-\frac{1+\lambda}{4}n}
+E_n^2(t)(1+t)^{-\frac{1+\lambda}{4}n},
\end{align*}
where we have used Lemma \ref{le-min}
(note that $\frac{1+\lambda}{2}n+\frac{1+\lambda}{2}-\lambda>1$
for all $n\ge2$ and $\lambda\in(-1,0)$) and
the following decay estimates on $\|Q(s)\|_{L^1}$ and $\|Q(s)\|$
(here and after, we use $D^j:=\partial_x^j$
and we may also write $\boldsymbol u$ as $u$ for simplicity):
\begin{align*}
\|Q_1(s)\|_{L^1}
&\lesssim \|uDv\|_{L^1}+\|vDu\|_{L^1}
\lesssim \|u\|\|Dv\|+\|v\|\|Du\|
\lesssim E_n^2(s)(1+s)^{-\frac{1+\lambda}{2}n-1},\\
\|Q_2(s)\|_{L^1}
&\lesssim \|uDu\|_{L^1}+\|vDv\|_{L^1}
\lesssim \|u\|\|Du\|+\|v\|\|Dv\|
\lesssim E_n^2(s)(1+s)^{-\frac{1+\lambda}{2}n-\frac{1+\lambda}{2}}.
\end{align*}
For $n=3$, we have
\begin{align*}
\|u(s)\|_{L^\infty}
&\lesssim \|Du\|^\frac{1}{2}\|D^2u\|^\frac{1}{2}
\lesssim E_n(s)(1+s)^{-\frac{1+\lambda}{4}n-\frac{1}{2}(1+\omega_2)},\\
\|v(s)\|_{L^\infty}
&\lesssim \|Dv\|^\frac{1}{2}\|D^2v\|^\frac{1}{2}
\lesssim E_n(s)(1+s)^{-\frac{1+\lambda}{4}n-\frac{3}{4}(1+\lambda)},\\
\|Du(s)\|_{L^\infty}
&\lesssim \|D^2u\|^\frac{1}{2}\|D^3u\|^\frac{1}{2}
\lesssim E_n(s)(1+s)^{-\frac{1+\lambda}{4}n-\frac{1}{2}(\omega_2+\omega_3)}, \\
\|Dv(s)\|_{L^\infty}
&\lesssim \|D^2v\|^\frac{1}{2}\|D^3v\|^\frac{1}{2}
\lesssim E_n(s)(1+s)^{-\frac{1+\lambda}{4}n-\frac{1}{2}(1+\lambda+\theta_3)},\\
\|D^2u(s)\|_{L^\infty}
&\lesssim \|D^3u\|^\frac{1}{2}\|D^4u\|^\frac{1}{2}
\lesssim E_n(s)(1+s)^{-\frac{1+\lambda}{8}n-\frac{1}{2}\omega_3}, \\
\|D^2v(s)\|_{L^\infty}
&\lesssim \|D^3v\|^\frac{1}{2}\|D^4v\|^\frac{1}{2}
\lesssim E_n(s)(1+s)^{-\frac{1+\lambda}{8}n-\frac{1}{2}\theta_3},
\end{align*}
and
\begin{align*}
\|Q_1(s)\|
&\lesssim \|uDv\|+\|vDu\|
\lesssim \|u\|_{L^\infty}\|Dv\|+\|v\|_{L^\infty}\|Du\|
\lesssim E_n^2(s)(1+s)^{-\frac{1+\lambda}{2}n-\frac{1+\lambda}{2}-\frac{1}{2}(1+\omega_2)}, \\
\|Q_2(s)\|
&\lesssim \|uDu\|+\|vDv\|
\lesssim \|u\|_{L^\infty}\|Du\|+\|v\|_{L^\infty}\|Dv\|
\lesssim E_n^2(s)(1+s)^{-\frac{1+\lambda}{2}n-\frac{5}{4}(1+\lambda)},\\
\|DQ_1(s)\|
&\lesssim \|DuDv\|+\|uD^2v\|+\|vD^2u\|
\lesssim E_n^2(s)(1+s)^{-\frac{1+\lambda}{2}n-\theta_{11}}, \\
\|DQ_2(s)\|
&\lesssim \|uD^2u\|+\|DuDu\|+\|vD^2v\|+\|DvDv\|
\lesssim E^2(s)(1+s)^{-\frac{1+\lambda}{2}n-\theta_{12}},
\end{align*}
where
\begin{align*}
\theta_{11}=&\min\Big\{
1+\frac{1+\lambda}{2}+\frac{\theta_3}{2},
1+\lambda+\frac{1}{2}(1+\omega_2),1+\frac{3}{4}(1+\lambda)
\Big\},\\
\theta_{12}=&\min\Big\{
\omega_2+\frac{1}{2}(1+\omega_2),1+\frac{1}{2}(\omega_2+\omega_3),
1+\lambda+\frac{3}{4}(1+\lambda),
1+\lambda+\frac{\theta_3}{2}\Big\}.
\end{align*}
Using the above estimates, we have
\begin{align*}
\|Dv(t)\|
\lesssim&
\|D\mathcal{G}_{11}(t,0)v_0\|+\|D\mathcal{G}_{12}(t,0)u_0\|
+\int_0^t\|D\mathcal{G}_{11}(t,s)Q_1(s)\|ds
+\int_0^t\|D\mathcal{G}_{12}(t,s)Q_2(s)\|ds
\\
\lesssim&
\varepsilon_0(1+t)^{-\frac{1+\lambda}{4}n-\frac{1+\lambda}{2}}
+\int_0^t
\Gamma^{\frac{n}{2}+1}(t,s)\cdot
(\|Q_1(s)\|_{L^1}+\|DQ_1(s)\|)ds
\\
&+\int_0^t
(1+s)^\lambda\cdot
\Gamma^{\frac{n}{2}+2}(t,s)\cdot
(\|Q_2(s)\|_{L^1}+\|DQ_2(s)\|)ds
\\
\lesssim&
\varepsilon_0(1+t)^{-\frac{1+\lambda}{4}n-\frac{1+\lambda}{2}}
+E_n^2(t)\int_0^t
\Gamma^{\frac{n}{2}+1}(t,s)\cdot
(1+s)^{-\frac{1+\lambda}{2}n-\min\{1,\theta_{11}\}}ds
\\
&+E_n^2(t)\int_0^t
(1+s)^\lambda\cdot
\Gamma^{\frac{n}{2}+2}(t,s)\cdot
(1+s)^{-\frac{1+\lambda}{2}n-\frac{1+\lambda}{2}}ds
\\
\lesssim&
\varepsilon_0(1+t)^{-\frac{1+\lambda}{4}n-\frac{1+\lambda}{2}}
+E_n^2(t)(1+t)^{-\frac{1+\lambda}{4}n-\frac{1+\lambda}{2}},
\end{align*}
provided that $\frac{1+\lambda}{2}n+\min\{1,\theta_{11}\}>1$.  Similarly, we also have
\begin{align*}
\|D^2v(t)\|
\lesssim&
\|D^2\mathcal{G}_{11}(t,0)v_0\|+\|D^2\mathcal{G}_{12}(t,0)u_0\|
+\int_0^t\|D^2\mathcal{G}_{11}(t,s)Q_1(s)\|ds
+\int_0^t\|D^2\mathcal{G}_{12}(t,s)Q_2(s)\|ds
\\
\lesssim&
\varepsilon_0(1+t)^{-\frac{1+\lambda}{4}n-(1+\lambda)}
+\int_0^t
\Gamma^{\frac{n}{2}+2}(t,s)\cdot
(\|Q_1(s)\|_{L^1}+\|D^2Q_1(s)\|)ds
\\
&+\int_0^t
(1+s)^\lambda\cdot
\Gamma^{\frac{n}{2}+3}(t,s)\cdot
(\|Q_2(s)\|_{L^1}+\|D^2Q_2(s)\|)ds
\\
\lesssim&
\varepsilon_0(1+t)^{-\frac{1+\lambda}{4}n-(1+\lambda)}
+E_n^2(t)\int_0^t
\Gamma^{\frac{n}{2}+2}(t,s)\cdot
(1+s)^{-\frac{1+\lambda}{2}n-\min\{1,\theta_{21}\}}ds
\\
&+E_n^2(t)\int_0^t
(1+s)^\lambda\cdot
\Gamma^{\frac{n}{2}+3}(t,s)\cdot
(1+s)^{-\frac{1+\lambda}{2}n-\min\{\frac{1+\lambda}{2},\theta_{22}\}}ds
\\
\lesssim&
\varepsilon_0(1+t)^{-\frac{1+\lambda}{4}n-(1+\lambda)}
+E_n^2(t)(1+t)^{-\frac{1+\lambda}{4}n-(1+\lambda)},
\end{align*}
provided that
\begin{equation} \label{eq-zlambda}
\begin{cases}
\frac{1+\lambda}{2}n+\min\{1,\theta_{21}\}>1, \quad
&\frac{1+\lambda}{2}n+\min\{1,\theta_{21}\}\ge
\frac{1+\lambda}{4}n+(1+\lambda),
\\
\frac{1+\lambda}{2}n+\min\{\frac{1+\lambda}{2},\theta_{22}\}-\lambda>1, \quad
&\frac{1+\lambda}{2}n+\min\{\frac{1+\lambda}{2},\theta_{22}\}-\lambda\ge
\frac{1+\lambda}{4}n+(1+\lambda),
\end{cases}
\end{equation}
where we have also used the following estimates
\begin{align*}
\|D^2Q_1(s)\|
&\lesssim \|uD^3v\|+\|DuD^2v\|+\|DvD^2u\|+\|vD^3u\|
\lesssim E_n^2(s)(1+s)^{-\frac{1+\lambda}{2}n-\theta_{21}}, \\
\|D^2Q_2(s)\|
&\lesssim \|uD^3u\|+\|DuD^2u\|+\|vD^3v\|+\|DvD^2v\|
\lesssim E_n^2(s)(1+s)^{-\frac{1+\lambda}{2}n-\theta_{22}},
\end{align*}
with
\begin{align*}
\theta_{21}&=\min\Big\{
\frac{1}{2}(1+\omega_2)+\theta_3,
1+\lambda+\frac{1}{2}(\omega_2+\omega_3),
1+\frac{1}{2}(1+\lambda+\theta_3),
\omega_3+\frac{3}{4}(1+\lambda)\Big\},
\\
\theta_{22}&=\min\Big\{
\frac{1}{2}(1+\omega_2)+\omega_3,
\omega_2+\frac{1}{2}(\omega_2+\omega_3),\theta_3+\frac{3}{4}(1+\lambda),
1+\lambda+\frac{1}{2}(1+\lambda+\omega_3)\Big\}.
\end{align*}

The decay estimates on
$\|\partial_x^\alpha v\|$ for $0\le|\alpha|\le [\frac{n}{2}]+1$
are based on the optimal decay estimates on
$\|\partial_x^\alpha G_{11}(t,s)\|$ and $\|\partial_x^\alpha G_{12}(t,s)\|$
in \eqref{eq-decay-G}.
However, the estimates on
$\|\partial_x^\alpha G_{21}(t,s)\|$ and $\|\partial_x^\alpha G_{22}(t,s)\|$
in \eqref{eq-decay-G} are insufficient for the optimal decay estimates
on $\|\partial_x^\alpha \boldsymbol u\|$ for $0\le|\alpha|\le [\frac{n}{2}]$.
In fact, we use the optimal decay estimates in \eqref{eq-decay-G-opt}
to show the decay estimates on
$\|\partial_x^\alpha \boldsymbol u\|$ for $0\le|\alpha|\le [\frac{n}{2}]$
in a similar way as $\|\partial_x^\alpha v\|$ for $1\le|\alpha|\le [\frac{n}{2}]+1$.
One can check that the condition on the estimate of
$\|\partial_x^k \boldsymbol u\|$ for $0\le k\le [\frac{n}{2}]$
is equivalent to the condition on the estimate of
$\|\partial_x^{k+1} v\|$.
For example,
\begin{align*}
\|\boldsymbol u(t)\|
\lesssim&
\|\mathcal{G}_{21}(t,0)v_0\|+\|\mathcal{G}_{22}(t,0)\boldsymbol u_0\|
+\int_0^t\|\mathcal{G}_{21}(t,s)Q_1(s)\|ds
+\int_0^t\|\mathcal{G}_{22}(t,s)Q_2(s)\|ds
\\
\lesssim&
\varepsilon_0(1+t)^{-\frac{1+\lambda}{4}n}
+\int_0^t
(1+t)^\lambda\cdot
\Gamma^{\frac{n}{2}+1}(t,s)\cdot
(\|Q_1(s)\|_{L^1}^l+\|DQ_1(s)\|^h)ds
\\
&+\int_0^t
(1+t)^\lambda(1+s)^\lambda\cdot
\Gamma^{\frac{n}{2}+2}(t,s)\cdot
(\|Q_2(s)\|_{L^1}^l+\|DQ_2(s)\|^h)ds
\\
\lesssim&
\varepsilon_0(1+t)^{-\frac{1+\lambda}{4}n-\frac{1-\lambda}{2}}
+E_n^2(t)\int_0^t
(1+t)^\lambda\cdot
\Gamma^{\frac{n}{2}+1}(t,s)\cdot
(1+s)^{-\frac{1+\lambda}{2}n-1}ds
\\
&+E_n^2(t)\int_0^t
(1+t)^\lambda(1+s)^\lambda\cdot
\Gamma^{\frac{n}{2}+2}(t,s)\cdot
(1+s)^{-\frac{1+\lambda}{2}n-\frac{1+\lambda}{2}}ds
\\
\lesssim&
\varepsilon_0(1+t)^{-\frac{1+\lambda}{4}n-\frac{1-\lambda}{2}}
+E_n^2(t)(1+t)^{-\frac{1+\lambda}{4}n-\frac{1-\lambda}{2}}.
\end{align*}

Further, we use the decay estimates in \eqref{eq-decay-G-can}
to show the decay estimates on
$\|\partial_x^\alpha \boldsymbol u\|$ for $[\frac{n}{2}]+1\le|\alpha|\le [\frac{n}{2}]+2$
since the regularity required in \eqref{eq-decay-G-can}
is one-order lower than that in \eqref{eq-decay-G-opt}.
We note that in this case the condition on the estimate of
$\|\partial_x^k \boldsymbol u\|$ for $[\frac{n}{2}]+1\le k\le [\frac{n}{2}]+2$
is similar to the condition on the estimate of
$\|\partial_x^k v\|$.
We have
\begin{align*}
\|D^3Q_1(s)\|
&\lesssim \|uD^4v\|+\|DuD^3v\|+\|D^2uD^2v\|+\|DvD^3u\|+\|vD^4u\|
\lesssim E_n^2(s)(1+s)^{-\frac{1+\lambda}{4}n-\theta_{31}},
\\
\|D^3Q_2(s)\|
&\lesssim \|uD^4u\|+\dots+\|D^2uD^2u\|+\|vD^4v\|+\dots+\|D^2vD^2v\|
\lesssim E_n^2(s)(1+s)^{-\frac{1+\lambda}{4}n-\theta_{32}},
\end{align*}
with
\begin{align*}
\theta_{31}=\min\Big\{&
\frac{1}{2}(1+\omega_2),\frac{1}{2}(\omega_2+\omega_3)+\frac{1+\lambda}{4}n+\theta_3,
\frac{1+\lambda}{8}n+1+\frac{\theta_3}{2},
\\
&\qquad\frac{1+\lambda}{4}n+\frac{1}{2}(1+\lambda+\theta_3)+\omega_3,
\frac{3}{4}(1+\lambda)\Big\},
\\
\theta_{32}=\min\Big\{&
\frac{1}{2}(1+\omega_2),\frac{1}{2}(\omega_2+\omega_3)+\frac{1+\lambda}{4}n+\omega_3,
\frac{1+\lambda}{8}n+1+\frac{\omega_3}{2},
\\
&\qquad\frac{3}{4}(1+\lambda),
\frac{1+\lambda}{4}n+\frac{1}{2}(1+\lambda+\theta_3)+\theta_3,
\frac{1+\lambda}{8}n+1+\lambda+\frac{\theta_3}{2}
\Big\}.
\end{align*}
Therefore, we arrive at
\begin{align*}
\|D^2u(t)\|
\lesssim&
\|D^2\mathcal{G}_{21}(t,0)v_0\|+\|D^2\mathcal{G}_{22}(t,0)u_0\|
+\int_0^t\|D^2\mathcal{G}_{21}(t,s)Q_1(s)\|ds
+\int_0^t\|D^2\mathcal{G}_{22}(t,s)Q_2(s)\|ds
\\
\lesssim&
\varepsilon_0(1+t)^{-\frac{1+\lambda}{4}n-(1+\lambda)+\lambda}
+\int_0^t(1+t)^\lambda\cdot
\Gamma^{\frac{n}{2}+3}(t,s)\cdot
(\|Q_1(s)\|_{L^1}+\|D^2Q_1(s)\|)ds
\\
&+\int_0^t
\Big(\frac{1+t}{1+s}\Big)^\lambda\cdot
\Gamma^{\frac{n}{2}+2}(t,s)\cdot
\Big(
(1+s)^{2\lambda}\cdot\Gamma^{2}(t,s)
+\frac{1}{(1+s)^{\lambda-1}}
+C_\kappa\Gamma^{\kappa}(t,s)
\Big)
\\
&\qquad
\cdot
(\|Q_2(s)\|_{L^1}+\|D^2Q_2(s)\|)ds
\\
\lesssim&
\varepsilon_0(1+t)^{-\frac{1+\lambda}{4}n-1}
\\
&+E_n^2(t)\int_0^t(1+t)^\lambda\cdot
\Gamma^{\frac{n}{2}+3}(t,s)\cdot
(1+s)^{-\frac{1+\lambda}{2}n-\min\{1,\theta_{21}\}}ds
\\
&+E_n^2(t)\int_0^t
\Big(\frac{1+t}{1+s}\Big)^\lambda\cdot
\Gamma^{\frac{n}{2}+2}(t,s)
\cdot
\Big(
(1+s)^{2\lambda}\cdot\Gamma^{2}(t,s)
+\frac{1}{(1+s)^{\lambda-1}}
+C_\kappa\Gamma^{\kappa}(t,s)
\Big)
\\
&\qquad\cdot
(1+s)^{-\frac{1+\lambda}{2}n-\min\{\frac{1+\lambda}{2},\theta_{22}\}}ds
\\
\lesssim&
\varepsilon_0(1+t)^{-\frac{1+\lambda}{4}n-\omega_2}
+E_n^2(t)(1+t)^{-\frac{1+\lambda}{4}n-\omega_2},
\end{align*}
provided that
\begin{equation} \label{eq-zlambda2}
\begin{cases}
\frac{1+\lambda}{2}n+\min\{1,\theta_{21}\}-\lambda\ge \frac{1+\lambda}{4}n+\omega_2,
\\
\frac{1+\lambda}{2}n+\min\{\frac{1+\lambda}{2},\theta_{22}\}\ge \frac{1+\lambda}{4}n+\omega_2.
\end{cases}
\end{equation}
Furthermore, we similarly have
\begin{align*}
\|D^3u(t)\|
\lesssim&
\|D^3\mathcal{G}_{21}(t,0)v_0\|+\|D^3\mathcal{G}_{22}(t,0)u_0\|
+\int_0^t\|D^3\mathcal{G}_{21}(t,s)Q_1(s)\|ds
+\int_0^t\|D^3\mathcal{G}_{22}(t,s)Q_2(s)\|ds
\\
\lesssim&
\varepsilon_0(1+t)^{-\frac{1+\lambda}{4}n-\frac{3}{2}(1+\lambda)+\lambda}
\\
&+\int_0^t(1+t)^\lambda\cdot
\Gamma^{\frac{n}{2}+4}(t,s)\cdot
(\|Q_1(s)\|_{L^1}+\|D^3Q_1(s)\|)ds
\\
&+\int_0^t
\Big(\frac{1+t}{1+s}\Big)^\lambda\cdot
\Gamma^{\frac{n}{2}+3}(t,s)\cdot
\Big(
(1+s)^{2\lambda}\cdot\Gamma^{2}(t,s)
+\frac{1}{(1+s)^{\lambda-1}}
+C_\kappa\Gamma^{\kappa}(t,s)
\Big)
\\
&\qquad\cdot
(\|Q_2(s)\|_{L^1}+\|D^3Q_2(s)\|)ds
\\
\lesssim&
\varepsilon_0(1+t)^{-\frac{1+\lambda}{4}n-\frac{3}{2}(1+\lambda)+\lambda}
\\
&+E_n^2(t)\int_0^t(1+t)^\lambda\cdot
\Gamma^{\frac{n}{2}+4}(t,s)\cdot
(1+s)^{-\min\{\frac{1+\lambda}{2}n+1,\frac{1+\lambda}{4}n+\theta_{31}\}}ds
\\
&+E_n^2(t)\int_0^t
\Big(\frac{1+t}{1+s}\Big)^\lambda\cdot
\Gamma^{\frac{n}{2}+3}(t,s)\cdot
\Big(
(1+s)^{2\lambda}\cdot\Gamma^{2}(t,s)
+\frac{1}{(1+s)^{\lambda-1}}
+C_\kappa\Gamma^{\kappa}(t,s)
\Big)
\\
&\qquad\cdot
(1+s)^{-\min\{\frac{1+\lambda}{2}n+\frac{1+\lambda}{2},\frac{1+\lambda}{4}n+\theta_{32}\}}ds
\\
\lesssim&
\varepsilon_0(1+t)^{-\frac{1+\lambda}{4}n-\omega_3}
+E_n^2(t)(1+t)^{-\frac{1+\lambda}{4}n-\omega_3},
\end{align*}
provided that
\begin{equation} \label{eq-zlambda3}
\begin{cases}
\frac{3}{2}(1+\lambda)-\lambda\ge \omega_3,
\\
\frac{1+\lambda}{4}n+\theta_{31}-\lambda\ge \frac{1+\lambda}{4}n+\omega_3,
\\
\frac{1+\lambda}{4}n+\theta_{32}\ge \frac{1+\lambda}{4}n+\omega_3.
\end{cases}
\end{equation}
Hence, the estimate on $\|D^3v\|$ is
\begin{align*}
\|D^3v(t)\|
\lesssim&
\|D^3\mathcal{G}_{11}(t,0)v_0\|+\|D^3\mathcal{G}_{12}(t,0)u_0\|
+\int_0^t\|D^3\mathcal{G}_{11}(t,s)Q_1(s)\|ds
+\int_0^t\|D^3\mathcal{G}_{12}(t,s)Q_2(s)\|ds
\\
\lesssim&
\varepsilon_0(1+t)^{-\frac{1+\lambda}{4}n-\frac{3}{2}(1+\lambda)}
\\
&+\int_0^t
\Gamma^{\frac{n}{2}+3}(t,s)\cdot
(\|Q_1(s)\|_{L^1}+\|D^3Q_1(s)\|)ds
\\
&+\int_0^t
(1+s)^\lambda\cdot
\Gamma^{\frac{n}{2}+4}(t,s)\cdot
(\|Q_2(s)\|_{L^1}+\|D^3Q_2(s)\|)ds
\\
\lesssim&
\varepsilon_0(1+t)^{-\frac{1+\lambda}{4}n-\frac{3}{2}(1+\lambda)}
\\
&+E_n^2(t)\int_0^t
\Gamma^{\frac{n}{2}+3}(t,s)\cdot
(1+s)^{-\min\{\frac{1+\lambda}{2}n+1,\frac{1+\lambda}{4}n+\theta_{31}\}}ds
\\
&+E_n^2(t)\int_0^t
(1+s)^\lambda\cdot
\Gamma^{\frac{n}{2}+4}(t,s)\cdot
(1+s)^{-\min\{\frac{1+\lambda}{2}n+\frac{1+\lambda}{2},\frac{1+\lambda}{4}n+\theta_{32}\}}ds
\\
\lesssim&
\varepsilon_0(1+t)^{-\frac{1+\lambda}{4}n-\theta_3}
+E_n^2(t)(1+t)^{-\frac{1+\lambda}{4}n-\theta_3},
\end{align*}
under the condition
\begin{equation} \label{eq-zlambda4}
\begin{cases}
\frac{3}{2}(1+\lambda)\ge \theta_3,
\\
\frac{1+\lambda}{4}n+\theta_{31}
+\frac{1+\lambda}{2}(\frac{n}{2}+3)-1
\ge
\frac{1+\lambda}{4}n+\theta_3,
\\
\frac{1+\lambda}{4}n+\theta_{32}-\lambda>1.
\end{cases}
\end{equation}

Combining the above conditions together, we fix
$\theta_3=-\frac{1+\lambda}{3}$, $\omega_3=\frac{1+\lambda}{4}\cdot3$,
and $\omega_2=\frac{1+\lambda}{4}n+\frac{1+\lambda}{2}$ for the case $n=3$.
We note that the restriction on $\omega_3$ is \eqref{eq-zlambda3}$_3$
such that $\omega_3\le\theta_{32}$
and $\frac{1+\lambda}{4}n+\theta_{32}$ is the decay rate of
$\|D^3Q_2(s)\|$, where the worst term (decaying slowest) is $\|vD^4v\|$
restricted by $\|v\|_{L^\infty}$.
For general dimension $n$, we have
\begin{align*}
\|D^{[\frac{n}{2}]+2}Q_2(s)\|\lesssim &
\|vD^{[\frac{n}{2}]+3}v(s)\|
\lesssim \|v\|_{L^\infty}\cdot E_n(s)
\\
\lesssim &
\begin{cases}
\|D^{[\frac{n}{2}]}v\|^\frac{1}{2}\cdot\|D^{[\frac{n}{2}]+1}v\|^\frac{1}{2}\cdot E_n(s),
\quad & \text{for~odd~}n,
\\
\|D^{[\frac{n}{2}]-1}v\|^\frac{1}{2}\cdot\|D^{[\frac{n}{2}]+1}v\|^\frac{1}{2}\cdot E_n(s),
\quad & \text{for~even~}n,
\end{cases}
\\
\lesssim &
\begin{cases}
E_n^2(s)\cdot\big((1+s)^{-\frac{1+\lambda}{4}n-\frac{1+\lambda}{2}[\frac{n}{2}]}
\cdot(1+s)^{-\frac{1+\lambda}{4}n-\frac{1+\lambda}{2}([\frac{n}{2}]+1)}\big)^\frac{1}{2},
\quad & \text{for~odd~}n,
\\
E_n^2(s)\cdot\big((1+s)^{-\frac{1+\lambda}{4}n-\frac{1+\lambda}{2}([\frac{n}{2}]-1)}
\cdot(1+s)^{-\frac{1+\lambda}{4}n-\frac{1+\lambda}{2}([\frac{n}{2}]+1)}\big)^\frac{1}{2},
\quad & \text{for~even~}n,
\end{cases}
\\
\approx &
E_n^2(s)\cdot(1+s)^{-\frac{1+\lambda}{4}n-\frac{1+\lambda}{4}n}.
\end{align*}
Therefore, it suffices to take $\omega_{[\frac{n}{2}]+2}=\frac{1+\lambda}{4}n$
for general dimension of $n$.
The condition \eqref{eq-zlambda}$_1$, which is necessary for the optimal decay
of $\|D^2v\|$, is $\frac{1+\lambda}{2}n+\omega_3+\frac{1+\lambda}{4}\cdot3>1$ for $n=3$,
and is $\frac{1+\lambda}{2}n+\frac{1+\lambda}{4}n+\frac{1+\lambda}{4}n>1$ for general $n\ge2$.
That is, $(1+\lambda)n>1$, which is equivalent to $\lambda\in(-\frac{n-1}{n},0)$.
The condition $\lambda\in(-\frac{n}{n+2},0)$
is stronger than $\lambda\in(-\frac{n-1}{n},0)$.

{\it Step 2}: High-order energy estimates.

We note that the condition \eqref{eq-apriori-energy} in Lemma \ref{le-energy}
under the a priori assumption \eqref{eq-apriori-n4} is
$\frac{1+\lambda}{4}n>-\lambda$,
which is true for $\lambda\in(-\frac{n}{n+4},0)$
and is false for $\lambda\in(-1,-\frac{n}{n+4})$.
Fortunately, the condition \eqref{eq-apriori-energy-G} in Lemma \ref{le-energy-G}
under the a priori assumption \eqref{eq-apriori-n4} is
$$
\frac{1+\lambda}{4}n+\min\Big\{\frac{1+\lambda}{2}-\lambda,
\omega_{[n/2]+2}\Big\}
=\min\Big\{\frac{1+\lambda}{4}n+\frac{1+\lambda}{2}-\lambda,
\frac{1+\lambda}{2}n
\Big\}
>-\lambda,
$$
which is true for all $\lambda\in(-\frac{n}{n+2},0)$.
Therefore, we can apply the high-order energy estimates of
Lemma \ref{eq-energy-high-G} to get
\begin{equation} \label{eq-zenergy-high}
\|\nabla(v,\boldsymbol u)\|_{H^{[\frac{n}{2}]+2}}^2
\lesssim \|\nabla(v_0,\boldsymbol u_0)\|_{H^{[\frac{n}{2}]+2}}^2
+\delta_0^2\omega^2(t),
\end{equation}
where $\delta_0\omega(t)=(1+t)^{
-\frac{1+\lambda}{2}n-\lambda}$
decays to zero.

{\it Step 3}: Closure of the a priori estimate \eqref{eq-apriori-n4}.

We now combine the above estimates and choose $\varepsilon_0>0$ and $\delta_0>0$
to be sufficiently small such that
$$C(\varepsilon_0+\delta_0^2+\delta_0\omega(t))\le \delta_0,$$
where $C>0$ is a universal constant.
It suffices to choose $C\delta_0\le 1/4$, and $C\varepsilon_0=\delta_0/2$,
and to consider the problem starting form $t_0$ such that $C\omega(t_0)\le 1/4$
since $\omega(t)$ decays to zero.
We see that the a priori estimate
\eqref{eq-apriori-n4} holds for all the time $t\in(0,+\infty)$.

Finally, we show that those estimates
($\|\partial_x^\alpha v\|$ with $0\le|\alpha|\le[\frac{n}{2}]+1$ and
$\|\partial_x^\alpha \boldsymbol u\|$ with $0\le|\alpha|\le[\frac{n}{2}]$)
are optimal.
We take the estimate on $\|v\|$ for example.
According to the optimal decay estimates in Lemma \ref{le-decay-G} and
the energy estimates in Step 1 before,
we choose the initial data $(v_0,\boldsymbol u_0)$ such that
$\|\mathcal{G}_{11}(t,0)v_0\|$ decays optimally, then we have
$$
\|v(t)\|
\gtrsim
\|\mathcal{G}_{11}(t,0)v_0\|-\|\mathcal{G}_{12}(t,0)\boldsymbol u_0\|
-\int_0^t\|\mathcal{G}_{11}(t,s)Q_1(s)\|ds
-\int_0^t\|\mathcal{G}_{12}(t,s)Q_2(s)\|ds,
$$
where $\|\mathcal{G}_{12}(t,0)\boldsymbol u_0\|$ decays faster than
$\|\mathcal{G}_{11}(t,0)v_0\|$,
and $\int_0^t\|\mathcal{G}_{11}(t,s)Q_1(s)\|ds
+\int_0^t\|\mathcal{G}_{12}(t,s)Q_2(s)\|ds$ decays
no slower than $\|\mathcal{G}_{11}(t,0)v_0\|$.
We note that $Q_1(t,x)$ and $Q_2(t,x)$ are quadratic,
and we rescale the initial data as $(\varepsilon_1v_0,\varepsilon_1\boldsymbol u_0)$
with $\varepsilon_1>0$ sufficiently small
such that neither
$\int_0^t\|\mathcal{G}_{11}(t,s)Q_1(s)\|ds$ nor
$\int_0^t\|\mathcal{G}_{12}(t,s)Q_2(s)\|ds$
is comparable with $\|\mathcal{G}_{11}(t,0)v_0\|$.
In fact, according to the proof in Step 1, we have
$$\|\mathcal{G}_{11}(t,0)\varepsilon_1v_0\|
\approx \varepsilon_0(1+t)^{-\frac{1+\lambda}{4}n},$$
and
$$\int_0^t\|\mathcal{G}_{11}(t,s)Q_1(s)\|ds
+\int_0^t\|\mathcal{G}_{12}(t,s)Q_2(s)\|ds
\lesssim E_n^2(t)(1+t)^{-\frac{1+\lambda}{4}n}
\lesssim \delta_0^2(1+t)^{-\frac{1+\lambda}{4}n}
\lesssim \varepsilon_0^2(1+t)^{-\frac{1+\lambda}{4}n},$$
even though they are nonlinear.
Therefore, $\|v(t)\|$ decays in the same order as
$\|\mathcal{G}_{11}(t,0)v_0\|$.
The proof is completed.
$\hfill\Box$

\subsection{Optimal $L^q$ decay estimates}

We now turn to the $L^1$-$L^q$ decay estimates of the nonlinear system \eqref{eq-vbdu}.

\begin{lemma} \label{le-decay-G-h}
For $q\in[2,\infty]$ and $\lambda\in(-1,0)$, then
\begin{align*}
&\|\partial_x^\alpha \mathcal{G}_{11}(t,s)\phi(x)\|_{L^q}\lesssim
\Gamma^{\gamma_{1,q}+|\alpha|}(t,s)\cdot
(\|\phi\|_{L^{1}}^l+\|\partial_x^{|\alpha|+\omega_{2,q}}\phi\|^h),
\\
&\|\partial_x^\alpha \mathcal{G}_{12}(t,s)\phi(x)\|_{L^q}\lesssim
(1+s)^\lambda\cdot
\Gamma^{\gamma_{1,q}+|\alpha|+1}(t,s)\cdot
(\|\phi\|_{L^{1}}^l+\|\partial_x^{|\alpha|+\omega_{2,q}}\phi\|^h),
\\
&\|\partial_x^\alpha \mathcal{G}_{21}(t,s)\phi(x)\|_{L^q}\lesssim
(1+t)^\lambda\cdot
\Gamma^{\gamma_{1,q}+|\alpha|+1}(t,s)\cdot
(\|\phi\|_{L^{1}}^l+\|\partial_x^{|\alpha|+\omega_{2,q}}\phi\|^h),
\\
&\|\partial_x^\alpha \mathcal{G}_{22}(t,s)\phi(x)\|_{L^q}\lesssim
\Big(\frac{1+t}{1+s}\Big)^\lambda\cdot
\Gamma^{\gamma_{1,q}+|\alpha|}(t,s)\cdot
(\|\phi\|_{L^{1}}^l+\|\partial_x^{|\alpha|+\omega_{2,q}}\phi\|^h),
\end{align*}
where $\gamma_{1,q}:=n(1-1/q)$,
and $\omega_{2,q}>\gamma_{2,q}:=n(1/2-1/q)$.
Furthermore, it holds
\begin{align*}
\|\partial_x^\alpha \mathcal{G}_{22}(t,s)\phi(x)\|_{L^q}\lesssim&
(1+t)^\lambda(1+s)^\lambda\cdot
\Gamma^{\gamma_{1,q}+|\alpha|+2}(t,s)\cdot
(\|\phi\|_{L^{1}}^l+\|\partial_x^{|\alpha|+1+\omega_{2,q}}\phi\|^h),
\\
\|\partial_x^\alpha \mathcal{G}_{22}(t,s)\phi(x)\|\lesssim&
\Big(\frac{1+t}{1+s}\Big)^\lambda\cdot
\Gamma^{\gamma_{1,q}+|\alpha|}(t,s)
\\
&\cdot
\Big(
(1+s)^{2\lambda}\cdot\Gamma^{2}(t,s)
+\frac{1}{(1+s)^{\lambda-1}}
+C_\kappa\Gamma^{\kappa}(t,s)
\Big)
\cdot
(\|\phi\|_{L^{1}}^l+\|\partial_x^{|\alpha|+1+\omega_{2,q}}\phi\|^h),
\end{align*}
where $\kappa\ge2$ can be chosen arbitrarily large
and $C_\kappa>0$ is a constant depending on $\kappa$.
\end{lemma}
{\it\bfseries Proof.}
These estimates are conclusions of Theorem \ref{th-linear} in Appendix.
$\hfill\Box$

\vskip2mm
We prove the optimal $L^q$ decay estimates Theorem \ref{th-nonlinear-Lp}
of the nonlinear system \eqref{eq-vbdu}.

\vskip2mm
{\it\bfseries Proof of Theorem \ref{th-nonlinear-Lp}.}
Since $\lambda\in(-\frac{n}{n+2},0)$ satisfies the condition in Theorem \ref{th-nonlinear},
we see that the a priori assumption \eqref{eq-apriori-n4}
in the proof of Theorem \ref{th-nonlinear} is valid,
which is based on the smallness of the initial data
$\|(v_0,\boldsymbol u_0)\|_{L^1\cap H^{[\frac{n}{2}]+3}}\le\varepsilon_0$.
Here under the stronger condition
$\|(v_0,\boldsymbol u_0)\|_{L^1\cap H^{[\frac{n}{2}]+k}}\le\varepsilon_0$,
we can enforce the decay estimates as follows.
Denote the new weighted energy function by
\begin{align*} \nonumber
F_n(\tilde t):=\sup_{t\in(0,\tilde t)}\Big\{&
\sum_{0\le|\alpha|\le [n/2]+1}
(1+t)^{\frac{1+\lambda}{4}n+\frac{1+\lambda}{2}|\alpha|}\|\partial_x^\alpha v\|,
\sum_{0\le|\alpha|\le [n/2]}
(1+t)^{\frac{1+\lambda}{4}n+\frac{1+\lambda}{2}(|\alpha|+1)-\lambda}
\|\partial_x^\alpha \boldsymbol u\|,
\\
&\sum_{|\alpha|=[n/2]+1}
(1+t)^{\frac{1+\lambda}{4}n+\omega_{|\alpha|}}
\|\partial_x^\alpha \boldsymbol u\|,
\sum_{[n/2]+2\le|\alpha|\le[n/2]+k-1}
(1+t)^{\frac{1+\lambda}{4}n+\theta_{|\alpha|}}
\|\partial_x^\alpha v\|,
\\
&\sum_{[n/2]+2\le|\alpha|\le[n/2]+k-1}
(1+t)^{\frac{1+\lambda}{4}n+\omega_{|\alpha|}}
\|\partial_x^\alpha \boldsymbol u\|,
\sum_{|\alpha|=[n/2]+k}
\|\partial_x^\alpha (v,\boldsymbol u)\|
\Big\},
\end{align*}
where $\omega_{|\alpha|}$ and $\theta_{|\alpha|}$ are constants depending on $n$ and $\lambda$.
We claim that under the small initial data condition
$\|(v_0,\boldsymbol u_0)\|_{L^1\cap H^{[\frac{n}{2}]+k}}\le\varepsilon_0$,
there holds
\begin{equation} \label{eq-apriori-Lq}
F_n(\tilde t)\lesssim \delta_0, \quad \forall \tilde t\in(0,T),
\end{equation}
where $\varepsilon_0>0$ and $\delta_0>0$ are small constants to be determined.

We take, for example,  the case $n=3$ again.
Note that for $n=3$, $\gamma_{2,q}=3(1/2-1/q)\le3/2<2$.
We take $k=3+[\gamma_{2,q}]=4$ and $\omega_{2,q}=2>\gamma_{2,q}$.
We prove the estimate on $\|\partial_x^\alpha v\|_{L^q}$ with $|\alpha|=1$
in \eqref{eq-nonlinear-Lp}.
According to the Duhamel principle \eqref{eq-Duhamel} and
the $L^1$-$L^q$ decay estimates of the Green matrix in Lemma \ref{le-decay-G-h},
we have
\begin{align*}
&\|Dv(t)\|_{L^q}
\\
\lesssim&
\|D\mathcal{G}_{11}(t,0)v_0\|_{L^q}+\|D\mathcal{G}_{12}(t,0)u_0\|_{L^q}
+\int_0^t\|D\mathcal{G}_{11}(t,s)Q_1(s)\|_{L^q}ds
+\int_0^t\|D\mathcal{G}_{12}(t,s)Q_2(s)\|_{L^q}ds
\\
\lesssim&
\varepsilon_0(1+t)^{-\frac{1+\lambda}{2}\gamma_{1,q}-\frac{1+\lambda}{2}}
+\int_0^t
\Gamma^{\gamma_{1,q}+1}(t,s)\cdot
(\|Q_1(s)\|_{L^1}+\|D^{1+\omega_{2,q}}Q_1(s)\|)ds
\\
&+\int_0^t
(1+s)^\lambda\cdot
\Gamma^{\gamma_{1,q}+2}(t,s)\cdot
(\|Q_2(s)\|_{L^1}+\|D^{1+\omega_{2,q}}Q_2(s)\|)ds
\\
\lesssim&
\varepsilon_0(1+t)^{-\frac{1+\lambda}{2}\gamma_{1,q}-\frac{1+\lambda}{2}}
+F_n^2(t)\int_0^t
\Gamma^{\gamma_{1,q}+1}(t,s)\cdot
(1+s)^{-\frac{1+\lambda}{2}n-\min\{1,\theta_{31}\}}ds
\\
&+F_n^2(t)\int_0^t
(1+s)^\lambda\cdot
\Gamma^{\gamma_{1,q}+2}(t,s)\cdot
(1+s)^{-\frac{1+\lambda}{2}n-\min\{\frac{1+\lambda}{2},\theta_{32}\}}ds
\\
\lesssim&
\varepsilon_0(1+t)^{-\frac{1+\lambda}{2}\gamma_{1,q}-\frac{1+\lambda}{2}}
+F_n^2(t)(1+t)^{-\frac{1+\lambda}{2}\gamma_{1,q}-\frac{1+\lambda}{2}},
\end{align*}
with $\omega_{2,q}=2>\gamma_{2,q}$ provided that
\begin{equation} \label{eq-zlambda-Lq}
\begin{cases}
\frac{1+\lambda}{2}n+\min\{1,\theta_{31}\}\ge
\frac{1+\lambda}{2}\gamma_{1,q}+\frac{1+\lambda}{2},
\qquad
&\frac{1+\lambda}{2}n+\min\{1,\theta_{31}\}>1,
\\
\frac{1+\lambda}{2}n+\min\{\frac{1+\lambda}{2},\theta_{32}\}-\lambda\ge
\frac{1+\lambda}{2}\gamma_{1,q}+\frac{1+\lambda}{2},
\qquad
&\frac{1+\lambda}{2}n+\min\{\frac{1+\lambda}{2},\theta_{32}\}-\lambda>1.
\end{cases}
\end{equation}
Here, $\frac{1+\lambda}{2}n+\theta_{31}$ and
$\frac{1+\lambda}{2}n+\theta_{32}$ are the decay rates of
$\|D^3Q_1\|$ and $\|D^3Q_2\|$
under the a priori assumption \eqref{eq-apriori-Lq}
(which is stronger than \eqref{eq-apriori-n4})
such that
\begin{align*}
\|D^3Q_1(s)\|
&\lesssim \|uD^4v\|+\|DuD^3v\|+\|D^2uD^2v\|+\|DvD^3u\|+\|vD^4u\|
\lesssim F_n^2(s)(1+s)^{-\frac{1+\lambda}{2}n-\theta_{31}},
\\
\|D^3Q_2(s)\|
&\lesssim \|uD^4u\|+\dots+\|D^2uD^2u\|+\|vD^4v\|+\dots+\|D^2vD^2v\|
\lesssim F_n^2(s)(1+s)^{-\frac{1+\lambda}{2}n-\theta_{32}},
\end{align*}
with
\begin{align*}
\theta_{31}=\min\Big\{&
\frac{1}{2}(1+\omega_2)+\theta_4,\frac{1}{2}(\omega_2+\omega_3)+\theta_3,
\omega_2+\frac{\theta_3+\theta_4}{2},
\\
&\qquad\frac{1}{2}(1+\lambda+\theta_3)+\omega_3,
\frac{3}{4}(1+\lambda)+\omega_4\Big\},
\\
\theta_{32}=\min\Big\{&
\frac{1}{2}(1+\omega_2)+\omega_4,\frac{1}{2}(\omega_2+\omega_3)+\omega_3,
\omega_2+\frac{1}{2}(\omega_3+\omega_4),
\\
&\qquad\frac{3}{4}(1+\lambda)+\theta_4,
\frac{1}{2}(1+\lambda+\theta_3)+\theta_3,
1+\lambda+\frac{1}{2}(\theta_3+\theta_4)
\Big\}.
\end{align*}
Similar to the proof of Theorem \ref{th-nonlinear},
where $\|D^3Q\|$ decays at the same rate as $\|(v,u)\|_{L^\infty}$
since the energy $\|D^4(v,u)\|$ is only bounded,
here $\|D^4Q\|$ decays at the same rate as $\|(v,u)\|_{L^\infty}$
due to the boundedness of $\|D^5(v,u)\|$ and
we take $\omega_4=\frac{1+\lambda}{4}n$.
Then \eqref{eq-zlambda-Lq} is equivalent to
$\frac{1+\lambda}{2}n+\frac{1+\lambda}{4}n+\frac{1+\lambda}{4}n>1$,
that is, $(1+\lambda)n>1$.
The condition $\lambda>-\frac{n}{n+2}$ is stronger than $(1+\lambda)n>1$ for all $n\ge2$.

The high-order energy estimate is similar to the Step 2 in the proof of Theorem
\ref{th-nonlinear}, where the restriction is
the condition \eqref{eq-apriori-energy-G} in Lemma \ref{le-energy-G}.
Now it reads as
$$
\frac{1+\lambda}{4}n+\min\Big\{\frac{1+\lambda}{2}-\lambda,
\omega_{[n/2]+k-2}\Big\}
=\min\Big\{\frac{1+\lambda}{4}n+\frac{1+\lambda}{2}-\lambda,
\frac{1+\lambda}{2}n\Big\}>-\lambda,
$$
under the a priori assumption \eqref{eq-apriori-Lq}
with $\omega_{[n/2]+k-2}=\frac{1+\lambda}{4}n$.
It suffices to set $\lambda\in(-\frac{n}{n+2},0)$.
$\hfill\Box$

\section{Time-weighted iteration scheme}
\label{sec-iteration}

In this section we develop a new technique which is the artful combination of the
time-weighted energy method and Green function method
to formulate the decay estimates of the over-damped Euler equation.
As shown in the above section, the Green function method is powerful in the optimal decay
estimates of the low-order energies but may have some troubles for the high-order energies.
Meanwhile, the classical weighted energy method is suitable for high-order energy estimates
but the decay rates are generally not optimal.
Therefore, we combine these two methods together.

Denote as before $b(t)=\frac{\mu}{(1+t)^\lambda}$ with $\mu>0$ and $\lambda\in[-1,0)$,
and
$$
Q_1(t,x)=-\boldsymbol u\cdot\nabla v-\varpi v\nabla\cdot \boldsymbol u,
\qquad
Q_2(t,x)=-(\boldsymbol u\cdot\nabla) \boldsymbol u-\varpi v\nabla v.
$$
We note that $Q_2$ is vector-valued, which should be written as $\boldsymbol Q_2$,
but we slightly abuse the notion of $Q_2$ for simplicity.
We may also write $\boldsymbol u$ as $u$ in the proof of this section for convenience.
Rewrite the nonlinear system into nonlinear wave equations
\begin{equation} \label{eq-wave-v}
\partial_t^2 v-\Delta v+b(t)\cdot\partial_t v=
\partial_tQ_1+b(t)\cdot Q_1-\nabla\cdot Q_2,
\end{equation}
and
\begin{equation} \label{eq-wave-u}
\partial_t^2 \boldsymbol u-\Delta \boldsymbol u+\partial_t (b(t)\cdot\boldsymbol u)=
\partial_t Q_2-\nabla Q_1.
\end{equation}

\subsection{Time-weighted energy estimates}

The main idea of the time-weighted iteration scheme is to
sacrifice the decay estimates of the low-order energies
(i.e., $\|\partial_t^j\partial_x^k(v,\boldsymbol u)\|$ with $j=0,1$ and $k+j=m\ge0$)
for better decay rates of high-order energies
(i.e., $\|\partial_t^j\partial_x^k(v,\boldsymbol u)\|$ with $j=0,1$ and $k+j=m+1$)
in the time-weighted energy estimates,
and the optimal decay rates of the basic energy $\|(v,\boldsymbol u)\|$
are closed through the Green function method,
where those better decays of high-order energies are necessary.

We have the following time-weighted energy estimates for $\lambda\in(-1,0)$
(the critical case of $\lambda=-1$ will be treated separately in next section).
Note that we replace the small negative constant in the classical time-weighted energy method
by a small positive constant $\delta$,
such that the high-order energies are decaying better
but the estimates on the low-order energies are absent.

\begin{lemma} \label{le-timew}
For any nonnegative integer $k$, $\lambda\in(-1,0)$,
$\delta\in(0,\frac{1+\lambda}{4})$, and $|\alpha|=k$, there hold
\begin{align} \nonumber
&\frac{d}{dt}\int E^v(\partial_t \partial_x^\alpha v,\nabla \partial_x^\alpha v,\partial_x^\alpha v)
+\int \big[
(1+t)^{1+\delta}|\partial_t \partial_x^\alpha v|^2
+(1+t)^{\lambda+\delta}|\nabla \partial_x^\alpha v|^2
\big]
\\ \label{eq-energy-v}
&\lesssim
\int (1+t)^{\delta-1}(\partial_x^\alpha v)^2
+\int \partial_x^\alpha(\partial_tQ_1+b(t)\cdot Q_1-\nabla\cdot Q_2)
\cdot((1+t)^{1+\lambda+\delta}\partial_t \partial_x^\alpha v
+\mu_1(1+t)^{\lambda+\delta} \partial_x^\alpha v),
\end{align}
and
\begin{align} \nonumber
&\frac{d}{dt}\int E^u(\partial_t \partial_x^\alpha\boldsymbol u,
\nabla \partial_x^\alpha\boldsymbol u,\partial_x^\alpha\boldsymbol u)
+\int \big[
(1+t)^{1-2\lambda+\delta}|\partial_t \partial_x^\alpha\boldsymbol u|^2
+(1+t)^{-\lambda+\delta}|\nabla \partial_x^\alpha\boldsymbol u|^2
\big]
\\ \label{eq-energy-u}
&\lesssim
\int (1+t)^{-2\lambda+\delta-1}
|\partial_x^\alpha\boldsymbol u|^2
+\int \partial_x^\alpha(\partial_t Q_2-\nabla Q_1)
\cdot((1+t)^{1-\lambda+\delta}\partial_t \partial_x^\alpha\boldsymbol u
+\mu_2(1+t)^{-\lambda+\delta} \partial_x^\alpha\boldsymbol u),
\end{align}
where $\mu_1>0$ and $\mu_2>0$ are constants and
\begin{align*}
E^v(\partial_t \partial_x^\alpha v,\nabla \partial_x^\alpha v,\partial_x^\alpha v)
&\approx
(1+t)^{1+\lambda+\delta}(|\partial_t \partial_x^\alpha v|^2+|\nabla \partial_x^\alpha v|^2)
+(1+t)^{\delta}(\partial_x^\alpha v)^2,
\\
E^u(\partial_t \partial_x^\alpha \boldsymbol u,\nabla \partial_x^\alpha \boldsymbol u,
\partial_x^\alpha \boldsymbol u)
&\approx
(1+t)^{1-\lambda+\delta}(|\partial_t \partial_x^\alpha \boldsymbol u|^2
+|\nabla \partial_x^\alpha \boldsymbol u|^2)
+(1+t)^{-2\lambda+\delta}|\partial_x^\alpha \boldsymbol u|^2.
\end{align*}
\end{lemma}
{\it\bfseries Proof.}
Multiplying \eqref{eq-wave-v} by
$(1+t)^{1+\lambda+\delta}\partial_t v+\mu_1(1+t)^{\lambda+\delta} v$
with $\delta\in(0,\frac{1+\lambda}{4})$ and $\mu_1>0$, we have
(similar to Proposition A.1 in Appendix A of \cite{Todoraova})
\begin{align*}
&\frac{d}{dt}\int \big[(1+t)^{1+\lambda+\delta}(|\partial_t v|^2+|\nabla v|^2)
+2\mu_1(1+t)^{\lambda+\delta} v\partial_tv
+(\mu_1 b(t)(1+t)^{\lambda+\delta}-(\lambda+\delta)\mu_1(1+t)^{\lambda+\delta-1})v^2\big]
\\
&\quad+\int \big[
(-(1+\lambda+\delta)(1+t)^{\lambda+\delta}+2b(t)(1+t)^{1+\lambda+\delta}
-2\mu_1(1+t)^{\lambda+\delta})|\partial_t v|^2
\\
&\qquad\qquad
+(-(1+\lambda+\delta)(1+t)^{\lambda+\delta}+2\mu_1(1+t)^{\lambda+\delta})|\nabla v|^2
\big]
\\
&\quad+\int ((\lambda+\delta)(\lambda+\delta-1)\mu_1(1+t)^{\lambda+\delta-2}
-\partial_t(\mu_1b(t)(1+t)^{\lambda+\delta}))v^2
\\
&=2\int (\partial_tQ_1+b(t)\cdot Q_1-\nabla\cdot Q_2)
\cdot((1+t)^{1+\lambda+\delta}\partial_t v+\mu_1(1+t)^{\lambda+\delta} v),
\end{align*}
which can be simplified as
\begin{align*}
&\frac{d}{dt}\int E^v(\partial_t v,\nabla v,v)
+\int \big[
(1+t)^{1+\delta}|\partial_t v|^2
+(2\mu_1-(1+\lambda+\delta))(1+t)^{\lambda+\delta}|\nabla v|^2
\big]
\\
&\lesssim
\int (1+t)^{\delta-1}v^2
+\int (\partial_tQ_1+b(t)\cdot Q_1-\nabla\cdot Q_2)
\cdot((1+t)^{1+\lambda+\delta}\partial_t v+\mu_1(1+t)^{\lambda+\delta} v),
\end{align*}
where
\begin{align*}
&E^v(\partial_t v,\nabla v,v)
\\
&:=(1+t)^{1+\lambda+\delta}(|\partial_t v|^2+|\nabla v|^2)
+2\mu_1(1+t)^{\lambda+\delta} v\partial_tv
+(\mu_1 b(t)(1+t)^{\lambda+\delta}-(\lambda+\delta)\mu_1(1+t)^{\lambda+\delta-1})v^2
\\
&\approx
(1+t)^{1+\lambda+\delta}(|\partial_t v|^2+|\nabla v|^2)
+(1+t)^{\delta}v^2.
\end{align*}
Here we fix $\mu_1$ such that $\mu_1\ge 1+\lambda+\delta>0$.

Next, multiplying \eqref{eq-wave-u} by
$(1+t)^{1-\lambda+\delta}\partial_t \boldsymbol u+\mu_2(1+t)^{-\lambda+\delta} \boldsymbol u$
with $\delta\in(0,\frac{1+\lambda}{4})$ and $\mu_2>0$, we have
\begin{align*}
&\frac{d}{dt}\int \big[
(1+t)^{1-\lambda+\delta}(|\partial_t \boldsymbol u|^2+|\nabla \boldsymbol u|^2)
+2\mu_2(1+t)^{-\lambda+\delta} \boldsymbol u\cdot\partial_t \boldsymbol u
\\
&\qquad\qquad
+(\mu_2 b(t)(1+t)^{-\lambda+\delta}-(-\lambda+\delta)\mu_2(1+t)^{-\lambda+\delta-1}
+b'(t)\cdot (1+t)^{1-\lambda+\delta})|\boldsymbol u|^2\big]
\\
&\quad+\int \big[
(-(1-\lambda+\delta)(1+t)^{-\lambda+\delta}+2b(t)(1+t)^{1-\lambda+\delta}
-2\mu_2(1+t)^{-\lambda+\delta})|\partial_t \boldsymbol u|^2
\\
&\qquad\qquad+(-(1-\lambda+\delta)(1+t)^{-\lambda+\delta}
+2\mu_2(1+t)^{-\lambda+\delta})|\nabla \boldsymbol u|^2
\big]
\\
&\quad+\int ((-\lambda+\delta)(-\lambda+\delta-1)\mu_2(1+t)^{-\lambda+\delta-2}
-\partial_t(\mu_2b(t)(1+t)^{-\lambda+\delta})
\\
&\qquad\qquad+2\mu_2 b'(t)(1+t)^{-\lambda+\delta}
-\partial_t(b'(t)(1+t)^{1-\lambda+\delta})
)|\boldsymbol u|^2
\\
&=2\int (\partial_t Q_2-\nabla Q_1)
\cdot((1+t)^{1-\lambda+\delta}\partial_t \boldsymbol u
+\mu_2(1+t)^{-\lambda+\delta} \boldsymbol u).
\end{align*}
We simplify the above equality as
\begin{align*}
&\frac{d}{dt}\int E^u(\partial_t \boldsymbol u,\nabla \boldsymbol u,\boldsymbol u)
+\int \big[
(1+t)^{1-2\lambda+\delta}|\partial_t \boldsymbol u|^2
+(2\mu_2-(1-\lambda+\delta))(1+t)^{-\lambda+\delta}|\nabla \boldsymbol u|^2
\big]
\\
&\lesssim
\int
((\mu_2-\lambda)\delta+2\lambda^2)(1+t)^{-2\lambda+\delta-1}
|\boldsymbol u|^2
+\int (\partial_t Q_2-\nabla Q_1)
\cdot((1+t)^{1-\lambda+\delta}\partial_t \boldsymbol u
+\mu_2(1+t)^{-\lambda+\delta} \boldsymbol u),
\end{align*}
where
\begin{align*}
E^u(\partial_t \boldsymbol u,\nabla \boldsymbol u,\boldsymbol u)
:=&
(1+t)^{1-\lambda+\delta}(|\partial_t \boldsymbol u|^2+|\nabla \boldsymbol u|^2)
+2\mu_2(1+t)^{-\lambda+\delta} \boldsymbol u\cdot\partial_t \boldsymbol u
\\
&+(\mu_2 b(t)(1+t)^{-\lambda+\delta}-(-\lambda+\delta)\mu_2(1+t)^{-\lambda+\delta-1}
+b'(t)\cdot (1+t)^{1-\lambda+\delta})|\boldsymbol u|^2
\\
\approx&
(1+t)^{1-\lambda+\delta}(|\partial_t \boldsymbol u|^2+|\nabla \boldsymbol u|^2)
+(1+t)^{-2\lambda+\delta}|\boldsymbol u|^2.
\end{align*}
We choose $\mu_2>0$ such that $\mu_2\ge 1-\lambda+\delta$.
Thus, the proof for the case of $k=0$ is completed.

Differentiating  $\partial_x^\alpha \eqref{eq-wave-v}$ and $\partial_x^\alpha \eqref{eq-wave-u}$,
and multiplying the resulting equations by
$(1+t)^{1+\lambda+\delta}\partial_t \partial_x^\alpha v
+\mu_1(1+t)^{\lambda+\delta} \partial_x^\alpha v$
and $(1+t)^{1-\lambda+\delta}\partial_t \partial_x^\alpha\boldsymbol u
+\mu_2(1+t)^{-\lambda+\delta} \partial_x^\alpha\boldsymbol u$
respectively,
we can prove \eqref{eq-energy-v} and \eqref{eq-energy-u}
in a similar procedure. The detail is omitted.
$\hfill\Box$

\begin{remark}
Compared with the multiplier method developed by Todorova and Yordanov \cite{Todoraova}
for the wave equation with variable coefficients
($b(t)=\frac{\mu}{(1+t)^\lambda}$ replaced by $\frac{\mu}{(1+|x|)^\alpha}$ with $\alpha\in(0,1)$)
and the weighted energy method employed by Pan \cite{PanXH-AA} for the
wave equation with under-damping with $\lambda\in(0,1)$,
here for over-damping with $\lambda\in[-1,0)$ we take the weights only dependent on time.
The reason is that for the over-damping case,
the simple weights depending on time can take advantage of the time-asymptotically growing
over-damping, which turns out to be sufficient for the closure of the decay estimates
for all $\lambda\in(-1,0)$.
\end{remark}

\begin{remark}
The energy estimates \eqref{eq-energy-v} and \eqref{eq-energy-u} are deduced by rewriting
both $v$ and $\boldsymbol u$ as solutions to time-dependent damped nonlinear wave equations.
This differs from the approach in \cite{PanXH-AA} for under-damping case,
where the estimates of $\boldsymbol u$ are formulated
according to the equation \eqref{eq-vbdu}$_2$.
Here for the over-damping case we cannot apply the above procedure in \cite{PanXH-AA}
since the estimates on $\|\partial_x^k\boldsymbol u\|$ depends on
at least one of $\|\partial_x^{k+1}v\|$ and  $\|\partial_x^{k+1}\boldsymbol u\|$,
and other efforts should be made for the closure of the weighted energy estimates.
\end{remark}

We define the following time-weighted energy functions
for $N\ge[\frac{n}{2}]+2$ and $0\le k\le N-1$,
\begin{align} \label{eq-energy-Phi}
\Phi_{k+1}(T):=\sup_{t\in(0,T)}\Big\{
\sum_{|\alpha|=k}\Big[ &
(1+t)^{1+\lambda+\delta}\int(|\partial_t \partial_x^\alpha v|^2+|\nabla \partial_x^\alpha v|^2) \\
&+(1+t)^{1-\lambda+\delta}\int(|\partial_t \partial_x^\alpha\boldsymbol u|^2
+|\nabla \partial_x^\alpha\boldsymbol u|^2)
\Big]
\Big\}^\frac{1}{2},
\end{align}
and
\begin{align} \nonumber
\Psi_{k+1}(T):=\sup_{t\in(0,T)}\Big\{
\sum_{|\alpha|=k}\Big[ & \int \Big(
(1+t)^{1+\delta}|\partial_t \partial_x^\alpha v|^2
+(1+t)^{\lambda+\delta}|\nabla \partial_x^\alpha v|^2
\Big)
\\ \label{eq-energy-Psi}
& +\int \big[
(1+t)^{1-2\lambda+\delta}|\partial_t \partial_x^\alpha\boldsymbol u|^2
+(1+t)^{-\lambda+\delta}|\nabla \partial_x^\alpha\boldsymbol u|^2
\big]
\Big]
\Big\}^\frac{1}{2},
\end{align}
which satisfies $\Psi_{k+1}^2(t)\ge (1+t)^{-1}\cdot\Phi_{k+1}^2(t)$.
We may assume that $\Phi_{k+1}(T)\ge \Phi_k(T)$ for all $k\ge1$ and $T$.
Otherwise, we can modify the definition of $\Phi_{k+1}(T)$.
The energy function $\Phi_{k+1}(T)$ is defined according to the time-weighted energy estimates
in Lemma \ref{le-timew},
but the decay estimates on $\|v\|$ and $\|\boldsymbol u\|$ are absent
and insufficient for the closure of the energy estimates.
Additionally, we define the following weighted energy function
\begin{align} \label{eq-energy-Psi0}
\Psi_0(T):=\sup_{t\in(0,T)}\Big\{&
(1+t)^{\frac{1+\lambda}{4}n}\|v\|,
(1+t)^{\frac{1+\lambda}{4}n+\frac{1-\lambda}{2}}
\|\boldsymbol u\|
\Big\}.
\end{align}
The energy estimates in $\Psi_0(T)$ will be closed through the Green function method
instead of the time-weighted energy method.
There holds
\begin{equation} \label{eq-initial}
\|(v_0,\boldsymbol u_0)\|_{H^N}
\approx \sum_{k=1}^{N}\Phi_k(0)+\Psi_0(0)
\approx \Phi_N(0)+\Psi_0(0).
\end{equation}
According to Sobolev embedding theorem, we have
\begin{align} \label{eq-Linfty}
(1+t)^\frac{1+\lambda+\delta}{2}\|\partial_x^j v\|_{L^\infty}
+(1+t)^\frac{1-\lambda+\delta}{2}\|\partial_x^j \boldsymbol u\|_{L^\infty}
\lesssim \max_{1\le k\le [\frac{n}{2}]+2}\Phi_k(t)
\lesssim \Phi_N(t),
\quad 0\le j\le 1, n\ge3,
\end{align}
and
\begin{align} \nonumber
&(1+t)^{\frac{1+\lambda}{2}+\frac{\delta}{4}}\|v\|_{L^\infty}
+(1+t)^{\frac{1-\lambda}{2}+\frac{1+\lambda+\delta}{4}}\|\boldsymbol u\|_{L^\infty}
+(1+t)^\frac{1+\lambda+\delta}{2}\|\partial_x v\|_{L^\infty}
+(1+t)^\frac{1-\lambda+\delta}{2}\|\partial_x \boldsymbol u\|_{L^\infty}
\\ \label{eq-Linftyn2}
&\lesssim \max_{1\le k\le [\frac{n}{2}]+2}\Phi_k(t)+\Psi_0(t)
\lesssim \Phi_N(t)+\Psi_0(t),
\quad n=2.
\end{align}

We have the following iteration scheme based on Lemma \ref{le-timew}.

\begin{lemma}[Time-weighted iteration scheme] \label{le-iteration}
For $\lambda\in(-1,0)$ and $\delta\in(0,\frac{1+\lambda}{4})$, there holds
\begin{align} \nonumber
\Phi_1^2(t)+\int_0^t \Psi_1^2(s)ds
\lesssim&
\Phi_1^2(0)+
\int_0^t (1+s)^{-1-\frac{1+\lambda}{2}n+\delta}\cdot\Psi_0^2(s)ds
\\ \nonumber
&+\int_0^t\int (\partial_tQ_1+b(s)\cdot Q_1-\nabla\cdot Q_2)
\cdot((1+s)^{1+\lambda+\delta}\partial_t v
+\mu_1(1+s)^{\lambda+\delta} v)ds
\\ \label{eq-it-0}
&+\int_0^t\int (\partial_t Q_2-\nabla Q_1)
\cdot((1+s)^{1-\lambda+\delta}\partial_t \boldsymbol u
+\mu_2(1+s)^{-\lambda+\delta} \boldsymbol u)ds.
\end{align}
For any integer $k\ge1$, there holds
\begin{align} \nonumber
\Phi_{k+1}^2(t)&+\int_0^t \Psi_{k+1}^2(s)ds
\lesssim
\Phi_{k+1}^2(0)+
\int_0^t (1+s)^{-1-\lambda}\cdot\Psi_k^2(s)ds
\\ \nonumber
&+\sum_{|\alpha|=k}
\int_0^t\int \partial_x^\alpha(\partial_tQ_1+b(s)\cdot Q_1-\nabla\cdot Q_2)
\cdot((1+s)^{1+\lambda+\delta}\partial_t \partial_x^\alpha v
+\mu_1(1+s)^{\lambda+\delta} \partial_x^\alpha v)ds
\\ \label{eq-it-1}
&+\sum_{|\alpha|=k}
\int_0^t\int \partial_x^\alpha(\partial_t Q_2-\nabla Q_1)
\cdot((1+s)^{1-\lambda+\delta}\partial_t \partial_x^\alpha\boldsymbol u
+\mu_2(1+s)^{-\lambda+\delta} \partial_x^\alpha\boldsymbol u)ds.
\end{align}
\end{lemma}
{\it\bfseries Proof.}
This is a simple conclusion of Lemma \ref{le-timew}
with the notations $\Phi_k(t)$, $\Psi_k(t)$, and $\Psi_0(t)$ defined by \eqref{eq-energy-Phi},
\eqref{eq-energy-Psi}, and \eqref{eq-energy-Psi0}.
$\hfill\Box$

\subsection{A priori estimates involving inhomogeneous terms}

We estimate the inhomogeneous terms in the inequalities \eqref{eq-energy-v}
and \eqref{eq-energy-u} in Lemma \ref{le-timew}.
We first consider the case of $k=0$
and in order to extend the proof to a general case of $k>0$
we should avoid directly using the energy estimates of
the second order derivatives (such as $\|\partial_t\nabla v\|$) in $\Phi_k(t)$,
since that would be $(k+2)$-th order derivatives for general $k>0$ and
cause trouble in the closure of the weighted energy estimates.

\begin{lemma} \label{le-Q}
There holds, for $\lambda\in(-1,0)$ and $\delta\in(0,\frac{1+\lambda}{4})$, that
\begin{align*}
&\int (\partial_tQ_1+b(t)\cdot Q_1-\nabla\cdot Q_2)
\cdot((1+t)^{1+\lambda+\delta}\partial_t v+\mu_1(1+t)^{\lambda+\delta} v)
\\
&\qquad+\int (\partial_t Q_2-\nabla Q_1)
\cdot((1+t)^{1-\lambda+\delta}\partial_t \boldsymbol u
+\mu_2(1+t)^{-\lambda+\delta} \boldsymbol u)
\\
&\lesssim
\partial_tJ_1(t)+
(\Psi_0(t)+\Phi_N(t))\cdot\Psi_1^2(t)
+\Phi_N(t)\cdot\Psi_0^2(t)\cdot(1+t)^{-1-\frac{1+\lambda}{4}},
\end{align*}
provided that $\|v\|_{L^\infty}\le \frac{1}{\gamma-1}$
(which is valid under the a priori assumption $\Phi_N(t)+\Psi_0(t)\le \delta_0$
with a small constant $\delta_0$),
where
$$
J_1(t)\lesssim \|v\|_{L^\infty}\cdot\Phi_1^2(t).
$$
\end{lemma}
{\it\bfseries Proof.}
The estimates of the two integrals are separated into two steps.

{\it Step I}. We first estimate the term involving $b(t)\cdot Q_1$ as follows
\begin{align*}
&\int b(t)\cdot Q_1
\cdot((1+t)^{1+\lambda+\delta}\partial_t v+\mu_1(1+t)^{\lambda+\delta} v)
\\
&\lesssim \int (|\boldsymbol u\cdot\nabla v|+|v\nabla\cdot \boldsymbol u|)
\cdot ((1+t)^{1+\delta}|\partial_t v|+(1+t)^{\delta} |v|)
\\
&\lesssim
\int (\|\nabla \boldsymbol u\|_{L^\infty}\cdot|v|
+\|\nabla v\|_{L^\infty}\cdot|\boldsymbol u|)
\cdot((1+t)^{1+\delta}|\partial_t v|+(1+t)^{\delta} |v|)
\\
&\lesssim \Phi_N(t)(1+t)^{-\frac{1-\lambda+\delta}{2}}(1+t)^{1+\delta}
\int |v||\partial_t v|
+\Phi_N(t)(1+t)^{-\frac{1-\lambda+\delta}{2}}(1+t)^{\delta}\int v^2
\\
&\quad +\Phi_N(t)(1+t)^{-\frac{1+\lambda+\delta}{2}}(1+t)^{1+\delta}
\int |\boldsymbol u||\partial_t v|
+\Phi_N(t)(1+t)^{-\frac{1+\lambda+\delta}{2}}(1+t)^{\delta}\int |\boldsymbol u||v|
\\
&\lesssim
\Phi_N(t)(1+t)^{-\frac{1-\lambda+\delta}{2}}(1+t)^{1+\delta}
\cdot\Psi_1(t)(1+t)^{-\frac{1+\delta}{2}}
\cdot\Psi_0(t)(1+t)^{-\frac{1+\lambda}{4}n}
\\
&\quad +\Phi_N(t)(1+t)^{-\frac{1-\lambda+\delta}{2}}(1+t)^{\delta}
\cdot\Psi_0^2(t)(1+t)^{-\frac{1+\lambda}{2}n}
\\
&\quad +\Phi_N(t)(1+t)^{-\frac{1+\lambda+\delta}{2}}(1+t)^{1+\delta}
\cdot\Psi_1(t)(1+t)^{-\frac{1+\delta}{2}}
\cdot\Psi_0(t)(1+t)^{-\frac{1+\lambda}{4}n-\frac{1-\lambda}{2}}
\\
&\quad +\Phi_N(t)(1+t)^{-\frac{1+\lambda+\delta}{2}}(1+t)^{\delta}
\cdot \Psi_0(t)(1+t)^{-\frac{1+\lambda}{4}n-\frac{1-\lambda}{2}}
\cdot\Psi_0(t)(1+t)^{-\frac{1+\lambda}{4}n}
\\
&\lesssim
\Phi_N(t)\Psi_1^2(t)
+\Phi_N(t)\Psi_0^2(t)(1+t)^{\lambda-\frac{1+\lambda}{2}n}
+\Phi_N(t)\Psi_0^2(t)(1+t)^{-\frac{1-\lambda}{2}-\frac{1+\lambda}{2}n+\frac{\delta}{2}}
\\
&\quad
+\Phi_N(t)\Psi_0^2(t)(1+t)^{-1-\frac{1+\lambda}{2}n}
+\Phi_N(t)\Psi_0^2(t)(1+t)^{-1-\frac{1+\lambda}{2}n+\frac{\delta}{2}},
\end{align*}
where
\begin{align*}
\begin{cases}
\frac{1-\lambda}{2}+\frac{1+\lambda}{2}n-\frac{\delta}{2}
\ge \frac{1-\lambda}{2}+\frac{1+\lambda}{2}+\frac{1+\lambda}{2}-\frac{\delta}{2}
\ge 1+\frac{1+\lambda}{4},
\\
1+\frac{1+\lambda}{2}n\ge 1+\frac{1+\lambda}{4},
\\
1+\frac{1+\lambda}{2}n-\frac{\delta}{2}
\ge 1+\frac{1+\lambda}{4},
\end{cases}
\end{align*}
for all $n\ge2$ and $\lambda\in(-1,0)$,
and
$$-\lambda+\frac{1+\lambda}{2}n
=1+\frac{1+\lambda}{2}(n-2)>1+\frac{1+\lambda}{4}$$
for $n\ge3$.
For the case of $n=2$, we modify the above estimate
(replacing the inequality \eqref{eq-Linfty} by \eqref{eq-Linftyn2}) as
\begin{align*}
&\int |v\nabla\cdot \boldsymbol u|\cdot (1+t)^{1+\delta}|\partial_t v|
\\
&\lesssim
\int \|v\|_{L^\infty}|\nabla\cdot \boldsymbol u|\cdot (1+t)^{1+\delta}|\partial_t v|
\\
&\lesssim
(\Psi_0(t)+\Phi_N(t))(1+t)^{-\frac{1+\lambda}{2}-\frac{\delta}{4}}
\cdot
(1+t)^{1+\delta}
\cdot
\Psi_1(t)(1+t)^{-\frac{1+\delta}{2}}
\cdot
\Psi_1(t)(1+t)^{-\frac{-\lambda+\delta}{2}}
\\
&=(\Psi_0(t)+\Phi_N(t))\Psi_1^2(t)(1+t)^{-\frac{\delta}{4}}.
\end{align*}

Next, we calculate the term involving $\partial_tQ_1$ as follows
\begin{align*}
&\int \partial_tQ_1
\cdot((1+t)^{1+\lambda+\delta}\partial_t v+\mu_1(1+t)^{\lambda+\delta} v)
\\
&= \int
(-\partial_t\boldsymbol u\cdot\nabla v-\varpi\partial_t v\nabla\cdot \boldsymbol u)
\cdot((1+t)^{1+\lambda+\delta}\partial_t v+\mu_1(1+t)^{\lambda+\delta} v)
\\
&\quad+\int
(-\boldsymbol u\cdot\nabla \partial_t v)
\cdot((1+t)^{1+\lambda+\delta}\partial_t v+\mu_1(1+t)^{\lambda+\delta} v)
\\
&\quad+\int
(-\varpi v\nabla\cdot \partial_t\boldsymbol u)
\cdot((1+t)^{1+\lambda+\delta}\partial_t v+\mu_1(1+t)^{\lambda+\delta} v)
\\
&=:I_{11}+I_{12}+I_{13}.
\end{align*}
We have
\begin{align*}
I_{11}&\lesssim
(\|\partial_t\boldsymbol u\|\cdot\|\nabla v\|_{L^\infty}
+\|\partial_t v\|\cdot\|\nabla\cdot \boldsymbol u\|_{L^\infty})
\cdot((1+t)^{1+\lambda+\delta}\|\partial_t v\|+(1+t)^{\lambda+\delta} \|v\|)
\\
&\lesssim
(\Psi_1(t)(1+t)^{-\frac{1-2\lambda+\delta}{2}}\cdot
\Phi_N(t)(1+t)^{-\frac{1+\lambda+\delta}{2}}
+\Psi_1(t)(1+t)^{-\frac{1+\delta}{2}}\cdot
\Phi_N(t)(1+t)^{-\frac{1-\lambda+\delta}{2}})
\\
&\qquad\cdot
((1+t)^{1+\lambda+\delta}\Psi_1(t)(1+t)^{-\frac{1+\delta}{2}}
+(1+t)^{\lambda+\delta} \Psi_0(t)(1+t)^{-\frac{1+\lambda}{4}n})
\\
&\lesssim
\Phi_N(t)\Psi_1^2(t)(1+t)^{-\frac{1+\delta}{2}+\frac{3}{2}\lambda}
+\Phi_N(t)\Psi_1^2(t)(1+t)^{-1-\frac{1+\lambda}{4}n+\frac{3}{2}\lambda}
+\Phi_N(t)\Psi_0^2(t)(1+t)^{-1-\frac{1+\lambda}{4}n+\frac{3}{2}\lambda}.
\end{align*}
The crucial point in the estimates of $I_{12}$ and $I_{13}$
is to avoid the direct estimates on $\nabla \partial_t v$
and $\nabla\cdot \partial_t\boldsymbol u$ through integration by parts such that
\begin{align*}
I_{12}&=
-\frac{1}{2}(1+t)^{1+\lambda+\delta}\int
\boldsymbol u\cdot\nabla (\partial_t v)^2
-\mu_1(1+t)^{\lambda+\delta}
\int
v\boldsymbol u\cdot\nabla \partial_t v
\\
&=\frac{1}{2}(1+t)^{1+\lambda+\delta}\int
(\nabla\cdot\boldsymbol u)(\partial_t v)^2
+\mu_1(1+t)^{\lambda+\delta}
\int
\partial_t v\cdot(\boldsymbol u\cdot\nabla v+v\nabla\cdot \boldsymbol u)
\\
&\lesssim
(1+t)^{1+\lambda+\delta}
\|\nabla\cdot\boldsymbol u\|_{L^\infty}\|\partial_t v\|^2
+(1+t)^{\lambda+\delta}
\|\partial_t v\|\cdot(\|\boldsymbol u\|\cdot\|\nabla v\|_{L^\infty}
+\|v\|\cdot\|\nabla\cdot \boldsymbol u\|_{L^\infty})
\\
&\lesssim
(1+t)^{1+\lambda+\delta}\cdot\Phi_N(t)(1+t)^{-\frac{1-\lambda+\delta}{2}}
\cdot \Psi_1^2(t)(1+t)^{-(1+\delta)}
\\
&\quad
+(1+t)^{\lambda+\delta}\cdot
\Psi_1(t)(1+t)^{-\frac{1+\delta}{2}}
\\
&\quad\quad
\cdot\big(
\Phi_0(t)(1+t)^{-\frac{1+\lambda}{4}n-\frac{1-\lambda}{2}}\cdot
\Phi_N(t)(1+t)^{-\frac{1+\lambda+\delta}{2}}
+\Phi_0(t)(1+t)^{-\frac{1+\lambda}{4}n}\cdot
\Phi_N(t)(1+t)^{-\frac{1-\lambda+\delta}{2}}
\big)
\\
&\lesssim
\Phi_N(t)\Psi_1^2(t)\cdot(1+t)^{-\frac{1+\lambda+\delta}{2}}
+\Phi_N(t)\Psi_1^2(t)(1+t)^{-1-\frac{1+\lambda}{4}n+\frac{3}{2}\lambda}
+\Phi_N(t)\Psi_0^2(t)(1+t)^{-1-\frac{1+\lambda}{4}n+\frac{3}{2}\lambda},
\end{align*}
and
\begin{align*}
I_{13}=
-(1+t)^{1+\lambda+\delta}\int
\varpi v\partial_t v\nabla\cdot \partial_t\boldsymbol u
-\mu_1(1+t)^{\lambda+\delta}\int
\varpi v^2\nabla\cdot \partial_t\boldsymbol u
=:I_{13}^1+I_{13}^2,
\end{align*}
where
\begin{align*}
I_{13}^2&=
\mu_1(1+t)^{\lambda+\delta}\int
\varpi\nabla(v^2)\cdot \partial_t\boldsymbol u
\\
&\lesssim
(1+t)^{\lambda+\delta}
\|v\|\cdot\|\nabla v\|_{L^\infty}\cdot\|\partial_t\boldsymbol u\|
\\
&\lesssim
(1+t)^{\lambda+\delta}\cdot
\Psi_0(t)(1+t)^{-\frac{1+\lambda}{4}n}
\cdot\Phi_N(t)(1+t)^{-\frac{1+\lambda+\delta}{2}}
\cdot\Psi_1(t)(1+t)^{-\frac{1-2\lambda+\delta}{2}}
\\
&\lesssim
\Phi_N(t)\Psi_1^2(t)(1+t)^{-1-\frac{1+\lambda}{4}n+\frac{3}{2}\lambda}
+\Phi_N(t)\Psi_0^2(t)(1+t)^{-1-\frac{1+\lambda}{4}n+\frac{3}{2}\lambda}.
\end{align*}
The treatment of $\nabla\cdot \partial_t \boldsymbol u$ in $I_{13}^1$
is to rewrite \eqref{eq-vbdu}$_1$ into
\begin{equation} \label{eq-znablau}
\nabla\cdot\boldsymbol u=-\frac{\partial_t v+\boldsymbol u\cdot\nabla v}{1+\varpi v},
\end{equation}
with $1+\varpi v\ge 1/2$ since $\|v\|_{L^\infty}\le1/(\gamma-1)$ and then
$$
\nabla\cdot\partial_t\boldsymbol u
=-\frac{\partial_t^2 v+\partial_t\boldsymbol u\cdot\nabla v
+\boldsymbol u\cdot\nabla \partial_t v}{1+\varpi v}
+\frac{\varpi(\partial_t v+\boldsymbol u\cdot\nabla v)\partial_t v}{(1+\varpi v)^2},
$$
similar to the proof of Lemma 2.2 in \cite{PanXH-AA}
but the most tricky parts and details are different.
Therefore,
\begin{align*}
I_{13}^1&=-(1+t)^{1+\lambda+\delta}\int
\varpi v\partial_t v\nabla\cdot \partial_t\boldsymbol u
\\
&=(1+t)^{1+\lambda+\delta}\int
\varpi v\partial_t v \cdot
\frac{\partial_t^2 v+\partial_t\boldsymbol u\cdot\nabla v
+\boldsymbol u\cdot\nabla \partial_t v}{1+\varpi v}
\\
&\quad
-(1+t)^{1+\lambda+\delta}\int
\varpi v\partial_t v \cdot
\frac{\varpi(\partial_t v+\boldsymbol u\cdot\nabla v)\partial_t v}{(1+\varpi v)^2}
\\
&=:I_{13}^{11}+I_{13}^{12},
\end{align*}
where
\begin{align*}
I_{13}^{12}
&\lesssim
(1+t)^{1+\lambda+\delta}
\|v\|_{L^\infty}\|\partial_t v\| \cdot
(\|\partial_t v\|_{L^\infty}+\|\boldsymbol u\|_{L^\infty}\|\nabla v\|_{L^\infty})\|\partial_t v\|
\\
&\lesssim
(1+t)^{1+\lambda+\delta}
\cdot\Phi_N(t)(1+t)^{-\frac{1+\lambda+\delta}{2}}
\cdot\Psi_1^2(t)(1+t)^{-(1+\delta)}
\\
&\lesssim
\Phi_N(t)\Psi_1^2(t)\cdot(1+t)^{-\frac{1-\lambda+\delta}{2}}.
\end{align*}
The estimate on $I_{13}^{11}$ is
\begin{align*}
I_{13}^{11}&=
(1+t)^{1+\lambda+\delta}\int
\varpi v\partial_t v \cdot
\frac{\partial_t^2 v}{1+\varpi v}
+
(1+t)^{1+\lambda+\delta}\int
\varpi v\partial_t v \cdot
\frac{\partial_t\boldsymbol u\cdot\nabla v}{1+\varpi v}
\\
&\qquad+
(1+t)^{1+\lambda+\delta}\int
\varpi v\partial_t v \cdot
\frac{\boldsymbol u\cdot\nabla \partial_t v}{1+\varpi v}
\\
&=
(1+t)^{1+\lambda+\delta}\frac{1}{2}\int
\varpi v \cdot
\frac{\partial_t(\partial_t v)^2}{1+\varpi v}
+
(1+t)^{1+\lambda+\delta}\int
\varpi v\partial_t v \cdot
\frac{\partial_t\boldsymbol u\cdot\nabla v}{1+\varpi v}
\\
&\qquad+
(1+t)^{1+\lambda+\delta}\frac{1}{2}\int
\varpi v \cdot
\frac{\boldsymbol u\cdot\nabla (\partial_t v)^2}{1+\varpi v}
\\
&=
\partial_t\Big(
(1+t)^{1+\lambda+\delta}\frac{1}{2}\int
\varpi v \cdot
\frac{(\partial_t v)^2}{1+\varpi v}\Big)
-\frac{1}{2}\int(\partial_t v)^2
\cdot\partial_t\Big(
\frac{\varpi v}{1+\varpi v}\cdot(1+t)^{1+\lambda+\delta}\Big)
\\
&\qquad+
(1+t)^{1+\lambda+\delta}\int
\varpi v\partial_t v \cdot
\frac{\partial_t\boldsymbol u\cdot\nabla v}{1+\varpi v}
-(1+t)^{1+\lambda+\delta}\frac{1}{2}\int(\partial_t v)^2\cdot
\Big(\nabla\cdot
\frac{\varpi v \boldsymbol u}{1+\varpi v}\Big)
\\
&=:\partial_tJ_1(t)+\tilde I_{13}^{11},
\end{align*}
where
\begin{equation} \label{eq-zJ1}
J_1(t):=(1+t)^{1+\lambda+\delta}\frac{1}{2}\int
\varpi v \cdot
\frac{(\partial_t v)^2}{1+\varpi v}
\lesssim (1+t)^{1+\lambda+\delta}\|v\|_{L^\infty}\|\partial_t v\|^2
\lesssim \|v\|_{L^\infty}\cdot\Phi_1^2(t).
\end{equation}
We see that $\tilde I_{13}^{11}$ are integrals only involving first order derivatives
and can be estimated in the similar way as $I_{11}$.
This completes the proof of the estimates involving $\partial_t Q_1$.

We now consider the term involving $-\nabla\cdot Q_2$ such that
\begin{align*}
&\int (-\nabla\cdot Q_2)
\cdot((1+t)^{1+\lambda+\delta}\partial_t v+\mu_1(1+t)^{\lambda+\delta} v)
\\
&=\int (
\sum_{j=1}^n\sum_{k=1}^n\partial_{x_k}u^j\cdot\partial_{x_j}u^k
+\varpi \nabla v\cdot\nabla v)
\cdot((1+t)^{1+\lambda+\delta}\partial_t v+\mu_1(1+t)^{\lambda+\delta} v)
\\
&\quad
+\int (\boldsymbol u\cdot\nabla)(\nabla\cdot \boldsymbol u)
\cdot((1+t)^{1+\lambda+\delta}\partial_t v+\mu_1(1+t)^{\lambda+\delta} v)
\\
&\quad
+\int \varpi v(\nabla\cdot\nabla v)
\cdot((1+t)^{1+\lambda+\delta}\partial_t v+\mu_1(1+t)^{\lambda+\delta} v)
\\
&=:I_{21}+I_{22}+I_{23}.
\end{align*}
Similar to $I_{11}$,
\begin{align*}
I_{21}
&\lesssim
(\|\nabla \boldsymbol u\|_{L^\infty}\|\nabla \boldsymbol u\|
+ \|\nabla v\|_{L^\infty}\|\nabla v\|)
\cdot((1+t)^{1+\lambda+\delta}\|\partial_t v\|+(1+t)^{\lambda+\delta} \|v\|)
\\
&\lesssim
\Phi_N(t)(1+t)^{-\frac{1+\lambda+\delta}{2}}
\cdot \Psi_1(t)(1+t)^{-\frac{\lambda+\delta}{2}}
\cdot ((1+t)^{1+\lambda+\delta}\Psi_1(t)(1+t)^{-\frac{1+\delta}{2}}
+(1+t)^{\lambda+\delta} \Psi_0(t)(1+t)^{-\frac{1+\lambda}{4}n})
\\
&\lesssim
\Phi_N(t)\Psi_1^2(t)\cdot(1+t)^{-\frac{\delta}{2}}
+\Phi_N(t)\Psi_1^2(t)
+\Phi_N(t)\Psi_0^2(t)\cdot(1+t)^{-1-\frac{1+\lambda}{2}n}.
\end{align*}
Integrating by parts implies that
\begin{align*}
I_{22}=&\int (\boldsymbol u\cdot\nabla)(\nabla\cdot \boldsymbol u)
\cdot((1+t)^{1+\lambda+\delta}\partial_t v+\mu_1(1+t)^{\lambda+\delta} v)
\\
=&-\int (\nabla\cdot \boldsymbol u)^2
\cdot((1+t)^{1+\lambda+\delta}\partial_t v+\mu_1(1+t)^{\lambda+\delta} v)
\\
&-\int (\nabla\cdot \boldsymbol u)\boldsymbol u
\cdot((1+t)^{1+\lambda+\delta}\nabla\partial_t v+\mu_1(1+t)^{\lambda+\delta}\nabla v),
\end{align*}
and
\begin{align*}
I_{23}&=\int \varpi v(\nabla\cdot\nabla v)
\cdot((1+t)^{1+\lambda+\delta}\partial_t v+\mu_1(1+t)^{\lambda+\delta} v)
\\
&=-\int \varpi |\nabla v|^2
\cdot(1+t)^{1+\lambda+\delta}\partial_t v
-\int \varpi v\nabla v\cdot
(1+t)^{1+\lambda+\delta}\nabla\partial_t v
-2\int \varpi |\nabla v|^2\cdot \mu_1(1+t)^{\lambda+\delta} v.
\end{align*}
All the above integrals not involving second order derivatives in $I_{22}$ and $I_{23}$
can be estimated as $I_{21}$,
except for
\begin{align*}
I_{22}^1:=&
-\int (\nabla\cdot \boldsymbol u)\boldsymbol u
\cdot(1+t)^{1+\lambda+\delta}\nabla\partial_t v,
\\
I_{23}^1:=&
-\int \varpi v\nabla v\cdot
(1+t)^{1+\lambda+\delta}\nabla\partial_t v,
\end{align*}
which need to be treated in the same procedure as $I_{13}^{1}$.
Specifically, we have according to \eqref{eq-znablau},
\begin{align*}
I_{22}^1:=&
-\int (\nabla\cdot \boldsymbol u)\boldsymbol u
\cdot(1+t)^{1+\lambda+\delta}\nabla\partial_t v
\\
=&
\int \frac{\partial_t v}{1+\varpi v}\boldsymbol u
\cdot(1+t)^{1+\lambda+\delta}\nabla\partial_t v
+\cdots
\\
=&
\frac{1}{2}\int \frac{\boldsymbol u}{1+\varpi v}
\cdot(1+t)^{1+\lambda+\delta}\nabla(\partial_t v)^2
+\cdots
\\
=&
-\frac{1}{2}\int (\nabla\cdot\frac{\boldsymbol u}{1+\varpi v})
\cdot(1+t)^{1+\lambda+\delta}(\partial_t v)^2
+\cdots
\end{align*}
where we only write down the cubic terms involving second order derivatives
and the integral in the last equality only involves first order derivatives.
According to \eqref{eq-wave-v} and integration by parts
\begin{align*}
I_{23}^1=
\int \varpi v\Delta v\cdot
(1+t)^{1+\lambda+\delta}\partial_t v+\cdots
=
(1+t)^{1+\lambda+\delta}\int \varpi v
\partial_t v\cdot\partial_t^2 v+\cdots
\end{align*}
whose most tricky part is the same as $I_{13}^{11}$ in $I_{13}^{1}$.
This completes the proof of the estimates involving $-\nabla\cdot Q_2$.

{\it Step II}.
We turn to show the estimates of the second integral of this lemma.
We may only focus on the terms involving second order derivatives
since the estimates on the others are similar to those in the first step of this proof.
We have
\begin{align*}
&\int (\partial_t Q_2-\nabla Q_1)
\cdot((1+t)^{1-\lambda+\delta}\partial_t \boldsymbol u
+\mu_2(1+t)^{-\lambda+\delta} \boldsymbol u)
\\
&=\int (-(\boldsymbol u\cdot\nabla)\partial_t \boldsymbol u)
\cdot((1+t)^{1-\lambda+\delta}\partial_t \boldsymbol u
+\mu_2(1+t)^{-\lambda+\delta} \boldsymbol u)
\\
&\quad+
\int (-\varpi v\nabla\partial_t v)
\cdot((1+t)^{1-\lambda+\delta}\partial_t \boldsymbol u
+\mu_2(1+t)^{-\lambda+\delta} \boldsymbol u)
\\
&\quad+\int ((\boldsymbol u\cdot\nabla)\nabla v)
\cdot((1+t)^{1-\lambda+\delta}\partial_t \boldsymbol u
+\mu_2(1+t)^{-\lambda+\delta} \boldsymbol u)+
\\
&\quad+\int (\varpi v\nabla(\nabla\cdot \boldsymbol u))
\cdot((1+t)^{1-\lambda+\delta}\partial_t \boldsymbol u
+\mu_2(1+t)^{-\lambda+\delta} \boldsymbol u)+\cdots
\\
&=I_{31}+I_{32}+I_{33}+I_{34}+\cdots
\end{align*}
We proceed as before such that
\begin{align*}
I_{31}&=
-\frac{1}{2}\int (1+t)^{1-\lambda+\delta}(\boldsymbol u\cdot\nabla)
|\partial_t \boldsymbol u|^2
-\int \mu_2(1+t)^{-\lambda+\delta} |\boldsymbol u|^2(\nabla\cdot\partial_t \boldsymbol u)
\\
&=
\frac{1}{2}\int (1+t)^{1-\lambda+\delta}(\nabla\cdot\boldsymbol u)
|\partial_t \boldsymbol u|^2
+\int \mu_2(1+t)^{-\lambda+\delta} (\nabla|\boldsymbol u|^2)
)\cdot\partial_t \boldsymbol u,
\end{align*}
and
\begin{align} \nonumber
I_{32}=&
-(1+t)^{1-\lambda+\delta}\int (\varpi v\nabla\partial_t v)
\cdot\partial_t \boldsymbol u
+\mu_2(1+t)^{-\lambda+\delta}\int \varpi \partial_t v\cdot
(\nabla\cdot(v\boldsymbol u))
\\ \label{eq-zI32}
=&
(1+t)^{1-\lambda+\delta}\int \varpi v\partial_t v
(\nabla\cdot\partial_t \boldsymbol u)+\cdots
\end{align}
where the integral in the last inequality of \eqref{eq-zI32} is in the same form as
$I_{13}^1$ but the signs are opposite (such that this one is a good term)
and the time-weight is stronger.
It suffices to modify the definition of $J_1(t)$ in \eqref{eq-zJ1}
by adding a negative integral, which does not affect the inequality
$J_1(t)\lesssim \|v\|_{L^\infty}\cdot\Phi_1^2(t)$ in \eqref{eq-zJ1}.
We also have
\begin{align*}
I_{33}=&
-(1+t)^{1-\lambda+\delta}\int (\nabla\cdot\boldsymbol u)\nabla v
\cdot\partial_t \boldsymbol u
-(1+t)^{1-\lambda+\delta}\int (\boldsymbol u\cdot \nabla v)
\cdot\partial_t (\nabla\cdot\boldsymbol u)
\\
&\qquad
-\int \mu_2(1+t)^{-\lambda+\delta} \nabla|\boldsymbol u|^2\cdot\nabla v
\\
=&
-(1+t)^{1-\lambda+\delta}\int (\boldsymbol u\cdot \nabla v)
\cdot\partial_t (\nabla\cdot\boldsymbol u)+\cdots
\\
=&
(1+t)^{1-\lambda+\delta}\int (\boldsymbol u\cdot \partial_t \boldsymbol u)
\cdot\partial_t (\nabla\cdot\boldsymbol u)+\cdots
\\
=&
(1+t)^{1-\lambda+\delta}\int \boldsymbol u\cdot
\Big(
\nabla\cdot(\partial_t \boldsymbol u\otimes\partial_t \boldsymbol u)
-\frac{1}{2}\nabla|\partial_t \boldsymbol u|^2
\Big)
+\cdots
\\
=&
-(1+t)^{1-\lambda+\delta}\int
(\partial_t \boldsymbol u\otimes\partial_t \boldsymbol u)\odot
(\nabla\boldsymbol u)
+(1+t)^{1-\lambda+\delta}\frac{1}{2}\int
|\partial_t \boldsymbol u|^2(\nabla\cdot \boldsymbol u)
+\cdots
\end{align*}
where ``$\odot$'' denotes the summation of all the element-wise product of two matrices
and we have used the following identity for a general vector-valued function $\boldsymbol\varphi$
(we take $\boldsymbol\varphi=\partial_t \boldsymbol u$)
\begin{equation} \label{eq-id}
(\nabla\cdot \boldsymbol\varphi)\boldsymbol\varphi
=\nabla\cdot(\boldsymbol\varphi\otimes\boldsymbol\varphi)
-\frac{1}{2}\nabla|\boldsymbol\varphi|^2.
\end{equation}
The last integral $I_{34}$ is estimated as follows
\begin{align*}
I_{34}=&
-(1+t)^{1-\lambda+\delta}\int \varpi(\nabla\cdot \boldsymbol u)\nabla v
\cdot\partial_t \boldsymbol u
-(1+t)^{1-\lambda+\delta}\int \varpi v(\nabla\cdot \boldsymbol u)\cdot
\partial_t (\nabla\cdot \boldsymbol u)
\\
&\quad
-\mu_2(1+t)^{-\lambda+\delta} \int \varpi (\nabla\cdot \boldsymbol u))
(\nabla\cdot (v\boldsymbol u))
+\cdots
\\
=&
-(1+t)^{1-\lambda+\delta}\int \varpi v(\nabla\cdot \boldsymbol u)\cdot
\partial_t (\nabla\cdot \boldsymbol u)
+\cdots
\\
=&(1+t)^{1-\lambda+\delta}\int \varpi v\big(\frac{\partial_t v}{1+\varpi v}\big)\cdot
\partial_t (\nabla\cdot \boldsymbol u)
+\cdots
\end{align*}
according to \eqref{eq-znablau}
similar to the treatment of $I_{13}^1$.
Here the integral in the last inequality of the estimate of $I_{34}$
is of the opposite sign compared with $I_{13}^1$
and hence is a good term.
The proof is completed.
$\hfill\Box$

\begin{remark}
From the decay estimates in the proof of Lemma \ref{le-Q},
we see that the inhomogeneous terms involving $b(t)\cdot Q_1$
and the terms involving $v\nabla v$ in $Q_2$ decay slowest
since $b(t)$ is time-asymptotically growing
and $v$ decays slower than $\boldsymbol u$.
\end{remark}

For general integer $k\ge1$, we proceed similarly to deduce the time-weighted energy estimates.
The following ``tame'' product estimate is needed.

\begin{lemma}[\cite{Guo,Tao}] \label{le-tame}
For $1<p<\infty$, $s\ge0$, there holds
$$
\|uv\|_{W^{s,p}}\lesssim
\|u\|_{L^\infty}\|v\|_{W^{s,p}}
+\|v\|_{L^\infty}\|u\|_{W^{s,p}},
$$
for functions $u$ and $v$ in $L^\infty\cap W^{s,p}$.
\end{lemma}

\begin{lemma} \label{le-Q2}
There holds, for integer $k\ge1$, $\lambda\in(-1,0)$, that
$\delta\in(0,\frac{1+\lambda}{4})$, and $|\alpha|=k$
\begin{align*}
&\int \partial_x^\alpha(\partial_tQ_1+b(t)\cdot Q_1-\nabla\cdot Q_2)
\cdot((1+t)^{1+\lambda+\delta}\partial_t \partial_x^\alpha v
+\mu_1(1+t)^{\lambda+\delta} \partial_x^\alpha v)
\\
&\qquad
+\int \partial_x^\alpha(\partial_t Q_2-\nabla Q_1)
\cdot((1+t)^{1-\lambda+\delta}\partial_t \partial_x^\alpha\boldsymbol u
+\mu_2(1+t)^{-\lambda+\delta} \partial_x^\alpha\boldsymbol u)
\\
&\lesssim
\partial_t J_{k+1}(t)+
(\Psi_0(t)+\Phi_N(t))\cdot\Psi_{k+1}^2(t)(1+t)^{-\frac{\delta}{4}}
+(\Psi_0(t)+\Phi_N(t))\cdot\Psi_k^2(t)(1+t)^{-(1+\lambda+\frac{\delta}{4})},
\end{align*}
under the assumption that $\|v\|_{L^\infty}\le \frac{1}{\gamma-1}$,
where
$$
J_{k+1}(t)\lesssim \|v\|_{L^\infty}\cdot\Phi_{k+1}^2(t).
$$
\end{lemma}
{\it\bfseries Proof.}
For $|\alpha|=k\ge1$ and $n\ge3$, we have
\begin{align*}
&\int b(t)\cdot \partial_x^\alpha Q_1
\cdot((1+t)^{1+\lambda+\delta}\partial_t \partial_x^\alpha v
+\mu_1(1+t)^{\lambda+\delta}\partial_x^\alpha v)
\\
&\lesssim \int (|\boldsymbol u\cdot\nabla \partial_x^\alpha v|
+\sum_{j=1}^k|\partial_x^j\boldsymbol u\cdot\nabla \partial_x^{k-j} v|
+|v\nabla\cdot \partial_x^\alpha\boldsymbol u|)
\cdot ((1+t)^{1+\delta}|\partial_t \partial_x^\alpha v|+(1+t)^{\delta} |\partial_x^\alpha v|)
\\
&\lesssim
((1+t)^{1+\delta}\|\partial_t \partial_x^\alpha v\|+
(1+t)^{\delta} \|\partial_x^\alpha v\|)
\cdot
(\|\boldsymbol u\|_{L^\infty}\cdot\|\nabla \partial_x^\alpha v\|
+\|v\|_{L^\infty}\cdot\|\nabla \partial_x^\alpha\boldsymbol u\|)
\\
&\lesssim
((1+t)^{1+\delta}\Psi_{k+1}(t)(1+t)^{-\frac{1+\delta}{2}}+
(1+t)^{\delta}\Psi_k(t)(1+t)^{-\frac{\lambda+\delta}{2}})
\\
&\qquad
\cdot
(\Phi_N(t)(1+t)^{-\frac{1-\lambda+\delta}{2}}
\cdot\Psi_{k+1}(t)(1+t)^{-\frac{\lambda+\delta}{2}}
+\Phi_N(t)(1+t)^{-\frac{1+\lambda+\delta}{2}}
\cdot\Psi_{k+1}(t)(1+t)^{-\frac{-\lambda+\delta}{2}})
\\
&\lesssim
\Phi_N(t)\Psi_{k+1}(t)
(\Psi_{k+1}(t)(1+t)^{-\frac{\delta}{2}}+
\Psi_k(t)(1+t)^{-\frac{1+\lambda+\delta}{2}})
\\
&\lesssim
\Phi_N(t)\Psi_{k+1}^2(t)(1+t)^{-\frac{\delta}{2}}
+\Phi_N(t)\Psi_k^2(t)(1+t)^{-(1+\lambda+\frac{\delta}{2})},
\end{align*}
where we have used \eqref{eq-Linfty} and Lemma \ref{le-tame}.
The case of $n=2$ follows similarly
according to \eqref{eq-Linftyn2} as follows
\begin{align*}
&\int b(t)\cdot \partial_x^\alpha Q_1
\cdot((1+t)^{1+\lambda+\delta}\partial_t \partial_x^\alpha v
+\mu_1(1+t)^{\lambda+\delta}\partial_x^\alpha v)
\\
&\lesssim
((1+t)^{1+\delta}\|\partial_t \partial_x^\alpha v\|+
(1+t)^{\delta} \|\partial_x^\alpha v\|)
\cdot
(\|\boldsymbol u\|_{L^\infty}\cdot\|\nabla \partial_x^\alpha v\|
+\|v\|_{L^\infty}\cdot\|\nabla \partial_x^\alpha\boldsymbol u\|)
\\
&\lesssim
((1+t)^{1+\delta}\Psi_{k+1}(t)(1+t)^{-\frac{1+\delta}{2}}+
(1+t)^{\delta}\Psi_k(t)(1+t)^{-\frac{\lambda+\delta}{2}})
\\
&\qquad
\cdot
((\Psi_0(t)+\Phi_N(t))(1+t)^{-\frac{1-\lambda}{2}-\frac{1+\lambda+\delta}{4}}
\cdot\Psi_{k+1}(t)(1+t)^{-\frac{\lambda+\delta}{2}}
\\
&\qquad\qquad
+(\Psi_0(t)+\Phi_N(t))(1+t)^{-\frac{1+\lambda}{2}-\frac{\delta}{4}}
\cdot\Psi_{k+1}(t)(1+t)^{-\frac{-\lambda+\delta}{2}})
\\
&\lesssim
(\Psi_0(t)+\Phi_N(t))\Psi_{k+1}(t)
(\Psi_{k+1}(t)(1+t)^{-\frac{\delta}{4}}+
\Psi_k(t)(1+t)^{-\frac{1+\lambda}{2}-\frac{\delta}{4}})
\\
&\lesssim
(\Psi_0(t)+\Phi_N(t))\cdot\Psi_{k+1}^2(t)(1+t)^{-\frac{\delta}{4}}
+(\Psi_0(t)+\Phi_N(t))\cdot\Psi_k^2(t)(1+t)^{-(1+\lambda+\frac{\delta}{4})}.
\end{align*}

The other integrals are treated in the same procedure as
in the proof of Lemma \ref{le-Q},
where all the terms involving the $(k+2)$-th order derivatives
are estimated through integration by parts
such that $\Psi_{k+2}(t)$ is not needed.
$\hfill\Box$

\subsection{Closure through Green function method}

We employ the Green function method to deduce the basic energy estimates
in $\Psi_0(t)$.

\begin{lemma} \label{le-Psi0}
There hold, for $\lambda\in(-1,0)$, that
\begin{align*}
\|v\|&\lesssim \|(v_0,\boldsymbol u_0)\|_{L^1\cap L^2}\cdot(1+t)^{-\frac{1+\lambda}{4}n}
+\Psi_0(t)\Phi_N(t)\cdot(1+t)^{-\frac{1+\lambda}{4}n},
\\
\|\boldsymbol u\|&\lesssim
\|(v_0,\boldsymbol u_0)\|_{L^1\cap H^1}\cdot(1+t)^{-\frac{1+\lambda}{4}n-\frac{1-\lambda}{2}}
+(\Phi_N(t)+\Psi_0(t))\Phi_N(t)\cdot(1+t)^{-\frac{1+\lambda}{4}n-\frac{1-\lambda}{2}}.
\end{align*}
\end{lemma}
{\it\bfseries Proof.}
The proof is similar to that of Theorem \ref{th-nonlinear}
but the a priori assumptions are different.
According to the Duhamel principle \eqref{eq-Duhamel}
and the decay estimates of the Green matrix $\mathcal{G}(t,s)$
in Lemma \ref{le-decay-G}, we have
\begin{align*}
\|v(t)\|
\lesssim&
\|\mathcal{G}_{11}(t,0)v_0\|+\|\mathcal{G}_{12}(t,0)\boldsymbol u_0\|
+\int_0^t\|\mathcal{G}_{11}(t,s)Q_1(s)\|ds
+\int_0^t\|\mathcal{G}_{12}(t,s)Q_2(s)\|ds
\\
\lesssim&
\|(v_0,\boldsymbol u_0)\|_{L^1\cap L^2}\cdot(1+t)^{-\frac{1+\lambda}{4}n}
+\int_0^t
\Gamma^{\frac{n}{2}}(t,s)\cdot
(\|Q_1(s)\|_{L^1}^l+\|Q_1(s)\|^h)ds
\\
&+\int_0^t
(1+s)^\lambda\cdot
\Gamma^{\frac{n}{2}+1}(t,s)\cdot
(\|Q_2(s)\|_{L^1}^l+\|Q_2(s)\|^h)ds
\\
\lesssim&
\|(v_0,\boldsymbol u_0)\|_{L^1\cap L^2}\cdot(1+t)^{-\frac{1+\lambda}{4}n}
+\Psi_0(t)\Phi_N(t)\int_0^t
\Gamma^{\frac{n}{2}}(t,s)\cdot
(1+s)^{-\frac{1+\lambda}{4}n-\frac{1-\lambda+\delta}{2}}ds
\\
&+\Psi_0(t)\Phi_N(t)\int_0^t
(1+s)^\lambda\cdot
\Gamma^{\frac{n}{2}+1}(t,s)\cdot
(1+s)^{-\frac{1+\lambda}{4}n-\frac{1+\lambda+\delta}{2}}ds
\\
\lesssim&
\|(v_0,\boldsymbol u_0)\|_{L^1\cap L^2}\cdot(1+t)^{-\frac{1+\lambda}{4}n}
+\Psi_0(t)\Phi_N(t)\cdot(1+t)^{-\frac{1+\lambda}{4}n},
\end{align*}
where we have used Lemma \ref{le-min} (note that
\begin{equation} \label{eq-zdelta}
\begin{cases}
\frac{1+\lambda}{4}n+\frac{1-\lambda+\delta}{2}
\ge \frac{1+\lambda}{2}+\frac{1-\lambda+\delta}{2}=1+\frac{\delta}{2}>1,
\\
\frac{1+\lambda}{4}n+\frac{1+\lambda+\delta}{2}-\lambda
\ge \frac{1+\lambda}{2}+\frac{1+\lambda+\delta}{2}-\lambda=1+\frac{\delta}{2}>1,
\end{cases}
\end{equation}
for all $n\ge2$ and $\lambda\in(-1,0)$)
and the following decay estimates on $\|Q(s)\|_{L^1}$ and $\|Q(s)\|$
(here and after, we use $D^j:=\partial_x^j$)
\begin{align*}
\|Q_1(s)\|_{L^1}
&\lesssim \|uDv\|_{L^1}+\|vDu\|_{L^1}
\lesssim \|u\|\|Dv\|+\|v\|\|Du\|
\\
&\lesssim
\Psi_0(s)(1+s)^{-\frac{1+\lambda}{4}n-\frac{1-\lambda}{2}}
\cdot \Phi_N(s)(1+s)^{-\frac{1+\lambda+\delta}{2}}
+\Psi_0(s)(1+s)^{-\frac{1+\lambda}{4}n}
\cdot \Phi_N(s)(1+s)^{-\frac{1-\lambda+\delta}{2}}
\\
&\lesssim
\Psi_0(s)\Phi_N(s)\cdot(1+s)^{-\frac{1+\lambda}{4}n-\frac{1-\lambda+\delta}{2}},
\\
\|Q_2(s)\|_{L^1}
&\lesssim \|uDu\|_{L^1}+\|vDv\|_{L^1}
\lesssim \|u\|\|Du\|+\|v\|\|Dv\|
\\
&\lesssim
\Psi_0(s)(1+s)^{-\frac{1+\lambda}{4}n-\frac{1-\lambda}{2}}
\cdot \Phi_N(s)(1+s)^{-\frac{1-\lambda+\delta}{2}}
+\Psi_0(s)(1+s)^{-\frac{1+\lambda}{4}n}
\cdot \Phi_N(s)(1+s)^{-\frac{1+\lambda+\delta}{2}}
\\
&\lesssim
\Psi_0(s)\Phi_N(s)\cdot(1+s)^{-\frac{1+\lambda}{4}n-\frac{1+\lambda+\delta}{2}}.
\end{align*}
The decay estimates on $\|Q_1\|$ and $\|Q_2\|$ are
at least at the same rates
as $\|Q_1\|_{L^1}$ and $\|Q_2\|_{L^1}$
since the estimates on $\|Dv\|_{L^\infty}$ and $\|D\boldsymbol u\|_{L^\infty}$
decay at the same rates as $\|Dv\|$ and $\|D\boldsymbol u\|$
according to \eqref{eq-Linfty} and \eqref{eq-Linftyn2}.

In order to deduce the optimal decay estimate on $\|\boldsymbol u\|$
we need to utilize the optimal decay estimate on $\mathcal{G}_{22}$
in \eqref{eq-decay-G-opt}, which needs $\|DQ_2\|$,
instead of \eqref{eq-decay-G}, which only needs $\|Q_2\|$.
We see that
\begin{align*}
\|DQ_2(s)\|
&\lesssim \|uD^2u\|+\|DuDu\|+\|vD^2v\|+\|DvDv\|
\\
&\lesssim \|u\|_{L^\infty}\|D^2u\|+\|Du\|_{L^\infty}\|Du\|
+\|v\|_{L^\infty}\|D^2v\|+\|Dv\|_{L^\infty}\|Dv\|
\\
&\lesssim
\begin{cases}
\Phi_N^2(s)(1+s)^{-(1-\lambda+\delta)}+
\Phi_N^2(s)(1+s)^{-(1+\lambda+\delta)}, ~ & n\ge3,
\\
(\Phi_N(s)+\Psi_0(s))(1+s)^{-\frac{1-\lambda}{2}-\frac{1+\lambda+\delta}{4}}
\cdot \Phi_N(s)(1+s)^{-\frac{1-\lambda+\delta}{2}}
\\
\quad+\Phi_N^2(s)(1+s)^{-(1-\lambda+\delta)}+\Phi_N^2(s)(1+s)^{-(1+\lambda+\delta)}
\\
\quad+(\Phi_N(s)+\Psi_0(s))(1+s)^{-\frac{1+\lambda}{2}-\frac{\delta}{4}}
\cdot \Phi_N(s)(1+s)^{-\frac{1+\lambda+\delta}{2}}, ~ & n=2,
\end{cases}
\\
&\lesssim
\begin{cases}
\Phi_N^2(s)(1+s)^{-(1+\lambda+\delta)}, ~ & n\ge3,
\\
(\Phi_N(s)+\Psi_0(s))\Phi_N(s)\cdot(1+s)^{-(1+\lambda+\frac{3}{4}\delta)}, ~ & n=2,
\end{cases}
\\
&\lesssim
(\Phi_N(s)+\Psi_0(s))\Phi_N(s)\cdot(1+s)^{-(1+\lambda+\frac{3}{4}\delta)}, \quad n\ge2,
\end{align*}
according to \eqref{eq-Linfty} and \eqref{eq-Linftyn2}.
Therefore, we have
\begin{align*}
\|\boldsymbol u(t)\|
\lesssim&
\|\mathcal{G}_{21}(t,0)v_0\|+\|\mathcal{G}_{22}(t,0)\boldsymbol u_0\|
+\int_0^t\|\mathcal{G}_{21}(t,s)Q_1(s)\|ds
+\int_0^t\|\mathcal{G}_{22}(t,s)Q_2(s)\|ds
\\
\lesssim&
\|(v_0,\boldsymbol u_0)\|_{L^1\cap H^1}\cdot(1+t)^{-\frac{1+\lambda}{4}n}
+\int_0^t
(1+t)^\lambda\cdot
\Gamma^{\frac{n}{2}+1}(t,s)\cdot
(\|Q_1(s)\|_{L^1}^l+\|Q_1(s)\|^h)ds
\\
&+\int_0^t
(1+t)^\lambda(1+s)^\lambda\cdot
\Gamma^{\frac{n}{2}+2}(t,s)\cdot
(\|Q_2(s)\|_{L^1}^l+\|DQ_2(s)\|^h)ds
\\
\lesssim&
\|(v_0,\boldsymbol u_0)\|_{L^1\cap H^1}\cdot
(1+t)^{-\frac{1+\lambda}{4}n-\frac{1-\lambda}{2}}
+\Psi_0(t)\Phi_N(t)\int_0^t
(1+t)^\lambda\cdot
\Gamma^{\frac{n}{2}+1}(t,s)\cdot
(1+s)^{-\frac{1+\lambda}{4}n-\frac{1-\lambda+\delta}{2}}ds
\\
&+(\Phi_N(t)+\Psi_0(t))\Phi_N(t)\int_0^t
(1+t)^\lambda(1+s)^\lambda\cdot
\Gamma^{\frac{n}{2}+2}(t,s)
\cdot(1+s)^{-(1+\lambda+\frac{3}{4}\delta)}ds
\\
\lesssim&
\|(v_0,\boldsymbol u_0)\|_{L^1\cap H^1}\cdot
(1+t)^{-\frac{1+\lambda}{4}n-\frac{1-\lambda}{2}}
+(\Phi_N(t)+\Psi_0(t))\Phi_N(t)\cdot(1+t)^{-\frac{1+\lambda}{4}n-\frac{1-\lambda}{2}},
\end{align*}
since
\begin{equation} \label{eq-zcom1}
\begin{cases}
\frac{1+\lambda}{4}n+\frac{1-\lambda+\delta}{2}-\lambda
\ge \frac{1+\lambda}{2}+\frac{1-\lambda+\delta}{2}-\lambda=1-\lambda+\frac{\delta}{2}>1,
\quad
&\frac{1+\lambda}{4}n+\frac{1-\lambda+\delta}{2}-\lambda
\ge \frac{1+\lambda}{4}n+\frac{1-\lambda}{2},
\\
1+\lambda+\frac{3}{4}\delta-2\lambda=1-\lambda+\frac{3}{4}\delta>1,
\quad
&1+\lambda+\frac{3}{4}\delta-2\lambda
\ge \frac{1+\lambda}{4}n+\frac{1-\lambda}{2},
\end{cases}
\end{equation}
for all $n\ge2$ and $\lambda\in(-1,0)$,
except that the last inequality in \eqref{eq-zcom1} is not true
for the case of $\frac{1+\lambda}{2}n>1-\lambda+\frac{3}{2}\delta$.
Fortunately, this case has already been proved in Theorem \ref{th-nonlinear}
by means of the Green function method
(for $\lambda\in(-\frac{n}{n+2},0)$, i.e., $\frac{1+\lambda}{2}n>-\lambda$,
which covers the exceptional case here).
The proof is completed.
$\hfill\Box$

\begin{remark}
The introducing of the positive constant $\delta$ plays an important role
in the closure of the optimal decay estimate of $\|v\|$
(especially for the case of $n=2$)
according to the condition \eqref{eq-zdelta}
in the proof of Lemma \ref{le-Psi0}.
\end{remark}

We combine the above time-weighted iteration scheme and Green function method
to close the decay estimates for $\lambda\in(-1,0)$.

\begin{proposition} \label{th-com-0}
For $n\ge2$, $N\ge[\frac{n}{2}]+2$ and $\lambda\in(-1,0)$,
there exists a constant $\varepsilon_0>0$ such that
the solution $(v,\boldsymbol u)$ of the nonlinear system \eqref{eq-vbdu}
corresponding to small initial data
$\|(v_0,\boldsymbol u_0)\|_{L^1\cap H^N}\le \varepsilon_0$
exists globally and satisfies
\begin{equation} \label{eq-com-0}
\begin{cases}
\|v(t)\|\lesssim (1+t)^{-\frac{1+\lambda}{4}n}, \\
\|\boldsymbol u(t)\|\lesssim (1+t)^{-\frac{1+\lambda}{4}n-\frac{1-\lambda}{2}}.
\end{cases}
\end{equation}
The above decay rates are optimal and consistent with the optimal decay rates
of the linearized hyperbolic system.
\end{proposition}
{\it\bfseries Proof.}
We claim that the following a priori decay estimate
\begin{equation} \label{eq-apriori-com}
\Phi_N(t)+\Psi_0(t)\le \delta_0,
\end{equation}
holds for all the time $t>0$,
under the small energy assumption of initial data
$\|(v_0,\boldsymbol u_0)\|_{L^1\cap H^N}\le \varepsilon_0$,
where $\varepsilon_0$ and $\delta_0$ are positive constants to be determined.
In fact, Lemma \ref{le-Psi0} tells us that
\begin{equation} \label{eq-zPsi0}
\Psi_0(T)\le \sup_{t\in(0,T)}\Big\{
(1+t)^{\frac{1+\lambda}{4}n}\|v\|,
(1+t)^{\frac{1+\lambda}{4}n+\frac{1-\lambda}{2}}\|\boldsymbol u\|
\Big\}
\lesssim
\varepsilon_0+\delta_0^2.
\end{equation}
Substituting the estimates of inhomogeneous terms in Lemma \ref{le-Q} and Lemma \ref{le-Q2}
into the time-weighted iteration scheme \eqref{eq-it-0} and \eqref{eq-it-1}
in Lemma \ref{le-iteration},
we have for integer $0\le k\le N-1$ that
\begin{align*}
&\Phi_1^2(t)+\int_0^t \Psi_1^2(s)ds
\\
&\quad\lesssim
\Phi_1^2(0)+J_1(t)+
\int_0^t (1+s)^{-1-\frac{1+\lambda}{2}n+\delta}\cdot\Psi_0^2(s)ds
+\delta_0\int_0^t \Psi_1^2(s)ds
+\delta_0\int_0^t (1+s)^{-1-\frac{1+\lambda}{4}}\cdot\Psi_0^2(s)ds,
\\
&\Phi_{k+1}^2(t)+\int_0^t \Psi_{k+1}^2(s)ds
\\
&\quad\lesssim
\Phi_{k+1}^2(0)+J_{k+1}(t)+
\int_0^t (1+s)^{-1-\lambda}\cdot\Psi_k^2(s)ds
+\delta_0\int_0^t \Psi_{k+1}^2(s)ds
+\delta_0\int_0^t (1+s)^{-1-\lambda}\cdot\Psi_k^2(s)ds,
\end{align*}
where
\begin{align*}
&J_1(t)\lesssim \|v\|_{L^\infty}\cdot\Phi_1^2(t)
\lesssim (\Phi_N(t)+\Psi_0(t))\cdot\Phi_1^2(t)
\lesssim \delta_0\Phi_1^2(t)
\\
&J_{k+1}(t)\lesssim \|v\|_{L^\infty}\cdot\Phi_{k+1}^2(t)
\lesssim (\Phi_N(t)+\Psi_0(t))\Phi_{k+1}^2(t)
\lesssim \delta_0\Phi_{k+1}^2(t).
\end{align*}
We note that $\delta_0$ and $\varepsilon_0$ are small such that
the above inequalities can be simplified as
\begin{align} \label{eq-zPhi0}
&\Phi_1^2(t)+\int_0^t \Psi_1^2(s)ds
\lesssim
\varepsilon_0^2+
\Psi_0^2(t)\int_0^t (1+s)^{-1-\frac{1+\lambda}{4}}ds,
\\ \label{eq-zPhi}
&\Phi_{k+1}^2(t)+\int_0^t \Psi_{k+1}^2(s)ds
\lesssim
\varepsilon_0^2+
\int_0^t (1+s)^{-1-\lambda}\cdot\Psi_k^2(s)ds.
\end{align}
Multiplying \eqref{eq-zPhi} by small positive constants
for $0\le k\le N-1$, summing the resulting inequalities up together with \eqref{eq-zPhi0},
we have
\begin{align*}
\sum_{1\le j\le N}\Phi_j^2(t)
\lesssim
\varepsilon_0^2+
\Psi_0^2(t)\int_0^t (1+s)^{-1-\frac{1+\lambda}{4}}ds
\lesssim \varepsilon_0^2+(\varepsilon_0+\delta_0^2)^2,
\end{align*}
according to the estimate \eqref{eq-zPsi0}.
Therefore,
$$
\Phi_N(t)+\Psi_0(t)
\lesssim \varepsilon_0+\delta_0^2
\le \delta_0,
$$
for positive constants $\varepsilon_0$ and $\delta_0$ small enough.

We can show that the decay estimates \eqref{eq-com-0} are optimal
in a similar way as the proof of Theorem \ref{th-nonlinear},
just replacing the estimates on $\|v\|$ and $\|\boldsymbol u\|$
by those in Lemma \ref{le-Psi0}.
The proof is completed.
$\hfill\Box$

\vskip2mm
{\it\bfseries Proof of Theorem \ref{th-com-0-in}.} It is immediately proved from
 Proposition \ref{th-com-0}.
$\hfill\Box$

\section{Critical case of $\lambda=-1$: optimal logarithmic decays}

This section is devoted to the critical case of $\lambda=-1$.
We show the optimal decay estimates such that
$\|v(t)\|$ decays as powers of $\ln(e+t)$,
that is, $\|v(t)\|\approx|\ln(e+t)|^{-\frac{n}{4}}$.

We start with the optimal decay estimates of the Green matrix for the
critical case of $\lambda=-1$, which are special cases of Lemma \ref{le-decay-G}.
Here we write it down for the sake of convenience.

\begin{lemma} \label{le-decay-G-1}
For $\lambda=-1$, there hold
\begin{align} \nonumber
&\|\partial_x^\alpha \mathcal{G}_{11}(t,s)\phi(x)\|\lesssim
\Big(1+\ln\Big(\frac{1+t}{1+s}\Big)\Big)^{-\frac{1}{2}(\frac{n}{2}+|\alpha|)}
(\|\phi\|_{L^1}^l+\|\partial_x^{|\alpha|}\phi\|^h),
\\ \nonumber
&\|\partial_x^\alpha \mathcal{G}_{12}(t,s)\phi(x)\|\lesssim
(1+s)^{-1}\cdot
\Big(1+\ln\Big(\frac{1+t}{1+s}\Big)\Big)^{-\frac{1}{2}(\frac{n}{2}+|\alpha|+1)}
(\|\phi\|_{L^1}^l+\|\partial_x^{|\alpha|}\phi\|^h),
\\ \nonumber
&\|\partial_x^\alpha \mathcal{G}_{21}(t,s)\phi(x)\|\lesssim
(1+t)^{-1}\cdot
\Big(1+\ln\Big(\frac{1+t}{1+s}\Big)\Big)^{-\frac{1}{2}(\frac{n}{2}+|\alpha|+1)}
(\|\phi\|_{L^1}^l+\|\partial_x^{|\alpha|}\phi\|^h),
\\ \label{eq-decay-G-1}
&\|\partial_x^\alpha \mathcal{G}_{22}(t,s)\phi(x)\|\lesssim
\Big(\frac{1+t}{1+s}\Big)^{-1}\cdot
\Big(1+\ln\Big(\frac{1+t}{1+s}\Big)\Big)^{-\frac{1}{2}(\frac{n}{2}+|\alpha|)}
(\|\phi\|_{L^1}^l+\|\partial_x^{|\alpha|}\phi\|^h).
\end{align}

Moreover,
\begin{align}
\label{eq-decay-G-opt-1}
&\|\partial_x^\alpha \mathcal{G}_{22}(t,s)\phi(x)\|\lesssim
(1+t)^{-1}(1+s)^{-1}\cdot
\Big(1+\ln\Big(\frac{1+t}{1+s}\Big)\Big)^{-\frac{1}{2}(\frac{n}{2}+|\alpha|+2)}
(\|\phi\|_{L^1}^l+\|\partial_x^{|\alpha|+1}\phi\|^h).
\end{align}
\end{lemma}
{\it\bfseries Proof.}
These estimates are simple conclusions of Theorem \ref{th-linear} in Appendix.
$\hfill\Box$

\vskip2mm
The following time decay estimate of the ``convolution'' type integral
of two critical time decay functions involving logarithm
plays an essential role in the Green function method for $\lambda=-1$.

\begin{lemma}[Logarithmic time decay functions] \label{le-min-1}
For $\beta>0$ and $\gamma>1$, there holds (we may assume that $t\ge1$)
\begin{equation} \label{eq-min-1}
\int_0^t
\Big(1+\ln\Big(\frac{1+t}{1+s}\Big)\Big)^{-\beta}
(1+s)^{-1}|\ln(e+s)|^{-\gamma}ds
\approx
\begin{cases}
|\ln(e+t)|^{-\min\{\beta,\gamma-1\}},  &\gamma>1,\\
\ln(\ln(e^e+t)),  &\gamma=1,\\
|\ln(e+t)|^{1-\gamma}, &\gamma<1.
\end{cases}
\end{equation}
\end{lemma}
{\it\bfseries Proof.}
For $\gamma\le1$, we have
\begin{align*}
\int_0^t
\Big(1+\ln\Big(\frac{1+t}{1+s}\Big)\Big)^{-\beta}
(1+s)^{-1}|\ln(e+s)|^{-\gamma}ds
&\lesssim
\int_0^t
(e+s)^{-1}|\ln(e+s)|^{-\gamma}ds
\\
&\lesssim
\begin{cases}
\ln(\ln(e^e+t)),  \quad &\gamma=1,
\\
|\ln(e+t)|^{1-\gamma}, \quad & \gamma<1,
\end{cases}
\end{align*}
and
\begin{align*}
&\int_0^t
\Big(1+\ln\Big(\frac{1+t}{1+s}\Big)\Big)^{-\beta}
(1+s)^{-1}|\ln(e+s)|^{-\gamma}ds
\\
&\gtrsim
\int_\frac{t}{2}^t
\Big(1+\ln\Big(\frac{1+t}{1+s}\Big)\Big)^{-\beta}
(1+s)^{-1}|\ln(e+s)|^{-\gamma}ds
\\
&\approx
\int_\frac{t}{2}^t
(e+s)^{-1}|\ln(e+s)|^{-\gamma}ds
\\
&\approx
\begin{cases}
\ln(\ln(e^e+t)),  \quad &\gamma=1,
\\
|\ln(e+t)|^{1-\gamma}, \quad & \gamma<1.
\end{cases}
\end{align*}

For $\gamma>1$, we calculate the integral
divided into $(0,t^\varepsilon)$ and $(t^\varepsilon,t)$,
where $\varepsilon\in(0,1)$ is a small constant to be determined,
as follows
\begin{align} \nonumber
&\int_{t^\varepsilon}^t
\Big(1+\ln\Big(\frac{1+t}{1+s}\Big)\Big)^{-\beta}
(1+s)^{-1}|\ln(e+s)|^{-\gamma}ds
\\ \nonumber
&\lesssim
\int_{t^\varepsilon}^t
(e+s)^{-1}|\ln(e+s)|^{-\gamma}ds
\\ \label{eq-zlog1}
&\approx
|\ln(e+t^\varepsilon)|^{-(\gamma-1)}
\approx
|\varepsilon \ln(e+t)|^{-(\gamma-1)},
\end{align}
and
\begin{align} \nonumber
&\int_0^{t^\varepsilon}
\Big(1+\ln\Big(\frac{1+t}{1+s}\Big)\Big)^{-\beta}
(1+s)^{-1}|\ln(e+s)|^{-\gamma}ds
\\ \nonumber
&\approx
\int_0^{t^\varepsilon}
\Big(1+\ln\Big(\frac{e+t}{e+s}\Big)\Big)^{-\beta}
(e+s)^{-1}|\ln(e+s)|^{-\gamma}ds
\\ \nonumber
&=\int_0^{t^\varepsilon}
\Big(1+\ln\Big(\frac{e+t}{e+s}\Big)\Big)^{-\beta}
d\Big(\frac{-1}{\gamma-1}|\ln(e+s)|^{-(\gamma-1)}\Big)
\\ \nonumber
&=\Big[\frac{-1}{\gamma-1}|\ln(e+s)|^{-(\gamma-1)}
\Big(1+\ln\Big(\frac{e+t}{e+s}\Big)\Big)^{-\beta}
\Big]_0^{t^\varepsilon}
\\ \label{eq-zlog}
&\quad+\int_0^{t^\varepsilon}
\frac{\beta}{\gamma-1}|\ln(e+s)|^{-(\gamma-1)}
(e+s)^{-1}
\Big(1+\ln\Big(\frac{e+t}{e+s}\Big)\Big)^{-\beta-1}ds.
\end{align}
Now we fix $\varepsilon>0$ to be sufficiently small such that
\begin{align*}
&\int_0^{t^\varepsilon}
\frac{\beta}{\gamma-1}|\ln(e+s)|^{-(\gamma-1)}
(e+s)^{-1}
\Big(1+\ln\Big(\frac{e+t}{e+s}\Big)\Big)^{-\beta-1}ds
\\
&\le
\frac{1}{2}
\int_0^{t^\varepsilon}
\Big(1+\ln\Big(\frac{e+t}{e+s}\Big)\Big)^{-\beta}
(e+s)^{-1}|\ln(e+s)|^{-\gamma}ds,
\end{align*}
one of whose sufficient conditions is
\begin{align*}
\frac{\beta}{\gamma-1}|\ln(e+s)|
\Big(1+\ln\Big(\frac{e+t}{e+s}\Big)\Big)^{-1}
\le
\frac{1}{2},
\quad \forall s\in(0,t^\varepsilon).
\end{align*}
It suffices to take
$\frac{\beta}{\gamma-1}\cdot\frac{\varepsilon}{1-\varepsilon}\le\frac{1}{2}$,
which is true for a small $\varepsilon\in(0,1)$.
Now \eqref{eq-zlog} reads as
\begin{align*}
&\int_0^{t^\varepsilon}
\Big(1+\ln\Big(\frac{1+t}{1+s}\Big)\Big)^{-\beta}
(1+s)^{-1}|\ln(e+s)|^{-\gamma}ds
\\
&\approx
\Big[\frac{-1}{\gamma-1}|\ln(e+s)|^{-(\gamma-1)}
\Big(1+\ln\Big(\frac{e+t}{e+s}\Big)\Big)^{-\beta}
\Big]_0^{t^\varepsilon}
\\
&\approx
|\ln(e+t)|^{-\beta}
-|\ln(e+t)|^{-(\gamma-1)-\beta}
\approx |\ln(e+t)|^{-\beta}.
\end{align*}
On the other hand, we can improve \eqref{eq-zlog1} as
$$
\int_{t^\varepsilon}^t
\Big(1+\ln\Big(\frac{1+t}{1+s}\Big)\Big)^{-\beta}
(1+s)^{-1}|\ln(e+s)|^{-\gamma}ds
\approx |\ln(e+t)|^{-(\gamma-1)},
$$
since
\begin{align*}
&\int_{t^\varepsilon}^t
\Big(1+\ln\Big(\frac{1+t}{1+s}\Big)\Big)^{-\beta}
(1+s)^{-1}|\ln(e+s)|^{-\gamma}ds
\\
&\gtrsim
\int_\frac{t}{2}^t
(e+s)^{-1}|\ln(e+s)|^{-\gamma}ds
\\
&\approx
|\ln(e+t/2)|^{-(\gamma-1)}
\approx
|\ln(e+t)|^{-(\gamma-1)}.
\end{align*}
The proof is completed.
$\hfill\Box$

\vskip2mm
We apply the time-weighted iteration scheme developed in Section \ref{sec-iteration}
to the critical case of $\lambda=-1$.

\begin{lemma} \label{le-timew-1}
For any nonnegative integer $k$, $\lambda=-1$,
$\delta\in(0,\frac{n}{2})$, and $|\alpha|=k$, there hold
\begin{align} \nonumber
&\frac{d}{dt}\int
E^v(\partial_t \partial_x^\alpha v,\nabla \partial_x^\alpha v,\partial_x^\alpha v)
+\int \big[
(1+t)\cdot|\ln(e+t)|^{\delta+1}|\partial_t \partial_x^\alpha v|^2
+(1+t)^{-1}\cdot|\ln(e+t)|^{\delta}|\nabla \partial_x^\alpha v|^2
\big]
\\ \nonumber
&\lesssim
\int (1+t)^{-1}\cdot|\ln(e+t)|^{\delta-1}(\partial_x^\alpha v)^2
\\ \label{eq-energy-v-1}
&\qquad
+\int \partial_x^\alpha(\partial_tQ_1+b(t)\cdot Q_1-\nabla\cdot Q_2)
\cdot(|\ln(e+t)|^{\delta+1}\partial_t \partial_x^\alpha v
+\mu_1(1+t)^{-1}\cdot|\ln(e+t)|^{\delta} \partial_x^\alpha v),
\end{align}
and
\begin{align} \nonumber
&\frac{d}{dt}\int
E^u(\partial_t \partial_x^\alpha\boldsymbol u,
\nabla \partial_x^\alpha\boldsymbol u,\partial_x^\alpha\boldsymbol u)
+\int \big[
(1+t)^{3}\cdot|\ln(e+t)|^{\delta}|\partial_t \partial_x^\alpha\boldsymbol u|^2
+(1+t)\cdot|\ln(e+t)|^{\delta}|\nabla \partial_x^\alpha\boldsymbol u|^2
\big]
\\ \nonumber
&\lesssim
\int (1+t)\cdot|\ln(e+t)|^{\delta}
|\partial_x^\alpha\boldsymbol u|^2
\\ \label{eq-energy-u-1}
&\qquad+\int \partial_x^\alpha(\partial_t Q_2-\nabla Q_1)
\cdot((1+t)^{2}\cdot|\ln(e+t)|^{\delta}\partial_t \partial_x^\alpha\boldsymbol u
+\mu_2(1+t)\cdot|\ln(e+t)|^{\delta} \partial_x^\alpha\boldsymbol u),
\end{align}
where $\mu_1>0$ and $\mu_2>0$ are constants and
\begin{align*}
E^v(\partial_t \partial_x^\alpha v,\nabla \partial_x^\alpha v,\partial_x^\alpha v)
&\approx
|\ln(e+t)|^{\delta+1}(|\partial_t \partial_x^\alpha v|^2+|\nabla \partial_x^\alpha v|^2)
+|\ln(e+t)|^{\delta}(\partial_x^\alpha v)^2,
\\
E^u(\partial_t \partial_x^\alpha \boldsymbol u,\nabla \partial_x^\alpha \boldsymbol u,
\partial_x^\alpha \boldsymbol u)
&\approx
(1+t)^{2}\cdot|\ln(e+t)|^{\delta}(|\partial_t \partial_x^\alpha \boldsymbol u|^2
+|\nabla \partial_x^\alpha \boldsymbol u|^2)
+(1+t)^{2}\cdot|\ln(e+t)|^{\delta}|\partial_x^\alpha \boldsymbol u|^2.
\end{align*}
\end{lemma}
{\it\bfseries Proof.}
This is proved by multiplying \eqref{eq-wave-v} by
$$|\ln(e+t)|^{\delta+1}\partial_t \partial_x^\alpha v
+\mu_1(1+t)^{-1}\cdot|\ln(e+t)|^{\delta} \partial_x^\alpha v$$
and multiplying \eqref{eq-wave-u} by
$$(1+t)^{2}\cdot|\ln(e+t)|^{\delta}\partial_t \partial_x^\alpha\boldsymbol u
+\mu_2(1+t)\cdot|\ln(e+t)|^{\delta} \partial_x^\alpha\boldsymbol u$$
with $\delta\in(0,\frac{n}{2})$ and $\mu_1,\mu_2>0$.
We note that the time-weight of $\partial_t \partial_x^\alpha v$ is
$|\ln(e+t)|^{\delta+1}$ instead of $|\ln(e+t)|^{\delta}$.
The reason is that the time-weights are chosen such that
$$
\partial_t(|\ln(e+t)|^{\delta+1})
\approx (1+t)^{-1}\cdot|\ln(e+t)|^{\delta},
$$
and
$$
\partial_t((1+t)^{2}\cdot|\ln(e+t)|^{\delta})
\approx
(1+t)\cdot|\ln(e+t)|^{\delta}.
$$
The rest of the proof is similar to Lemma \ref{le-timew}.
We omit the details.
$\hfill\Box$

\vskip2mm
We define the following time-weighted energies
for the critical case of $\lambda=-1$, $N\ge[\frac{n}{2}]+2$ and $0\le k\le N-1$,
\begin{align} \nonumber
\Phi_{k+1}(T):=\sup_{t\in(0,T)}\Big\{&
\sum_{|\alpha|=k}\Big[
|\ln(e+t)|^{\delta+1}\int(|\partial_t \partial_x^\alpha v|^2+|\nabla \partial_x^\alpha v|^2)
\\ \label{eq-energy-Phi-1}
&+(1+t)^{2}\cdot|\ln(e+t)|^{\delta}\int(|\partial_t \partial_x^\alpha\boldsymbol u|^2
+|\nabla \partial_x^\alpha\boldsymbol u|^2)
\Big]
\Big\}^\frac{1}{2},
\end{align}
and
\begin{align} \nonumber
\Psi_{k+1}(T):=\sup_{t\in(0,T)}\Big\{&
\sum_{|\alpha|=k}\Big[\int \big[
(1+t)\cdot|\ln(e+t)|^{\delta+1}|\partial_t \partial_x^\alpha v|^2
+(1+t)^{-1}\cdot|\ln(e+t)|^{\delta}|\nabla \partial_x^\alpha v|^2
\big]
\\ \label{eq-energy-Psi-1}
&+\int \big[
(1+t)^{3}\cdot|\ln(e+t)|^{\delta}|\partial_t \partial_x^\alpha\boldsymbol u|^2
+(1+t)\cdot|\ln(e+t)|^{\delta}|\nabla \partial_x^\alpha\boldsymbol u|^2
\big]
\Big]
\Big\}^\frac{1}{2}.
\end{align}
We may assume that $\Phi_{k+1}(T)\ge \Phi_k(T)$ for all $k\ge1$ and $T$.
Similar to the case of $\lambda\in(-1,0)$, here for $\lambda=-1$
the energy $\Phi_{k+1}(T)$ is defined according to the time-weighted energy estimates
in Lemma \ref{le-timew-1},
but the decay estimates on $\|v\|$ and $\|\boldsymbol u\|$ are absent.
Therefore, we define the following weighted energy
\begin{align} \label{eq-energy-Psi0-1}
\Psi_0(T):=\sup_{t\in(0,T)}\Big\{&
|\ln(e+t)|^{\frac{n}{4}}\|v\|,
(1+t)\cdot |\ln(e+t)|^{\frac{n}{4}+\frac{1}{2}}
\|\boldsymbol u\|
\Big\}.
\end{align}
The energy estimates in $\Psi_0(T)$ will be closed through the Green function method
instead of the time-weighted energy method.
There still holds
\begin{equation} \label{eq-initial-1}
\|(v_0,\boldsymbol u_0)\|_{H^N}
\approx \sum_{k=1}^{N}\Phi_k(0)+\Psi_0(0)
\approx \Phi_N(0)+\Psi_0(0).
\end{equation}
According to Sobolev embedding theorem, we have
\begin{align} \nonumber
&|\ln(e+t)|^{\frac{\delta+1}{2}}\|\partial_x^j v\|_{L^\infty}
+(1+t)\cdot|\ln(e+t)|^{\frac{\delta}{2}}\|\partial_x^j \boldsymbol u\|_{L^\infty}
\\ \label{eq-Linfty-1}
&\lesssim \max_{1\le k\le [\frac{n}{2}]+2}\Phi_k(t)
\lesssim \Phi_N(t),
\qquad 0\le j\le 1, n\ge3,
\end{align}
and
\begin{align} \nonumber
&|\ln(e+t)|^{\frac{1}{2}+\frac{\delta}{4}}\|v\|_{L^\infty}
+(1+t)\cdot|\ln(e+t)|^{\frac{1}{2}+\frac{\delta}{4}}\|\boldsymbol u\|_{L^\infty}
\\ \nonumber
&\qquad+|\ln(e+t)|^{\frac{\delta+1}{2}}\|\partial_x v\|_{L^\infty}
+(1+t)\cdot|\ln(e+t)|^{\frac{\delta}{2}}\|\partial_x \boldsymbol u\|_{L^\infty}
\\ \label{eq-Linftyn2-1}
&\lesssim \max_{1\le k\le [\frac{n}{2}]+2}\Phi_k(t)+\Psi_0(t)
\lesssim \Phi_N(t)+\Psi_0(t),
\qquad n=2.
\end{align}

We have the following iteration scheme based on Lemma \ref{le-timew-1}
for the critical case of $\lambda=-1$.

\begin{lemma}[Time-weighted iteration scheme] \label{le-iteration-1}
For $\lambda=-1$ and $\delta\in(0,\frac{n}{2})$, there holds
\begin{align} \nonumber
&\Phi_1^2(t)+\int_0^t \Psi_1^2(s)ds
\\ \nonumber
&\lesssim
\Phi_1^2(0)+
\int_0^t (1+s)^{-1}\cdot|\ln(e+s)|^{\delta-1-\frac{n}{2}}\cdot\Psi_0^2(s)ds
\\ \nonumber
&\quad+\int_0^t\int (\partial_tQ_1+b(s)\cdot Q_1-\nabla\cdot Q_2)
\cdot(|\ln(e+s)|^{\delta+1}\partial_t v
+\mu_1(1+s)^{-1}\cdot|\ln(e+s)|^{\delta} v)ds
\\ \label{eq-it-0-1}
&\quad+\int_0^t\int (\partial_t Q_2-\nabla Q_1)
\cdot((1+s)^{2}\cdot|\ln(e+s)|^{\delta}\partial_t \boldsymbol u
+\mu_2(1+s)\cdot|\ln(e+s)|^{\delta} \boldsymbol u)ds,
\end{align}
and for any integer $k\ge1$, there holds
\begin{align} \nonumber
&\Phi_{k+1}^2(t)+\int_0^t \Psi_{k+1}^2(s)ds
\\ \nonumber
&\lesssim
\Phi_{k+1}^2(0)+
\int_0^t \Psi_k^2(s)ds
\\ \nonumber
&+\sum_{|\alpha|=k}
\int_0^t\int \partial_x^\alpha(\partial_tQ_1+b(s)\cdot Q_1-\nabla\cdot Q_2)
\cdot(|\ln(e+s)|^{\delta+1}\partial_t \partial_x^\alpha v
+\mu_1(1+s)^{-1}\cdot|\ln(e+s)|^{\delta} \partial_x^\alpha v)ds
\\ \label{eq-it-1-1}
&+\sum_{|\alpha|=k}
\int_0^t\int \partial_x^\alpha(\partial_t Q_2-\nabla Q_1)
\cdot((1+s)^{2}\cdot|\ln(e+s)|^{\delta}\partial_t \partial_x^\alpha\boldsymbol u
+\mu_2(1+s)\cdot|\ln(e+s)|^{\delta} \partial_x^\alpha\boldsymbol u)ds.
\end{align}
\end{lemma}
{\it\bfseries Proof.}
These are conclusions of Lemma \ref{le-timew-1}
with the notations $\Phi_k(t)$, $\Psi_k(t)$, and $\Psi_0(t)$ defined by \eqref{eq-energy-Phi-1},
\eqref{eq-energy-Psi-1}, and \eqref{eq-energy-Psi0-1}.
We note that
\begin{align*}
\int (1+t)^{-1}\cdot|\ln(e+t)|^{\delta-1}|v|^2
&\lesssim (1+t)^{-1}\cdot|\ln(e+t)|^{\delta-1-\frac{n}{2}}\cdot
\|v\|^2\cdot|\ln(e+t)|^{\frac{n}{2}}
\\
&\lesssim (1+t)^{-1}\cdot|\ln(e+t)|^{\delta-1-\frac{n}{2}}\cdot\Psi_0^2(t),
\\
\int (1+t)\cdot|\ln(e+t)|^{\delta}|\boldsymbol u|^2
&\lesssim
(1+t)^{-1}\cdot|\ln(e+t)|^{\delta-1-\frac{n}{2}}
\cdot \|\partial_x^\alpha\boldsymbol u\|^2\cdot
(1+t)^{2}\cdot|\ln(e+t)|^{1+\frac{n}{2}}
\\
&\lesssim
(1+t)^{-1}\cdot|\ln(e+t)|^{\delta-1-\frac{n}{2}}\cdot\Psi_0^2(t).
\end{align*}
The proof is completed.
$\hfill\Box$

The inhomogeneous terms in the inequalities \eqref{eq-energy-v-1}
and \eqref{eq-energy-u-1} in Lemma \ref{le-timew-1}
are estimated in a similar way as Lemma \ref{le-Q} and Lemma \ref{le-Q2}.

\begin{lemma} \label{le-Q-1}
There holds, for $\lambda=-1$ and $\delta\in(0,\frac{n}{2})$, that
\begin{align*}
&\int (\partial_tQ_1+b(t)\cdot Q_1-\nabla\cdot Q_2)
\cdot(|\ln(e+t)|^{\delta+1}\partial_t v
+\mu_1(1+t)^{-1}\cdot|\ln(e+t)|^{\delta} v)
\\
&\quad+\int (\partial_t Q_2-\nabla Q_1)
\cdot((1+t)^{2}\cdot|\ln(e+t)|^{\delta}\partial_t \boldsymbol u
+\mu_2(1+t)\cdot|\ln(e+t)|^{\delta} \boldsymbol u)
\\
&\lesssim
\partial_tJ_1(t)+
(\Psi_0(t)+\Phi_N(t))\cdot\Psi_1^2(t)
+\Phi_N(t)\cdot\Psi_0^2(t)\cdot(1+t)^{-1}\cdot|\ln(e+t)|^{-\frac{n}{4}},
\end{align*}
provided that $\|v\|_{L^\infty}\le \frac{1}{\gamma-1}$
(which is valid under the a priori assumption $\Phi_N(t)+\Psi_0(t)\le \delta_0$
with a small constant $\delta_0$),
where
$$
J_1(t)\lesssim \|v\|_{L^\infty}\cdot\Phi_1^2(t).
$$
\end{lemma}
{\it\bfseries Proof.}
Noticing that the only difference between this lemma and Lemma \ref{le-Q}
is the time-weights,
we can prove the above decay estimates in the same way as before.
Here we omit the details.
$\hfill\Box$

\begin{lemma} \label{le-Q2-1}
There holds, for integer $k\ge1$, $\lambda=-1$,
$\delta\in(0,\frac{n}{2})$, and $|\alpha|=k$,
\begin{align*}
&\int \partial_x^\alpha(\partial_tQ_1+b(t)\cdot Q_1-\nabla\cdot Q_2)
\cdot(|\ln(e+t)|^{\delta+1}\partial_t \partial_x^\alpha v
+\mu_1(1+t)^{-1}\cdot|\ln(e+t)|^{\delta} \partial_x^\alpha v)
\\
&+\int \partial_x^\alpha(\partial_t Q_2-\nabla Q_1)
\cdot((1+t)^{2}\cdot|\ln(e+t)|^{\delta}\partial_t \partial_x^\alpha\boldsymbol u
+\mu_2(1+t)\cdot|\ln(e+t)|^{\delta} \partial_x^\alpha\boldsymbol u)
\\
&\lesssim
\partial_t J_{k+1}(t)+
(\Psi_0(t)+\Phi_N(t))\cdot\Psi_{k+1}^2(t)\cdot|\ln(e+t)|^{-\frac{\delta}{4}}
+(\Psi_0(t)+\Phi_N(t))\cdot\Psi_k^2(t)\cdot|\ln(e+t)|^{-\frac{\delta}{4}},
\end{align*}
under the assumption that $\|v\|_{L^\infty}\le \frac{1}{\gamma-1}$,
where
$$
J_{k+1}(t)\lesssim \|v\|_{L^\infty}\cdot\Phi_{k+1}^2(t).
$$
\end{lemma}
{\it\bfseries Proof.}
This is proved in a similar way as Lemma \ref{le-Q2}
since the differences only lie in the time-weights.
$\hfill\Box$

The basic energy decay estimates in $\Psi_0(t)$ are deduced by means of
the Green function method.

\begin{lemma} \label{le-Psi0-1}
There hold for $\lambda=-1$ and $n\ge7$ that
\begin{align*}
\|v\|&\lesssim \|(v_0,\boldsymbol u_0)\|_{L^1\cap L^2}\cdot|\ln(e+t)|^{-\frac{n}{4}}
+\Psi_0(t)\Phi_N(t)\cdot|\ln(e+t)|^{-\frac{n}{4}},
\\
\|\boldsymbol u\|&\lesssim
\|(v_0,\boldsymbol u_0)\|_{L^1\cap H^1}\cdot(1+t)^{-1}\cdot|\ln(e+t)|^{-\frac{n}{4}-\frac{1}{2}}
+(\Phi_N(t)+\Psi_0(t))\Phi_N(t)\cdot(1+t)^{-1}\cdot|\ln(e+t)|^{-\frac{n}{4}-\frac{1}{2}}.
\end{align*}
\end{lemma}
{\it\bfseries Proof.}
According to the Duhamel principle \eqref{eq-Duhamel}
and the decay estimates of the Green matrix $\mathcal{G}(t,s)$
in Lemma \ref{le-decay-G-1}, we have
\begin{align*}
\|v(t)\|
\lesssim&
\|\mathcal{G}_{11}(t,0)v_0\|+\|\mathcal{G}_{12}(t,0)\boldsymbol u_0\|
+\int_0^t\|\mathcal{G}_{11}(t,s)Q_1(s)\|ds
+\int_0^t\|\mathcal{G}_{12}(t,s)Q_2(s)\|ds
\\
\lesssim&
\|(v_0,\boldsymbol u_0)\|_{L^1\cap L^2}\cdot|\ln(e+t)|^{-\frac{n}{4}}
+\int_0^t
\Gamma^{\frac{n}{2}}(t,s)\cdot
(\|Q_1(s)\|_{L^1}^l+\|Q_1(s)\|^h)ds
\\
&+\int_0^t
(1+s)^{-1}\cdot
\Gamma^{\frac{n}{2}+1}(t,s)\cdot
(\|Q_2(s)\|_{L^1}^l+\|Q_2(s)\|^h)ds
\\
\lesssim&
\|(v_0,\boldsymbol u_0)\|_{L^1\cap L^2}\cdot|\ln(e+t)|^{-\frac{n}{4}}
\\
&+\Psi_0(t)\Phi_N(t)\int_0^t
\Big(1+\ln\Big(\frac{1+t}{1+s}\Big)\Big)^{-\frac{n}{4}}
\cdot
(1+s)^{-1}\cdot|\ln(e+s)|^{-\frac{n}{4}-\frac{\delta}{2}}ds
\\
&+\Psi_0(t)\Phi_N(t)\int_0^t
(1+s)^{-1}\cdot
\Big(1+\ln\Big(\frac{1+t}{1+s}\Big)\Big)^{-\frac{1}{2}(\frac{n}{2}+1)}
\cdot
|\ln(e+s)|^{-\frac{n}{4}-\frac{\delta+1}{2}}ds
\\
\lesssim&
\|(v_0,\boldsymbol u_0)\|_{L^1\cap L^2}\cdot|\ln(e+t)|^{-\frac{n}{4}}
+\Psi_0(t)\Phi_N(t)\cdot|\ln(e+t)|^{-\frac{n}{4}},
\end{align*}
where we have used Lemma \ref{le-min-1} (note that
\begin{equation}
\begin{cases}
\frac{n}{4}+\frac{\delta}{2}>1,
\quad \frac{n}{4}+\frac{\delta}{2}-1\ge \frac{n}{4},
\\
\frac{n}{4}+\frac{\delta+1}{2}>1,
\quad \frac{n}{4}+\frac{\delta+1}{2}-1\ge \frac{n}{4},
\end{cases}
\end{equation}
for $n\ge5$ and $\delta\in(2,\frac{n}{2})$)
and the following decay estimates on $\|Q(s)\|_{L^1}$ and $\|Q(s)\|$
(we use $D^j:=\partial_x^j$)
\begin{align*}
\|Q_1(s)\|_{L^1}
&\lesssim \|uDv\|_{L^1}+\|vDu\|_{L^1}
\lesssim \|u\|\|Dv\|+\|v\|\|Du\|
\\
&\lesssim
\Psi_0(s)(1+s)^{-1}\cdot|\ln(e+s)|^{-\frac{n}{4}-\frac{1}{2}}
\cdot \Phi_N(s)|\ln(e+s)|^{-\frac{\delta+1}{2}}
\\
&\qquad
+\Psi_0(s)|\ln(e+s)|^{-\frac{n}{4}}
\cdot \Phi_N(s)(1+s)^{-1}\cdot|\ln(e+s)|^{-\frac{\delta}{2}}
\\
&\lesssim
\Psi_0(s)\Phi_N(s)\cdot
(1+s)^{-1}\cdot|\ln(e+s)|^{-\frac{n}{4}-\frac{\delta}{2}},
\\
\|Q_2(s)\|_{L^1}
&\lesssim \|uDu\|_{L^1}+\|vDv\|_{L^1}
\lesssim \|u\|\|Du\|+\|v\|\|Dv\|
\\
&\lesssim
\Psi_0(s)(1+s)^{-1}\cdot|\ln(e+s)|^{-\frac{n}{4}-\frac{1}{2}}
\cdot \Phi_N(s)
(1+s)^{-1}\cdot|\ln(e+s)|^{-\frac{\delta}{2}}
\\
&\qquad
+\Psi_0(s)|\ln(e+s)|^{-\frac{n}{4}}
\cdot \Phi_N(s)|\ln(e+s)|^{-\frac{\delta+1}{2}}
\\
&\lesssim
\Psi_0(s)\Phi_N(s)\cdot|\ln(e+s)|^{-\frac{n}{4}-\frac{\delta+1}{2}}.
\end{align*}
The decay estimates on $\|Q_1\|$ and $\|Q_2\|$ are
at least at the same rates
as $\|Q_1\|_{L^1}$ and $\|Q_2\|_{L^1}$
since the estimates on $\|Dv\|_{L^\infty}$ and $\|D\boldsymbol u\|_{L^\infty}$
decay at the same rates as $\|Dv\|$ and $\|D\boldsymbol u\|$
according to \eqref{eq-Linfty-1}.

We estimate $\|DQ_2\|$ for $n\ge3$ as follows
\begin{align*}
\|DQ_2(s)\|
&\lesssim \|uD^2u\|+\|DuDu\|+\|vD^2v\|+\|DvDv\|
\\
&\lesssim \|u\|_{L^\infty}\|D^2u\|+\|Du\|_{L^\infty}\|Du\|
+\|v\|_{L^\infty}\|D^2v\|+\|Dv\|_{L^\infty}\|Dv\|
\\
&\lesssim
\Phi_N^2(s)(1+s)^{-2}\cdot|\ln(e+s)|^{-\delta}+
\Phi_N^2(s)\cdot|\ln(e+s)|^{-\delta-1}
\\
&\lesssim
\Phi_N^2(s)\cdot|\ln(e+s)|^{-\delta-1},
\end{align*}
according to \eqref{eq-Linfty-1}.
Therefore, we have
\begin{align*}
\|\boldsymbol u(t)\|
\lesssim&
\|\mathcal{G}_{21}(t,0)v_0\|+\|\mathcal{G}_{22}(t,0)\boldsymbol u_0\|
+\int_0^t\|\mathcal{G}_{21}(t,s)Q_1(s)\|ds
+\int_0^t\|\mathcal{G}_{22}(t,s)Q_2(s)\|ds
\\
\lesssim&
\|(v_0,\boldsymbol u_0)\|_{L^1\cap H^1}
\cdot(1+t)^{-1}\cdot|\ln(e+t)|^{-\frac{n}{4}-\frac{1}{2}}
\\
&\quad
+\int_0^t
(1+t)^{-1}\cdot
\Gamma^{\frac{n}{2}+1}(t,s)\cdot
(\|Q_1(s)\|_{L^1}^l+\|Q_1(s)\|^h)ds
\\
&\quad+\int_0^t
(1+t)^{-1}(1+s)^{-1}\cdot
\Gamma^{\frac{n}{2}+2}(t,s)\cdot
(\|Q_2(s)\|_{L^1}^l+\|DQ_2(s)\|^h)ds
\\
\lesssim&
\|(v_0,\boldsymbol u_0)\|_{L^1\cap H^1}\cdot
(1+t)^{-1}\cdot|\ln(e+t)|^{-\frac{n}{4}-\frac{1}{2}}
\\
&\quad
+\Psi_0(t)\Phi_N(t)\int_0^t
(1+t)^{-1}\cdot
\Gamma^{\frac{n}{2}+1}(t,s)\cdot
(1+s)^{-1}\cdot|\ln(e+s)|^{-\frac{n}{4}-\frac{\delta}{2}}ds
\\
&\quad+(\Phi_N(t)+\Psi_0(t))\Phi_N(t)\int_0^t
(1+t)^{-1}(1+s)^{-1}\cdot
\Gamma^{\frac{n}{2}+2}(t,s)
\cdot|\ln(e+s)|^{-\delta-1}ds
\\
\lesssim&
\|(v_0,\boldsymbol u_0)\|_{L^1\cap H^1}\cdot
(1+t)^{-1}\cdot|\ln(e+t)|^{-\frac{n}{4}-\frac{1}{2}}
+(\Phi_N(t)+\Psi_0(t))\Phi_N(t)\cdot
(1+t)^{-1}\cdot|\ln(e+t)|^{-\frac{n}{4}-\frac{1}{2}},
\end{align*}
since
\begin{equation}
\begin{cases}
\frac{n}{4}+\frac{\delta}{2}>1, \quad
&\frac{n}{4}+\frac{\delta}{2}-1\ge \frac{n}{4}+\frac{1}{2},
\\
\delta+1>1, \quad
&\delta+1-1\ge \frac{n}{4}+\frac{1}{2},
\end{cases}
\end{equation}
for $n\ge7$ and $\delta\in(3,\frac{n}{2})$.
The proof is completed.
$\hfill\Box$

\begin{remark}
The restriction of $n\ge7$ comes from the imperfect decay estimate of $\|Q_1\|_1$,
which lays a barrier on the decay estimates of $\|(v,\boldsymbol u)\|$.
From the view of the optimal decay estimates of the linearized hyperbolic system,
it is supposed that both
$\|\boldsymbol u\partial_x v\|$ and $\|v\partial_x \boldsymbol u\|$
decay as $(1+t)^{-1}\cdot|\ln(e+t)|^{-\frac{n}{2}-1}$.
We note that here in the proof of Lemma \ref{le-Q2-1},
the estimate on $\|\boldsymbol u\partial_x v\|$ decays as
$(1+t)^{-1}\cdot|\ln(e+t)|^{-\frac{n}{4}-\frac{1}{2}-\frac{\delta+1}{2}}$,
which is close to the expected optimal decays since $\delta\in(0,\frac{n}{2})$;
while the estimate on $\|v\partial_x \boldsymbol u\|$ decays at
$(1+t)^{-1}\cdot|\ln(e+t)|^{-\frac{n}{4}-\frac{\delta}{2}}$,
which has at least a gap of $|\ln(e+t)|^{-1}$ decay to
the expected optimal decays.
\end{remark}

We combine the above time-weighted iteration scheme and Green function method
to close the decay estimates for $\lambda=-1$.

\begin{proposition} \label{th-com-1}
For $n\ge7$, $N\ge[\frac{n}{2}]+2$ and $\lambda=-1$,
there exists a constant $\varepsilon_0>0$ such that
the solution $(v,\boldsymbol u)$ of the nonlinear system \eqref{eq-vbdu}
corresponding to small initial data
$\|(v_0,\boldsymbol u_0)\|_{L^1\cap H^N}\le \varepsilon_0$
exists globally and satisfies
\begin{equation} \label{eq-com-1}
\begin{cases}
\|v(t)\|\lesssim |\ln(e+t)|^{-\frac{n}{4}}, \\
\|\boldsymbol u(t)\|\lesssim (1+t)^{-1}\cdot|\ln(e+t)|^{-\frac{n}{4}-\frac{1}{2}}.
\end{cases}
\end{equation}
The above decay rates are optimal and consistent with the optimal decay rates
of the linearized hyperbolic system.
\end{proposition}
{\it\bfseries Proof.}
The outline of this proof is similar to Proposition \ref{th-com-0}
for the case of $\lambda\in(-1,0)$.
We claim that the following a priori decay estimate
\begin{equation} \label{eq-apriori-com-1}
\Phi_N(t)+\Psi_0(t)\le \delta_0,
\end{equation}
holds for all the time $t>0$,
under the small energy assumption of initial data
$\|(v_0,\boldsymbol u_0)\|_{L^1\cap H^N}\le \varepsilon_0$,
where $\varepsilon_0$ and $\delta_0$ are positive constants to be determined.
Lemma \ref{le-Psi0-1} tells us that for $n\ge7$
\begin{equation} \label{eq-zPsi0-1}
\Psi_0(T)\le \sup_{t\in(0,T)}\Big\{
|\ln(e+t)|^{\frac{n}{4}}\|v\|,
(1+t)\cdot|\ln(e+t)|^{\frac{n}{4}+\frac{1}{2}}\|\boldsymbol u\|
\Big\}
\lesssim
\varepsilon_0+\delta_0^2.
\end{equation}
According to the time-weighted iteration scheme \eqref{eq-it-0-1} and \eqref{eq-it-1-1}
in Lemma \ref{le-iteration-1} and the estimates of inhomogeneous terms
in Lemma \ref{le-Q-1} and Lemma \ref{le-Q2-1},
we have for integer $0\le k\le N-1$ that
\begin{align*}
&\Phi_1^2(t)+\int_0^t \Psi_1^2(s)ds
\\
&\quad\lesssim
\Phi_1^2(0)+J_1(t)+
\int_0^t (1+s)^{-1}\cdot|\ln(e+s)|^{\delta-1-\frac{n}{2}}\cdot\Psi_0^2(s)ds
\\
&\qquad
+\delta_0\int_0^t \Psi_1^2(s)ds
+\delta_0\int_0^t (1+s)^{-1}\cdot|\ln(e+s)|^{-\frac{n}{4}}
\cdot\Psi_0^2(s)ds,
\\
&\Phi_{k+1}^2(t)+\int_0^t \Psi_{k+1}^2(s)ds
\\
&\quad\lesssim
\Phi_{k+1}^2(0)+J_{k+1}(t)+
\int_0^t \Psi_k^2(s)ds
+\delta_0\int_0^t \Psi_{k+1}^2(s)ds
+\delta_0\int_0^t \Psi_k^2(s)ds,
\end{align*}
where
\begin{align*}
&J_1(t)\lesssim \|v\|_{L^\infty}\cdot\Phi_1^2(t)
\lesssim (\Phi_N(t)+\Psi_0(t))\cdot\Phi_1^2(t)
\lesssim \delta_0\Phi_1^2(t)
\\
&J_{k+1}(t)\lesssim \|v\|_{L^\infty}\cdot\Phi_{k+1}^2(t)
\lesssim (\Phi_N(t)+\Psi_0(t))\Phi_{k+1}^2(t)
\lesssim \delta_0\Phi_{k+1}^2(t).
\end{align*}
We simplify the above inequalities as (note that $\delta_0$ and $\varepsilon_0$ are small)
\begin{align} \label{eq-zPhi0-1}
&\Phi_1^2(t)+\int_0^t \Psi_1^2(s)ds
\lesssim
\varepsilon_0^2+
\Psi_0^2(t)\int_0^t (1+s)^{-1}\cdot|\ln(e+s)|^{\max\{\delta-1-\frac{n}{2},-\frac{n}{4}\}}ds,
\\ \label{eq-zPhi-1}
&\Phi_{k+1}^2(t)+\int_0^t \Psi_{k+1}^2(s)ds
\lesssim
\varepsilon_0^2+
\int_0^t (1+s)^{-1-\lambda}\cdot\Psi_k^2(s)ds.
\end{align}
Multiplying \eqref{eq-zPhi-1} by small positive constants
for $0\le k\le N-1$, summing the resulting inequalities up together with \eqref{eq-zPhi0-1},
we have
\begin{align*}
\sum_{1\le j\le N}\Phi_j^2(t)
\lesssim
\varepsilon_0^2+
\Psi_0^2(t)\int_0^t (1+s)^{-1}\cdot|\ln(e+s)|^{\max\{\delta-1-\frac{n}{2},-\frac{n}{4}\}}ds
\lesssim \varepsilon_0^2+(\varepsilon_0+\delta_0^2)^2,
\end{align*}
according to the estimate \eqref{eq-zPsi0-1}
and $\max\{\delta-1-\frac{n}{2},-\frac{n}{4}\}<-1$ for $n\ge5$.
Therefore,
$$
\Phi_N(t)+\Psi_0(t)
\lesssim \varepsilon_0+\delta_0^2
\le \delta_0,
$$
for positive constants $\varepsilon_0$ and $\delta_0$ small enough.

The optimal property of the decay estimates \eqref{eq-com-1}
follows from the estimates on $\|v\|$ and $\|\boldsymbol u\|$
in Lemma \ref{le-Psi0-1} through a similar procedure as in
the proof of Theorem \ref{th-nonlinear}.
$\hfill\Box$

\vskip2mm
{\it\bfseries Proof of Theorem \ref{th-com-1-in}.}
The critical case of $\lambda=-1$ is proved in Proposition \ref{th-com-1}.
$\hfill\Box$

\begin{appendices}
\section{Time-dependent damped wave equations}

The optimal decay estimates of the time-dependent damped wave equations
\eqref{eq-Pv} and \eqref{eq-Pu} with over-damping $\lambda\in[-1,0)$
are formulated in the similar procedure to the under-damping case $\lambda\in[0,1)$ in \cite{Ji-Mei-1},
but modifications should be made.
Here we sketch the main line of the diagonalization scheme
developed by Wirth \cite{Wirth-JDE06,Wirth-JDE07}
and exact decay behavior of the fundamental solutions.
We would highlight the differences between these two cases.

The Fourier transforms of the time-dependent damped wave equations
\eqref{eq-Pv} and \eqref{eq-Pu} are
\begin{equation} \label{eq-Pv-F}
\begin{cases}
\partial_t^2\hat v+|\xi|^2\hat v+b(t)\partial_t\hat v=0,\\
\hat v(0,\xi)=\hat v_1(\xi), \quad \partial_t\hat v(0,\xi)=\hat v_2(\xi),
\end{cases}
\end{equation}
and
\begin{equation} \label{eq-Pu-F}
\begin{cases}
\partial_t^2\hat u+|\xi|^2\hat u+\partial_t(b(t)\hat u)=0,\\
\hat u(0,\xi)=\hat u_1(\xi), \quad \partial_t\hat u(0,\xi)=\hat u_2(\xi),
\end{cases}
\end{equation}
where $b(t)=\frac{\mu}{(1+t)^\lambda}$ with $\mu>0$ and $\lambda\in[-1,0)$.
The solutions can be represented in the form
\begin{align} \label{eq-Phiv}
\hat v(t,\xi)=\Phi_1^v(t,0,\xi)\hat v_1(\xi)+\Phi_2^v(t,0,\xi)\hat v_2(\xi),
\\ \label{eq-Phiu}
\hat u(t,\xi)=\Phi_1^u(t,0,\xi)\hat u_1(\xi)+\Phi_2^u(t,0,\xi)\hat u_2(\xi),
\end{align}
with Fourier multipliers $\Phi_j^v(t,s,\xi)$ and $\Phi_j^u(t,s,\xi)$, $j=1,2$,
representing the evolution of initial data starting from $s\le t$.
Let
\begin{align*}
\tilde v(t,\xi):=e^{\frac{1}{2}\int_0^tb(\tau)d\tau}\hat v(t,\xi),\\
\tilde u(t,\xi):=e^{\frac{1}{2}\int_0^tb(\tau)d\tau}\hat u(t,\xi).
\end{align*}
Then the equations in \eqref{eq-Pv-F} and \eqref{eq-Pu-F} are transformed into
\begin{align} \label{eq-Pv-tilde}
\partial_t^2\tilde v+\Big(|\xi|^2-\frac{1}{4}b^2(t)-\frac{1}{2}b'(t)\Big)\tilde v=0,
\\ \label{eq-Pu-tilde}
\partial_t^2\tilde u+\Big(|\xi|^2-\frac{1}{4}b^2(t)+\frac{1}{2}b'(t)\Big)\tilde u=0.
\end{align}
For simplicity, we denote
$$
m_v(t,\xi):=|\xi|^2-\frac{1}{4}b^2(t)-\frac{1}{2}b'(t),
\quad m_u(t,\xi):=|\xi|^2-\frac{1}{4}b^2(t)+\frac{1}{2}b'(t).
$$
Note that $|b'(t)|\approx \frac{1}{(1+t)^{1+\lambda}}$
is dominated by $b^2(t)\approx\frac{1}{(1+t)^{2\lambda}}$ as $\lambda\in[-1,0)$.
However, we will show that the difference between $m_v(t,\xi)$ and $m_u(t,\xi)$
leads to a faster decay of
the solution $u(t,x)$ of \eqref{eq-Pu} than
the solution $v(t,x)$ of \eqref{eq-Pv}.

We employ the diagonalization method developed by Wirth \cite{Wirth-JDE06,Wirth-JDE07}
and we pay more attention to the exact asymptotic behavior of different frequencies.
For the sake of simplicity, we only write down the analysis and diagonalization of
the problem \eqref{eq-Pv-tilde} and then we state the difference between the two problems.
The phase-time space $(t,\xi)$ of the problem \eqref{eq-Pv-tilde}
is divided into the following parts:
\begin{align*}
Z_\mathrm{hyp}^v:&=\{(t,\xi);\sqrt{|m_v(t,\xi)|}\ge N_vb(t), m_v(t,\xi)\ge0\}, \\
Z_\mathrm{pd}^v:&=\{(t,\xi);\varepsilon_vb(t)\le
\sqrt{|m_v(t,\xi)|}\le N_vb(t), m_v(t,\xi)\ge0\}, \\
Z_\mathrm{red}^v:&=\{(t,\xi);\sqrt{|m_v(t,\xi)|}\le \varepsilon_vb(t)\}, \\
Z_\mathrm{ell}^v:&=\{(t,\xi);\sqrt{|m_v(t,\xi)|}\ge \varepsilon_vb(t),
m_v(t,\xi)\le0, t\ge t_\mathrm{ell}^v\},
\end{align*}
where $\varepsilon_v>0$ is small and $N_v>\varepsilon_v$, $t_\mathrm{ell}^v>0$.
There remains a bounded part
$\{(t,\xi);\sqrt{|m_v(t,\xi)|}\ge \varepsilon_vb(t),
m_v(t,\xi)\le0, t\in(0,t_\mathrm{ell}^v)\}$
which is of no influence.
The treatment of the zones,
$Z_\mathrm{hyp}^v$, $Z_\mathrm{pd}^v$, $Z_\mathrm{red}^v$, and $Z_\mathrm{ell}^v$
is similar to that in \cite{Wirth-JDE07},
here we present the treatment of the elliptic zone $Z_\mathrm{ell}^v$
in detail since this part will determine the decay rates of solutions.

Note that the elliptic zone $Z_\mathrm{ell}^v$ is expanding.
For any fixed constant $c_0\in(0,\mu/2)$, we would call
\begin{align*}
\text{high~frequencies:}& ~(t,\xi)\in Z_\mathrm{hyp}^v, ~ \text{or~other~mixed~zones},\\
\text{low~frequencies:}& ~(t,\xi) \in Z_\mathrm{ell}^v, ~ |\xi|\le c_0,
\end{align*}
where mixed zones are $Z_\mathrm{pd}^v$, $Z_\mathrm{red}^v$,
and $Z_\mathrm{ell}^v$ with $|\xi|\ge c_0$.

In the elliptic zone $Z_\mathrm{ell}^v$, we let $D_t:=-i\partial_t$
and $V:=(\sqrt{|m_v(t,\xi)|}\tilde v, D_t\tilde v)^\mathrm{T}$,
where $(\cdot)^\mathrm{T}$ is the transpose of a matrix or vector.
Then the equation \eqref{eq-Pv-tilde} is converted into
\begin{equation} \label{eq-V}
D_tV=
\begin{pmatrix}
\frac{D_t\sqrt{|m_v(t,\xi)|}}{\sqrt{|m_v(t,\xi)|}} & \sqrt{|m_v(t,\xi)|}\\
-\sqrt{|m_v(t,\xi)|} &
\end{pmatrix}
V=:A(t,\xi)V.
\end{equation}
Let
$$
M=
\begin{pmatrix}
i & -i\\
1 & 1
\end{pmatrix}
, \quad M^{-1}=\frac{1}{2}
\begin{pmatrix}
-i & 1\\
i & 1
\end{pmatrix}.
$$
Then
\begin{equation} \label{eq-DtA}
\mathcal{D}_t-A(t,\xi)=M(\mathcal{D}_t-\mathcal{D}(t,\xi)-R(t,\xi))M^{-1},
\end{equation}
where
$$
\mathcal{D}_t=
\begin{pmatrix}
D_t & \\
 & D_t
\end{pmatrix}
, \quad \mathcal{D}(t,\xi)=
\begin{pmatrix}
-i\sqrt{|m_v(t,\xi)|} & \\
 & i\sqrt{|m_v(t,\xi)|}
\end{pmatrix}
, \quad R(t,\xi)={\textstyle\frac{D_t\sqrt{|m_v(t,\xi)|}}{2\sqrt{|m_v(t,\xi)|}}}
\begin{pmatrix}
1 & -1\\
-1 & 1
\end{pmatrix}.
$$
The diagonalization method developed by Wirth \cite{Wirth-JDE06,Wirth-JDE07}
is to proceed a step further,
\begin{equation} \label{eq-DtD}
(\mathcal{D}_t-\mathcal{D}(t,\xi)-R(t,\xi))N_1(t,\xi)=
N_1(t,\xi)(\mathcal{D}_t-\mathcal{D}(t,\xi)-F_0(t,\xi)-R_1(t,\xi)),
\end{equation}
with
$$
N^{(1)}(t,\xi)={\textstyle\frac{iD_t\sqrt{|m_v(t,\xi)|}}{2|m_v(t,\xi)|}}
\begin{pmatrix}
 & 1\\
-1 &
\end{pmatrix},
\qquad
F_0(t,\xi)={\textstyle\frac{D_t\sqrt{|m_v(t,\xi)|}}{2\sqrt{|m_v(t,\xi)|}}}
\begin{pmatrix}
1 & \\
 & 1
\end{pmatrix},
$$
and $N_1(t,\xi)=I+N^{(1)}(t,\xi)$,
$$
R_1(t,\xi)=-(I+N^{(1)}(t,\xi))^{-1}(D_tN^{(1)}(t,\xi)-R(t,\xi)N^{(1)}(t,\xi)
+N^{(1)}(t,\xi)F_0(t,\xi)).
$$
Now one can verify that $\|R_1(t,\xi)\|_{\max}\lesssim \frac{1}{(1+t)^{2-\lambda}}$,
whose integral with respect to time over any interval $(s,t)$ is uniformly bounded.

The following asymptotic analysis
will be used to show the optimal decay rates
of the solutions $\hat v(t,\xi)$ and $\hat u(t,\xi)$
for equations \eqref{eq-Pv-F} and \eqref{eq-Pu-F}.
Note that for the over-damping case $\lambda\in[-1,0)$, we have $b'(t)\ge0$,
which is slightly different from the under-damping case $\lambda\in[0,1)$.

\begin{lemma} \label{le-decay-est}
For $(t,\xi)\in Z_\mathrm{ell}^v$, there holds
\begin{equation} \label{eq-est-v}
\begin{cases}
\displaystyle
\sqrt{|m_v(t,\xi)|}+
{\textstyle\frac{\partial_t\sqrt{|m_v(t,\xi)|}}{2\sqrt{|m_v(t,\xi)|}}}
-\frac{b(t)}{2}
\le -|\xi|^2\frac{C_1}{b(t)}+\frac{b'(t)}{b(t)}+|r_v(t,\xi)|,
\\[3mm] \displaystyle
\sqrt{|m_v(t,\xi)|}+
{\textstyle\frac{\partial_t\sqrt{|m_v(t,\xi)|}}{2\sqrt{|m_v(t,\xi)|}}}
-\frac{b(t)}{2}
\ge -|\xi|^2\frac{C_2}{b(t)}+\frac{b'(t)}{b(t)}-|r_v(t,\xi)|,
\end{cases}
\end{equation}
and for $(t,\xi)\in Z_\mathrm{ell}^u$
(the definition of zones in the phase-time space corresponding to $\tilde u$
is completely similar to that of $\tilde v$), there holds
\begin{equation} \label{eq-est-u}
\begin{cases}
\displaystyle
\sqrt{|m_u(t,\xi)|}+
{\textstyle\frac{\partial_t\sqrt{|m_u(t,\xi)|}}{2\sqrt{|m_u(t,\xi)|}}}
-\frac{b(t)}{2}
\le -|\xi|^2\frac{C_3}{b(t)}+|r_u(t,\xi)|,
\\[3mm] \displaystyle
\sqrt{|m_u(t,\xi)|}+
{\textstyle\frac{\partial_t\sqrt{|m_u(t,\xi)|}}{2\sqrt{|m_u(t,\xi)|}}}
-\frac{b(t)}{2}
\ge -|\xi|^2\frac{C_4}{b(t)}-|r_u(t,\xi)|,
\end{cases}
\end{equation}
where $|r_v(t,\xi)|\lesssim \frac{1}{(1+t)^{2-\lambda}}$ and
$|r_u(t,\xi)|\lesssim \frac{1}{(1+t)^{2-\lambda}}$
such that the integrals of $|r_v(t,\xi)|$ and $|r_u(t,\xi)|$
with respect to time are uniformly bounded.
\end{lemma}
{\it \bfseries Proof.}
Recall that
$$
m_v(t,\xi):=|\xi|^2-\frac{1}{4}b^2(t)-\frac{1}{2}b'(t),
\quad m_u(t,\xi):=|\xi|^2-\frac{1}{4}b^2(t)+\frac{1}{2}b'(t),
$$
and in the elliptic zone $Z_\mathrm{ell}^v$ or $Z_\mathrm{ell}^u$,
$m_v(t,\xi)<0$ and $\sqrt{|m_v(t,\xi)|}\ge \varepsilon_v b(t)$,
or $m_u(t,\xi)<0$ and $\sqrt{|m_u(t,\xi)|}\ge \varepsilon_u b(t)$, respectively.
Then we have
$|m_v(t,\xi)|=\frac{1}{4}b^2(t)+\frac{1}{2}b'(t)-|\xi|^2\ge \varepsilon_v^2 b^2(t)$,
$|m_v(t,\xi)|\le \frac{1}{4}b^2(t)+\frac{1}{2}b'(t)\le\frac{1}{2}b^2(t)$
since $|b'(t)|$ is dominated by $b^2(t)$ and the elliptic zone is defined
within $t\ge t_\mathrm{ell}^v$ which can be chosen large.
Therefore,
\begin{align*}
&\sqrt{|m_v(t,\xi)|}+
{\textstyle\frac{\partial_t\sqrt{|m_v(t,\xi)|}}{2\sqrt{|m_v(t,\xi)|}}}
-\frac{b(t)}{2}
\\
=&\frac{|m_v(t,\xi)|^2-\frac{1}{4}b^2(t)}{\sqrt{|m_v(t,\xi)|}+\frac{b(t)}{2}}
+\frac{\frac{1}{2}b(t)b'(t)+\frac{1}{2}b''(t)}{4(\frac{1}{4}b^2(t)+\frac{1}{2}b'(t)-|\xi|^2)}
\\
=&\frac{-|\xi|^2}{\sqrt{|m_v(t,\xi)|}+\frac{b(t)}{2}}
+\frac{\frac{1}{2}b'(t)}{\sqrt{|m_v(t,\xi)|}+\frac{b(t)}{2}}
+\frac{\frac{1}{2}b(t)b'(t)}{4(\frac{1}{4}b^2(t)+\frac{1}{2}b'(t)-|\xi|^2)}
+\frac{\frac{1}{2}b''(t)}{4(\frac{1}{4}b^2(t)+\frac{1}{2}b'(t)-|\xi|^2)},
\end{align*}
and
\begin{align*}
&\Big(\sqrt{|m_v(t,\xi)|}+
{\textstyle\frac{\partial_t\sqrt{|m_v(t,\xi)|}}{2\sqrt{|m_v(t,\xi)|}}}
-\frac{b(t)}{2}\Big)
-\Big(\frac{-|\xi|^2}{\sqrt{|m_v(t,\xi)|}+\frac{b(t)}{2}}
+\frac{b'(t)}{b(t)}
\Big)
\\
=&\Big(
\frac{\frac{1}{2}b'(t)}{\sqrt{|m_v(t,\xi)|}+\frac{b(t)}{2}}
-\frac{\frac{1}{2}b'(t)}{b(t)}\Big)
+\Big(
\frac{\frac{1}{2}b(t)b'(t)}{4(\frac{1}{4}b^2(t)+\frac{1}{2}b'(t)-|\xi|^2)}
-\frac{\frac{1}{2}b(t)b'(t)}{b^2(t)}
\Big)
+\frac{\frac{1}{2}b''(t)}{4(\frac{1}{4}b^2(t)+\frac{1}{2}b'(t)-|\xi|^2)}
\\
=&:\bar r_v(t,\xi).
\end{align*}
We estimate $\bar r_v(t,\xi)$ as follows
\begin{align*}
|\bar r_v(t,\xi)|
\le& \Big|
\frac{\frac{1}{2}b'(t)}{\sqrt{|m_v(t,\xi)|}+\frac{b(t)}{2}}
-\frac{\frac{1}{2}b'(t)}{\frac{b(t)}{2}+\frac{b(t)}{2}}
\Big|
\\
&+\Big|
\frac{\frac{1}{2}b(t)b'(t)}{4(\frac{1}{4}b^2(t)+\frac{1}{2}b'(t)-|\xi|^2)}
-\frac{\frac{1}{2}b(t)b'(t)}{b^2(t)}
\Big|
+\Big|
\frac{\frac{1}{2}b''(t)}{4(\frac{1}{4}b^2(t)+\frac{1}{2}b'(t)-|\xi|^2)}
\Big|
\\
\lesssim&
\frac{\frac{1}{2}b'(t)\big|\frac{1}{2}b(t)-\sqrt{|m_v(t,\xi)|}\big|}%
{(\sqrt{|m_v(t,\xi)|}+\frac{b(t)}{2})b(t)}
+\frac{\frac{1}{2}b(t)b'(t)|-2b'(t)+4|\xi|^2|}%
{4(\frac{1}{4}b^2(t)+\frac{1}{2}b'(t)-|\xi|^2)b^2(t)}
+\frac{|b''(t)|}{b^2(t)}
\\
\lesssim&
\frac{b'(t)}{b^2(t)}
\cdot \frac{|-\frac{1}{2}b'(t)+|\xi|^2|}{\frac{1}{2}b(t)+\sqrt{|m_v(t,\xi)|}}
+\frac{|b'(t)|^2}{b^3(t)}
+\frac{b'(t)}{b^2(t)}\cdot\frac{|\xi|^2}{b(t)}
+\frac{|b''(t)|}{b^2(t)}
\\
\lesssim&
\frac{b'(t)}{b^2(t)}\cdot\frac{|\xi|^2}{b(t)}
+\frac{|b'(t)|^2}{b^3(t)}
+\frac{|b''(t)|}{b^2(t)}.
\end{align*}
By noticing that $\frac{|b'(t)|}{b^2(t)}\lesssim \frac{1}{(1+t)^{1-\lambda}}$,
which tends to zero as $t\to\infty$,
we find that $\bar r_v(t,\xi)$ can be split into
$$\bar r_v(t,\xi)=|\xi|^2\frac{1}{b(t)}\cdot \omega(t,\xi)+r_v(t,\xi),$$
with
$$|r_v(t,\xi)|\lesssim \frac{|b'(t)|^2}{b^3(t)}+\frac{|b''(t)|}{b^2(t)}
\lesssim\frac{1}{(1+t)^{2-\lambda}},$$
and
$$\frac{-|\xi|^2}{\sqrt{|m_v(t,\xi)|}+\frac{b(t)}{2}}
+|\xi|^2\frac{1}{b(t)}\cdot \omega(t,\xi)
\approx -|\xi|^2\frac{1}{b(t)}$$
since $|\omega(t,\xi)|\lesssim \frac{1}{(1+t)^{1-\lambda}}$
and we can choose $t_\mathrm{ell}^v$ large enough
such that $|\omega(t,\xi)|\le 1/4$.
The proof of \eqref{eq-est-u} follows similarly.
$\hfill\Box$

According to the asymptotic analysis of the frequencies, we can
formulate the following estimates.
We note that for the over-damping case $\lambda\in[-1,0)$,
the elliptic zone $Z_\mathrm{ell}^v$ is expanding,
which differs from the shrinking elliptic zone for under-damping case.

\begin{lemma} \label{le-Phi}
The multiplies $\Phi_j^v(t,s,\xi)$ and $\Phi_j^u(t,s,\xi)$, $j=1,2$,
in the equations \eqref{eq-Phiv} and \eqref{eq-Phiu} have the following estimates:
there exist $c_0>0$, $\varepsilon\in(0,1/2)$, $C>0$, and $T_0\ge0$
(only depending on $\mu$ and $\lambda$) such that

(i) For $(t,\xi)\in Z_\mathrm{ell}^v$, $0\le s\le t$, and $|\xi|\le c_0$, there hold
\begin{equation} \label{eq-Phiv-low}
|\Phi_1^v(t,s,\xi)|\lesssim e^{-C|\xi|^2\int_s^t\frac{1}{b(\tau)}d\tau},
\quad
|\Phi_2^v(t,s,\xi)|\lesssim \frac{1}{b(s)}\cdot e^{-C|\xi|^2\int_s^t\frac{1}{b(\tau)}d\tau};
\end{equation}
for $(t,\xi)\in Z_\mathrm{hyp}^v$ and $0\le s\le t$, there holds
$$
|\Phi_1^v(t,s,\xi)|+|\xi||\Phi_2^v(t,s,\xi)|
\lesssim e^{-(\frac{1}{2}-\varepsilon)\int_s^tb(\tau)d\tau};
$$
and for $(t,\xi)\not\in Z_\mathrm{hyp}^v$ with $0\le s\le t$ and $|\xi|\ge c_0$,
there hold
\begin{align*}
|\Phi_1^v(t,s,\xi)|&\lesssim e^{-C|\xi|^2\int_{\max\{s,t_\xi^v\}}^t\frac{1}{b(\tau)}d\tau
-(\frac{1}{2}-\varepsilon)\int_s^{\max\{s,t_\xi^v\}}b(\tau)d\tau},
\\
|\Phi_2^v(t,s,\xi)|&\lesssim {\textstyle\frac{1}{b(\max\{s,t_\xi^v\})}}\cdot
e^{-C|\xi|^2\int_{\max\{s,t_\xi^v\}}^t\frac{1}{b(\tau)}d\tau
-(\frac{1}{2}-\varepsilon)\int_s^{\max\{s,t_\xi^v\}}b(\tau)d\tau},
\end{align*}
where $t_\xi^v:=\sup\{t;(t,\xi)\in Z_\mathrm{hyp}^v\}$.

(ii) For $(t,\xi)\in Z_\mathrm{ell}^u$, $0\le s\le t$, and $|\xi|\le c_0$, there hold
\begin{equation} \label{eq-Phiu-low}
|\Phi_1^u(t,s,\xi)|\lesssim \frac{b(s)}{b(t)}\cdot
e^{-C|\xi|^2\int_s^t\frac{1}{b(\tau)}d\tau},
\quad
|\Phi_2^u(t,s,\xi)|\lesssim \frac{1}{b(t)}\cdot
e^{-C|\xi|^2\int_s^t\frac{1}{b(\tau)}d\tau};
\end{equation}
for $(t,\xi)\in Z_\mathrm{hyp}^u$ and $0\le s\le t$, there holds
$$
|\Phi_1^u(t,s,\xi)|+|\xi||\Phi_2^u(t,s,\xi)|
\lesssim e^{-(\frac{1}{2}-\varepsilon)\int_s^tb(\tau)d\tau};
$$
and for $(t,\xi)\not\in Z_\mathrm{hyp}^u$ with $0\le s\le t$ and $|\xi|\ge c_0$,
there hold
\begin{align*}
|\Phi_1^u(t,s,\xi)|&\lesssim
{\textstyle\frac{b(\max\{s,t_\xi^u\})}{b(t)}}\cdot
e^{-C|\xi|^2\int_{\max\{s,t_\xi^u\}}^t\frac{1}{b(\tau)}d\tau
-(\frac{1}{2}-\varepsilon)\int_s^{\max\{s,t_\xi^u\}}b(\tau)d\tau},
\\
|\Phi_2^u(t,s,\xi)|&\lesssim
{\textstyle\frac{1}{b(t)}}\cdot
e^{-C|\xi|^2\int_{\max\{s,t_\xi^u\}}^t\frac{1}{b(\tau)}d\tau
-(\frac{1}{2}-\varepsilon)\int_s^{\max\{s,t_\xi^u\}}b(\tau)d\tau},
\end{align*}
where $t_\xi^u:=\sup\{t;(t,\xi)\in Z_\mathrm{hyp}^u\}$.

(iii) For $(t,\xi)\in Z_\mathrm{ell}^v$, $T_0\le s\le t$, and $|\xi|\le c_0$,
the estimate \eqref{eq-Phiv-low} is optimal:
\begin{equation} \label{eq-Phiv-low-bk}
|\Phi_1^v(t,s,\xi)|\gtrsim
e^{-C |\xi|^2\int_s^t\frac{1}{b(\tau)}d\tau},
\quad
|\Phi_2^v(t,s,\xi)|\gtrsim
\frac{1}{b(s)}\cdot
e^{-C |\xi|^2\int_s^t\frac{1}{b(\tau)}d\tau},
\end{equation}
with another universal constant $C>0$.

(iv) For $(t,\xi)\in Z_\mathrm{ell}^u$, $T_0\le s\le t$, and $|\xi|\le c_0$,
the estimate \eqref{eq-Phiu-low} is optimal:
\begin{equation} \label{eq-phiu-low-bk}
|\Phi_1^u(t,s,\xi)|\gtrsim \frac{b(s)}{b(t)}\cdot
e^{-C|\xi|^2\int_s^t\frac{1}{b(\tau)}d\tau},
\quad
|\Phi_2^u(t,s,\xi)|\gtrsim \frac{1}{b(t)}\cdot
e^{-C|\xi|^2\int_s^t\frac{1}{b(\tau)}d\tau},
\end{equation}
with another universal constant $C>0$.
\end{lemma}
{\it \bfseries Proof.}
The estimates (i) with $s=0$ was proved by Wirth in Theorem 17 of \cite{Wirth-JDE07}.
Here we focus on the exact decay estimates of $\Phi_j^v(t,s,\xi)$ with $0\le s\le t$
for the application to nonlinear system \eqref{eq-vbdu}
since $\Phi_j^v(t,s,\xi)$ behaves different from
$\Phi_j^v(t-s,0,\xi)$.
The above estimates are proved in a similar way as
Lemma 2.3 in \cite{Ji-Mei-1} for the under-damping case.
Noticing that the elliptic zone $Z_\mathrm{ell}^v$ is expanding with respect to time,
for the mixed part $(t,\xi)\not\in Z_\mathrm{hyp}^v$ with $0\le s\le t$ and $|\xi|\ge c_0$,
we apply the estimates \eqref{eq-Phiv-low} to $\Phi_j^v(t,s,\xi)$ if $s\ge t_\xi^v$.
$\hfill\Box$

The above frequency analysis is used to show the following optimal decay estimates
of the wave equations \eqref{eq-Pv} and \eqref{eq-Pu}.
Note that the time decay function $\Gamma(t,s)$ is defined in \eqref{eq-Gamma}.

\begin{theorem}[Optimal decay rates of linear wave equations] \label{th-wave}
Let $v(t,x)$ and $u(t,x)$ be the solutions of the
Cauchy problems \eqref{eq-Pv} and \eqref{eq-Pu}
corresponding to initial data $(v(s,x),\partial_tv(s,x))$
and $(u(s,x),\partial_tu(s,x))$ starting from the time $s$, respectively.
For $q\in[2,\infty]$, $1\le p,r\le 2$ and $\lambda\in[-1,0)$, we have
\begin{align} \nonumber
\|\partial_x^\alpha v\|_{L^q}\lesssim &
\Gamma^{\gamma_{p,q}+|\alpha|}(t,s)\cdot
\Big(\big\|v(s,\cdot)\big\|_{L^{p}}^l
+\big\|\partial_x^{|\alpha|+\omega_{r,q}}v(s,\cdot)\big\|_{L^{r}}^h\Big)
\\ \label{eq-optimal-v}
&+(1+s)^\lambda\cdot\Gamma^{\gamma_{p,q}+|\alpha|}(t,s)\cdot
\Big(\big\|\partial_tv(s,\cdot)\big\|_{L^{p}}^l
+\big\|\partial_x^{|\alpha|-1+\omega_{r,q}}\partial_tv(s,\cdot)\big\|_{L^{r}}^h\Big),
\end{align}
and
\begin{align} \nonumber
\|\partial_x^\alpha u\|_{L^q}\lesssim &
\Big(\frac{1+t}{1+s}\Big)^\lambda\cdot
\Gamma^{\gamma_{p,q}+|\alpha|}(t,s)\cdot
\Big(\big\|u(s,\cdot)\big\|_{L^{p}}^l
+\big\|\partial_x^{|\alpha|+\omega_{r,q}}u(s,\cdot)\big\|_{L^{r}}^h\Big)
\\ \label{eq-optimal-u}
&+(1+t)^\lambda\cdot\Gamma^{\gamma_{p,q}+|\alpha|}(t,s)\cdot
\Big(\big\|\partial_tu(s,\cdot)\big\|_{L^{p}}^l
+\big\|\partial_x^{|\alpha|-1+\omega_{r,q}}\partial_tu(s,\cdot)\big\|_{L^{r}}^h\Big),
\end{align}
where $\gamma_{p,q}:=n(1/{p}-1/q)$,
and $\omega_{r,q}>\gamma_{r,q}$ for $(r,q)\ne(2,2)$ and $\omega_{2,2}=0$.

The decay estimates \eqref{eq-optimal-v} and \eqref{eq-optimal-u} are optimal
for all $t\ge s\ge 0$.
Moreover, there exists a $T_0\ge0$ such that the
decay estimates \eqref{eq-optimal-v} and \eqref{eq-optimal-u} are
element-by-element optimal
for all $\frac{t}{2}\ge s\ge T_0$.
\end{theorem}

\begin{corollary}
Let $v(t,x)$ and $u(t,x)$ be the solutions of the
Cauchy problems \eqref{eq-Pv} and \eqref{eq-Pu}
corresponding to initial data $(v(0,x),\partial_tv(0,x))$
and $(u(0,x),\partial_tu(0,x))$ respectively.

(i) For $q\in[2,\infty]$, $1\le p,r\le 2$ and $\lambda\in(-1,0)$, we have
\begin{align} \label{eq-optimal-v-st}
\|\partial_x^\alpha v\|_{L^q}\lesssim &
(1+t)^{-\frac{1+\lambda}{2}(\gamma_{p,q}+|\alpha|)}\cdot
\Big(\big\|(v(0,\cdot),\partial_tv(0,\cdot))\big\|_{L^{p}}^l
+\big\|(\partial_x^{|\alpha|+\omega_{r,q}}v(0,\cdot),
\partial_x^{|\alpha|-1+\omega_{r,q}}\partial_tv(0,\cdot))\big\|_{L^{r}}^h\Big),
\end{align}
and
\begin{align} \label{eq-optimal-u-st}
\|\partial_x^\alpha u\|_{L^q}\lesssim &
(1+t)^{-\frac{1+\lambda}{2}(\gamma_{p,q}+|\alpha|)+\lambda}\cdot
\Big(\big\|(u(0,\cdot),\partial_tu(0,\cdot))\big\|_{L^{p}}^l
+\big\|(\partial_x^{|\alpha|+\omega_{r,q}}u(0,\cdot),
\partial_x^{|\alpha|-1+\omega_{r,q}}\partial_tu(0,\cdot))\big\|_{L^{r}}^h\Big),
\end{align}
where $\gamma_{p,q}:=n(1/{p}-1/q)$,
and $\omega_{r,q}>\gamma_{r,q}$ for $(r,q)\ne(2,2)$ and $\omega_{2,2}=0$.

(ii) For $q\in[2,\infty]$, $1\le p,r\le 2$ and $\lambda=-1$, we have
\begin{align} \nonumber
\|\partial_x^\alpha v\|_{L^q}\lesssim &
|\ln(e+t)|^{-\frac{1}{2}(\gamma_{p,q}+|\alpha|)}
\\ \label{eq-optimal-v-st-1}
&\qquad\cdot
\Big(\big\|(v(0,\cdot),\partial_tv(0,\cdot))\big\|_{L^{p}}^l
+\big\|(\partial_x^{|\alpha|+\omega_{r,q}}v(0,\cdot),
\partial_x^{|\alpha|-1+\omega_{r,q}}\partial_tv(0,\cdot))\big\|_{L^{r}}^h\Big),
\end{align}
and
\begin{align} \nonumber
\|\partial_x^\alpha u\|_{L^q}\lesssim &
(1+t)^{-1}\cdot
|\ln(e+t)|^{-\frac{1}{2}(\gamma_{p,q}+|\alpha|)}
\\ \label{eq-optimal-u-st-1}
&\qquad\cdot
\Big(\big\|(u(0,\cdot),\partial_tu(0,\cdot))\big\|_{L^{p}}^l
+\big\|(\partial_x^{|\alpha|+\omega_{r,q}}u(0,\cdot),
\partial_x^{|\alpha|-1+\omega_{r,q}}\partial_tu(0,\cdot))\big\|_{L^{r}}^h\Big),
\end{align}
where $\gamma_{p,q}:=n(1/{p}-1/q)$,
and $\omega_{r,q}>\gamma_{r,q}$ for $(r,q)\ne(2,2)$ and $\omega_{2,2}=0$.

The decay estimates \eqref{eq-optimal-v-st}, \eqref{eq-optimal-u-st},
\eqref{eq-optimal-v-st-1} and \eqref{eq-optimal-u-st-1} are optimal.
\end{corollary}

\begin{remark}
The decay estimate \eqref{eq-optimal-v-st} for $s=0$
was first proved by Wirth \cite{Wirth-JDE07}
by developing a perfect diagonalization method.
For the application to nonlinear systems, we need to consider the
evolution of initial data starting from any $s\ge0$ to $t\ge s$
since the damping is time-dependent.
One of the main difficulties caused by the time-dependent damping is that
the evolution of the initial data starting from $s\ge0$ to $t\ge s$
is completely different form that starting from $0$ to $t-s$,
as can be seen from the estimates \eqref{eq-optimal-v} and \eqref{eq-optimal-u}.
\end{remark}

\begin{remark}
The two Cauchy problems \eqref{eq-Pv} and \eqref{eq-Pu}
decay with different rates.
We note that the function
$$
\varphi(t,x):=
\begin{cases}
\displaystyle
\frac{1}{(1+t)^{\frac{1+\lambda}{2}n}}%
e^{-\frac{\mu(1+\lambda)|x|^2}{4(1+t)^{1+\lambda}}}, \quad &\lambda\in(-1,0),
\\[4mm]
\displaystyle
\frac{1}{|\ln(e+t)|^{\frac{n}{2}}}%
e^{-\frac{\mu|x|^2}{4\ln(e+t)}}, \quad &\lambda=-1,
\end{cases}
$$
which satisfies $\frac{\mu}{(1+t)^\lambda} \partial_t \varphi=\Delta \varphi$,
is an asymptotic profile of \eqref{eq-Pv},
while $\psi(t,x):=\varphi(t,x)/(\frac{\mu}{(1+t)^\lambda})$,
which satisfies $\partial_t(\frac{\mu}{(1+t)^\lambda}\psi)=\Delta \psi$,
is a good asymptotic profile of \eqref{eq-Pu},
and $\psi(t,x)$ decays faster than $\varphi(t,x)$.
\end{remark}

{\it\bfseries Proof of Theorem \ref{th-wave}.}
The estimate \eqref{eq-optimal-v-st} for $s=0$ was proved by Wirth \cite{Wirth-JDE07}.
Here we focus on the influence of $s$ and show that $u(t,x)$
decays optimally faster than $v(t,x)$.
Note that $(\frac{1+t}{1+s})^\lambda$ decays to zero since $\lambda\in[-1,0)$
and
$$\int_s^t\frac{1}{b(\tau)}d\tau=\ln\Big(\frac{1+t}{1+s}\Big)$$
for the critical case $\lambda=-1$.
The results are proved through the same procedure as Proposition 2.1 in \cite{Ji-Mei-1}
according to the optimal decay estimates on the Fourier multiplies $\Phi_j^v(t,s,\xi)$
and $\Phi_j^u(t,s,\xi)$ in Lemma \ref{le-Phi}.
$\hfill\Box$

\section{Time-dependent damped linear system}

We next show the optimal decay estimates of the linear hyperbolic system \eqref{eq-vu}.

\begin{theorem}[Optimal decay rates of linear hyperbolic system] \label{th-linear-A}
Let $(v(t,x),u(t,x))$ be the solution of the linear hyperbolic system \eqref{eq-vu}
(the third equation of $\boldsymbol w(t,x)$ is neglected as it decays super-exponentially)
corresponding to the initial data $(v(s,x),u(s,x))$
starting from time $s$.
There exists a universal constant $T_0\ge0$ such that
for $q\in[2,\infty]$, $1\le p,r\le 2$, $\lambda\in[-1,0)$, and $t\ge s\ge T_0$, we have
\begin{align} \nonumber
\|\partial_x^\alpha v\|_{L^q}\lesssim &
\Gamma^{\gamma_{p,q}+|\alpha|}(t,s)\cdot
\Big(\big\|v(s,\cdot)\big\|_{L^{p}}^l
+\big\|\partial_x^{|\alpha|+\omega_{r,q}}v(s,\cdot)\big\|_{L^{r}}^h\Big)
\\ \label{eq-linear-v}
&+(1+s)^\lambda\cdot
\Gamma^{\gamma_{p,q}+|\alpha|+1}(t,s)\cdot
\Big(\big\| u(s,\cdot)\big\|_{L^{p}}^l
+\big\|
\partial_x^{|\alpha|+\omega_{r,q}} u(s,\cdot)\big\|_{L^{r}}^h\Big),
\end{align}
and
\begin{align} \nonumber
\|\partial_x^\alpha u\|_{L^q}\lesssim &
\Big(\frac{1+t}{1+s}\Big)^\lambda\cdot
\Gamma^{\gamma_{p,q}+|\alpha|}(t,s)\cdot
\Big(\big\|u(s,\cdot)\big\|_{L^{p}}^l
+\big\|\partial_x^{|\alpha|+\omega_{r,q}}u(s,\cdot)\big\|_{L^{r}}^h\Big)
\\ \label{eq-linear-u}
&+(1+t)^\lambda\cdot
\Gamma^{\gamma_{p,q}+|\alpha|+1}(t,s)\cdot
\Big(\big\|v(s,\cdot)\big\|_{L^{p}}^l
+\big\|\partial_x^{|\alpha|+\omega_{r,q}} v(s,\cdot)\big\|_{L^{r}}^h\Big),
\end{align}
where $\gamma_{p,q}:=n(1/{p}-1/q)$,
and $\omega_{r,q}>\gamma_{r,q}$ for $(r,q)\ne(2,2)$ and $\omega_{2,2}=0$.

Moreover, $u(t,x)$ decays faster than \eqref{eq-linear-u}
if we assume one-order higher regularity as follows,
\begin{align} \nonumber
\|\partial_x^\alpha u\|_{L^q}\lesssim &
(1+t)^\lambda\cdot
\Gamma^{\gamma_{p,q}+|\alpha|+1}(t,s)\cdot
\Big(\big\|v(s,\cdot)\big\|_{L^{p}}^l
+\big\|\partial_x^{|\alpha|+1+\omega_{r,q}}v(s,\cdot)\big\|_{L^{r}}^h\Big)
\\ \label{eq-linear-u-opt}
&+(1+t)^\lambda(1+s)^\lambda\cdot
\Gamma^{\gamma_{p,q}+|\alpha|+2}(t,s)\cdot
\Big(\big\| u(s,\cdot)\big\|_{L^{p}}^l
+\big\|\partial_x^{|\alpha|+1+\omega_{r,q}} u(s,\cdot)\big\|_{L^{r}}^h\Big).
\end{align}

The decay estimate \eqref{eq-linear-u} is improved by cancellation without
one-order higher regularity as follows
\begin{align} \nonumber
\|\partial_x^\alpha u(t,\cdot)\|_{L^q}
\lesssim &
(1+t)^\lambda\cdot
\Gamma^{\gamma_{p,q}+|\alpha|+1}(t,s)
\cdot
\Big(\big\|v(s,\cdot)\big\|_{L^{p}}^l
+\big\|\partial_x^{|\alpha|+\omega_{r,q}}v(s,\cdot)\big\|_{L^{r}}^h\Big)
\\ \nonumber
&+(1+t)^\lambda(1+s)^\lambda\cdot
\Gamma^{\gamma_{p,q}+|\alpha|+2}(t,s)
\cdot
\Big(\big\| u(s,\cdot)\big\|_{L^{p}}^l
+\big\|\partial_x^{|\alpha|+\omega_{r,q}} u(s,\cdot)\big\|_{L^{r}}^h\Big)
\\ \nonumber
&+\Big(\frac{1+t}{1+s}\Big)^\lambda\cdot
\Gamma^{\gamma_{p,q}+|\alpha|}(t,s)
\cdot
\Big(\frac{1}{(1+s)^{1-\lambda}}
+e^{-\varepsilon_u((1+t)^{1-\lambda}-(1+s)^{1-\lambda})}\Big)
\\ \label{eq-linear-u-can}
&\qquad\cdot
\Big(\big\| u(s,\cdot)\big\|_{L^{p}}^l
+\big\|\partial_x^{|\alpha|+\omega_{r,q}} u(s,\cdot)\big\|_{L^{r}}^h\Big),
\end{align}
where $\varepsilon_u>0$ is a constant.

The decay estimate \eqref{eq-linear-v} is element-by-element optimal for all
$\frac{t}{2}\ge s\ge T_0$;
the decay estimate \eqref{eq-linear-u-opt} is optimal with respect to $v(s,x)$
for all $\frac{t}{2}\ge s\ge T_0$;
the decay estimates \eqref{eq-linear-v} and \eqref{eq-linear-u-opt} are optimal
for all $t\ge s\ge0$ such that
\begin{align*}
\|\partial_x^\alpha v\|_{L^q}\approx &
\Gamma^{\gamma_{p,q}+|\alpha|}(t,s)\cdot
\Big(\big\|v(s,\cdot)\big\|_{L^{p}}^l
+\big\|\partial_x^{|\alpha|+\omega_{r,q}}v(s,\cdot)\big\|_{L^{r}}^h
+\big\| u(s,\cdot)\big\|_{L^{p}}^l
+\big\|
\partial_x^{|\alpha|+\omega_{r,q}} u(s,\cdot)\big\|_{L^{r}}^h\Big),
\end{align*}
and
\begin{align*} \nonumber
\|\partial_x^\alpha u\|_{L^q}\approx &
(1+t)^\lambda\cdot
\Gamma^{\gamma_{p,q}+|\alpha|+1}(t,s)
\\
&\qquad\cdot
\Big(\big\|v(s,\cdot)\big\|_{L^{p}}^l
+\big\|\partial_x^{|\alpha|+1+\omega_{r,q}}v(s,\cdot)\big\|_{L^{r}}^h
+\big\| u(s,\cdot)\big\|_{L^{p}}^l
+\big\|\partial_x^{|\alpha|+1+\omega_{r,q}} u(s,\cdot)\big\|_{L^{r}}^h\Big),
\end{align*}
for some nontrivial initial data.
\end{theorem}

Theorem \ref{th-linear-A} implies the optimal decay estimates of the
Green matrix $\mathcal{G}(t,s)$ in \eqref{eq-Duhamel}.

\begin{theorem} \label{th-linear}
For $q\in[2,\infty]$, $1\le p,r\le 2$,
$t\ge s\ge T_0$ ($T_0$ is the universal constant in Theorem \ref{th-linear-A}),
and $\lambda\in[-1,0)$, we have
\begin{align*}
&\|\partial_x^\alpha \mathcal{G}_{11}(t,s)\phi(x)\|_{L^q}\lesssim
\Gamma^{\gamma_{p,q}+|\alpha|}(t,s)\cdot
(\|\phi\|_{L^{p}}^l+\|\partial_x^{|\alpha|+\omega_{r,q}}\phi\|_{L^{r}}^h),
\\
&\|\partial_x^\alpha \mathcal{G}_{12}(t,s)\phi(x)\|_{L^q}\lesssim
(1+s)^\lambda\cdot
\Gamma^{\gamma_{p,q}+|\alpha|+1}(t,s)\cdot
(\|\phi\|_{L^{p}}^l+\|\partial_x^{|\alpha|+\omega_{r,q}}\phi\|_{L^{r}}^h),
\\
&\|\partial_x^\alpha \mathcal{G}_{21}(t,s)\phi(x)\|_{L^q}\lesssim
(1+t)^\lambda\cdot
\Gamma^{\gamma_{p,q}+|\alpha|+1}(t,s)\cdot
(\|\phi\|_{L^{p}}^l+\|\partial_x^{|\alpha|+\omega_{r,q}}\phi\|_{L^{r}}^h),
\\
&\|\partial_x^\alpha \mathcal{G}_{22}(t,s)\phi(x)\|_{L^q}\lesssim
\Big(\frac{1+t}{1+s}\Big)^\lambda\cdot
\Gamma^{\gamma_{p,q}+|\alpha|}(t,s)\cdot
(\|\phi\|_{L^{p}}^l+\|\partial_x^{|\alpha|+\omega_{r,q}}\phi\|_{L^{r}}^h),
\end{align*}
where $\gamma_{p,q}:=n(1/{p}-1/q)$,
and $\omega_{r,q}>\gamma_{r,q}$ for $(r,q)\ne(2,2)$ and $\omega_{2,2}=0$.
Furthermore,
\begin{align*}
\|\partial_x^\alpha \mathcal{G}_{22}(t,s)\phi(x)\|_{L^q}\lesssim&
(1+t)^\lambda(1+s)^\lambda\cdot
\Gamma^{\gamma_{p,q}+|\alpha|+2}(t,s)\cdot
(\|\phi\|_{L^{p}}^l+\|\partial_x^{|\alpha|+1+\omega_{r,q}}\phi\|_{L^{r}}^h),
\\
\|\partial_x^\alpha \mathcal{G}_{22}(t,s)\phi(x)\|\lesssim&
\Big(\frac{1+t}{1+s}\Big)^\lambda\cdot
\Gamma^{\gamma_{p,q}+|\alpha|}(t,s)
\\
&\cdot
\Big(
(1+s)^{2\lambda}\cdot\Gamma^{2}(t,s)
+\frac{1}{(1+s)^{\lambda-1}}
+C_\kappa\Gamma^{\kappa}(t,s)
\Big)
\cdot
(\|\phi\|_{L^{p}}^l+\|\partial_x^{|\alpha|+1+\omega_{r,q}}\phi\|_{L^{r}}^h),
\end{align*}
where $\kappa\ge2$ can be chosen arbitrarily large
and $C_\kappa>0$ is a constant depending on $\kappa$.
\end{theorem}
{\it\bfseries Proof.}
These estimates are conclusions of Theorem \ref{th-linear-A}.
The last estimate is proved according to
\eqref{eq-linear-u-can} and the following inequality
$$
\Gamma^{-\kappa}(t,s)\cdot e^{-\varepsilon_u((1+t)^{1-\lambda}-(1+s)^{1-\lambda})}
\lesssim C_\kappa,
$$
since the super-exponential function decays faster than any algebraical decays.
$\hfill\Box$

\begin{corollary}
Let $(v(t,x),u(t,x))$ be the solution of the linear hyperbolic system \eqref{eq-vu}
(the third equation of $\boldsymbol w(t,x)$ is neglected as it decays super-exponentially)
corresponding to the initial data $(v(0,x),u(0,x))$.

For $q\in[2,\infty]$, $1\le p,r\le 2$, and $\lambda\in(-1,0)$, we have
\begin{align*}
\|\partial_x^\alpha v\|_{L^q}\approx &
(1+t)^{-\frac{1+\lambda}{2}(\gamma_{p,q}+|\alpha|)}
\cdot
\Big(\big\|v(0,\cdot)\big\|_{L^{p}}^l
+\big\|\partial_x^{|\alpha|+\omega_{r,q}}v(0,\cdot)\big\|_{L^{r}}^h
+\big\| u(0,\cdot)\big\|_{L^{p}}^l
+\big\|
\partial_x^{|\alpha|+\omega_{r,q}} u(0,\cdot)\big\|_{L^{r}}^h\Big),
\end{align*}
and
\begin{align*} \nonumber
\|\partial_x^\alpha u\|_{L^q}\approx &
(1+t)^{-\frac{1+\lambda}{2}(\gamma_{p,q}+|\alpha|)-\frac{1-\lambda}{2}}
\\
&\qquad\cdot
\Big(\big\|v(0,\cdot)\big\|_{L^{p}}^l
+\big\|\partial_x^{|\alpha|+1+\omega_{r,q}}v(0,\cdot)\big\|_{L^{r}}^h
+\big\| u(0,\cdot)\big\|_{L^{p}}^l
+\big\|\partial_x^{|\alpha|+1+\omega_{r,q}} u(0,\cdot)\big\|_{L^{r}}^h\Big).
\end{align*}

For $q\in[2,\infty]$, $1\le p,r\le 2$, and $\lambda=-1$, we have
\begin{align*}
\|\partial_x^\alpha v\|_{L^q}\approx &
|\ln(e+t)|^{-\frac{1}{2}(\gamma_{p,q}+|\alpha|)}
\cdot
\Big(\big\|v(0,\cdot)\big\|_{L^{p}}^l
+\big\|\partial_x^{|\alpha|+\omega_{r,q}}v(0,\cdot)\big\|_{L^{r}}^h
+\big\| u(0,\cdot)\big\|_{L^{p}}^l
+\big\|
\partial_x^{|\alpha|+\omega_{r,q}} u(0,\cdot)\big\|_{L^{r}}^h\Big),
\end{align*}
and
\begin{align*} \nonumber
\|\partial_x^\alpha u\|_{L^q}\approx &
(1+t)^{-1}\cdot
|\ln(e+t)|^{-\frac{1}{2}(\gamma_{p,q}+|\alpha|+1)}
\\
&\qquad\cdot
\Big(\big\|v(0,\cdot)\big\|_{L^{p}}^l
+\big\|\partial_x^{|\alpha|+1+\omega_{r,q}}v(0,\cdot)\big\|_{L^{r}}^h
+\big\| u(0,\cdot)\big\|_{L^{p}}^l
+\big\|\partial_x^{|\alpha|+1+\omega_{r,q}} u(0,\cdot)\big\|_{L^{r}}^h\Big).
\end{align*}
\end{corollary}

\begin{remark}
The general solutions of the
wave equation \eqref{eq-Pu} (satisfied by $u(t,x)$) decay optimally
faster than those solutions of \eqref{eq-Pv} (satisfied by $v(t,x)$);
while in the linear system \eqref{eq-vu},
$u(t,x)$ decays even faster.
\end{remark}

\begin{remark}
The decay estimate \eqref{eq-linear-u} for $u$ in the linear system \eqref{eq-vu}
derived from the optimal decay
estimate \eqref{eq-optimal-u} in Theorem \ref{th-wave}
is not optimal here since the initial data $u(0,x)=u_0(x)$ and
$\partial_t u(0,x)=\Lambda v_0(x)-\mu u_0(x)$ are not independent.
Cancellation occurs and the decay rate increases as in \eqref{eq-linear-u-opt}.
However, the estimate \eqref{eq-linear-u} is still of importance
in the decay estimates of the nonlinear system \eqref{eq-vbdu}
since the regularity required is one-order lower than in the estimate \eqref{eq-linear-u-opt}.
\end{remark}

\begin{proposition} \label{th-linear-st}
Let $(v(t,x),u(t,x))$ be the solution of the linear system \eqref{eq-vu}
corresponding to the initial data $(v(s,x),u(s,x))$
starting from time $s$.

For $q\in[2,\infty]$, $1\le p,r\le 2$, $\lambda\in[-1,0)$, and $t\ge s\ge T_0$,
where $T_0\ge0$ is the constant in Lemma \ref{le-Phi}, we have
\begin{align} \nonumber
\|\partial_x^\alpha v\|_{L^q}\lesssim &
\Gamma^{\gamma_{p,q}+|\alpha|}(t,s)\cdot
\Big(\big\|v(s,\cdot)\big\|_{L^{p}}^l
+\big\|\partial_x^{|\alpha|+\omega_{r,q}}v(s,\cdot)\big\|_{L^{r}}^h\Big)
\\ \label{eq-linear-v-st}
&+(1+s)^\lambda\cdot
\Gamma^{\gamma_{p,q}+|\alpha|+1}(t,s)\cdot
\Big(\big\| u(s,\cdot)\big\|_{L^{p}}^l
+\big\|
\partial_x^{|\alpha|+\omega_{r,q}} u(s,\cdot)\big\|_{L^{r}}^h\Big),
\end{align}
and
\begin{align} \nonumber
\|\partial_x^\alpha u\|_{L^q}\lesssim &
\Big(\frac{1+t}{1+s}\Big)^\lambda\cdot
\Gamma^{\gamma_{p,q}+|\alpha|}(t,s)\cdot
\Big(\big\|u(s,\cdot)\big\|_{L^{p}}^l
+\big\|\partial_x^{|\alpha|+\omega_{r,q}}u(s,\cdot)\big\|_{L^{r}}^h\Big)
\\ \label{eq-linear-u-st}
&+(1+t)^\lambda\cdot
\Gamma^{\gamma_{p,q}+|\alpha|+1}(t,s)\cdot
\Big(\big\|v(s,\cdot)\big\|_{L^{p}}^l
+\big\|\partial_x^{|\alpha|+\omega_{r,q}} v(s,\cdot)\big\|_{L^{r}}^h\Big),
\end{align}
where $\gamma_{p,q}:=n(1/{p}-1/q)$,
and $\omega_{r,q}>\gamma_{r,q}$ for $(r,q)\ne(2,2)$ and $\omega_{2,2}=0$.
Moreover, $u(t,x)$ decays faster if we assume one-order higher regularity as follows,
\begin{align} \nonumber
\|\partial_x^\alpha u\|_{L^q}\lesssim &
(1+t)^\lambda\cdot
\Gamma^{\gamma_{p,q}+|\alpha|+1}(t,s)\cdot
\Big(\big\|v(s,\cdot)\big\|_{L^{p}}^l
+\big\|\partial_x^{|\alpha|+1+\omega_{r,q}}v(s,\cdot)\big\|_{L^{r}}^h\Big)
\\ \label{eq-linear-u-opt-st}
&+(1+t)^\lambda(1+s)^\lambda\cdot
\Gamma^{\gamma_{p,q}+|\alpha|+2}(t,s)\cdot
\Big(\big\| u(s,\cdot)\big\|_{L^{p}}^l
+\big\|\partial_x^{|\alpha|+1+\omega_{r,q}} u(s,\cdot)\big\|_{L^{r}}^h\Big).
\end{align}

The decay estimate \eqref{eq-linear-v-st} is element-by-element optimal for all
$\frac{t}{2}\ge s\ge T_0$;
the decay estimate \eqref{eq-linear-u-opt-st} is optimal with respect to $v(s,x)$
for all $\frac{t}{2}\ge s\ge T_0$;
the decay estimates \eqref{eq-linear-v-st} and \eqref{eq-linear-u-opt-st} are optimal
for all $t\ge s\ge0$.
\end{proposition}
{\it\bfseries Proof.}
This is proved in a similar way as Proposition 2.2 in \cite{Ji-Mei-1}
based on the optimal decay estimates of the linear wave equations in Theorem \ref{th-wave}
and the optimal decay estimates of the Fourier multiplies in Lemma \ref{le-Phi}.
$\hfill\Box$

We improve the decay estimates \eqref{eq-linear-u-st} on
$\|\partial_x^\alpha u(t,\cdot)\|_{L^q}$ in Proposition \ref{th-linear-st}
by taking advantage of the cancellation between the initial data
$u(s,x)$ and $\partial_t u(s,x)=\Lambda v(s,x)-b(s)u(s,x)$
if we regard $u(t,x)$ as a solution of the wave equation \eqref{eq-Pu}.

\begin{proposition}[Decay rates improved by cancellation]
\label{th-linear-st-can}
Let $(v(t,x),u(t,x))$ be the solution of the linear system \eqref{eq-vu}
corresponding to the initial data $(v(s,x),u(s,x))$
starting from the time $s$.
Then for $q\in[2,\infty]$, $1\le p,r\le 2$, and $\lambda\in[-1,0)$,
and for $t\ge s\ge T_0$ ($T_0\ge0$ is the constant in Lemma \ref{le-Phi}), we have
\begin{align} \nonumber
\|\partial_x^\alpha u(t,\cdot)\|_{L^q}
\lesssim &
(1+t)^\lambda\cdot
\Gamma^{\gamma_{p,q}+|\alpha|+1}(t,s)
\cdot
\Big(\big\|v(s,\cdot)\big\|_{L^{p}}^l
+\big\|\partial_x^{|\alpha|+\omega_{r,q}}v(s,\cdot)\big\|_{L^{r}}^h\Big)
\\ \nonumber
&+(1+t)^\lambda(1+s)^\lambda\cdot
\Gamma^{\gamma_{p,q}+|\alpha|+2}(t,s)
\cdot
\Big(\big\| u(s,\cdot)\big\|_{L^{p}}^l
+\big\|\partial_x^{|\alpha|+\omega_{r,q}} u(s,\cdot)\big\|_{L^{r}}^h\Big)
\\ \nonumber
&+\Big(\frac{1+t}{1+s}\Big)^\lambda\cdot
\Gamma^{\gamma_{p,q}+|\alpha|}(t,s)
\cdot
\Big(\frac{1}{(1+s)^{1-\lambda}}
+e^{-\varepsilon_u((1+t)^{1-\lambda}-(1+s)^{1-\lambda})}\Big)
\\ \label{eq-linear-u-st-can}
&\qquad\cdot
\Big(\big\| u(s,\cdot)\big\|_{L^{p}}^l
+\big\|\partial_x^{|\alpha|+\omega_{r,q}} u(s,\cdot)\big\|_{L^{r}}^h\Big),
\end{align}
where $\varepsilon_u>0$ is the constant in the definition of different zones in the
phase-time space, $\gamma_{p,q}:=n(1/{p}-1/q)$,
and $\omega_{r,q}>\gamma_{r,q}$ for $(r,q)\ne(2,2)$ and $\omega_{2,2}=0$.
The decay estimate \eqref{eq-linear-u-st-can} is optimal
with respect to $v(s,x)$ for all $\frac{t}{2}\ge s\ge T_0$.
\end{proposition}
{\it\bfseries Proof.}
The outline of the proof is similar to Proposition 2.3 in \cite{Ji-Mei-1}.
Here we omit the details.
$\hfill\Box$

{\it\bfseries Proof of Theorem \ref{th-linear-A}.}
The optimal decay estimates \eqref{eq-linear-v} and \eqref{eq-linear-u-opt} are proved
in Proposition \ref{th-linear-st}
and the decay estimate \eqref{eq-linear-u-can} improved by cancellation
is proved in Proposition \ref{th-linear-st-can}.
$\hfill\Box$

\end{appendices}

\

{\bf Acknowledgement}.
This work was done when the first author visited McGill University
supported by China Scholarship Council (CSC) for
the senior visiting scholar program.
He would like to express his sincere thanks for the hospitality
of McGill University and CSC.
The research of the first author was supported by NSFC Grant No.~11701184 and CSC No. 201906155021.
The research of the second author was supported in part
by NSERC Grant RGPIN 354724-16, and FRQNT Grant No. 2019-CO-256440.


\begin{thebibliography}{10}

\bibitem {Burq-Raugel-Schlag-2015}
R. Burq, G. Raugel, and W. Schlag,
Long time dynamics for damped Klein-Gordon equations,
{\it Ann. Sci. \'Ec. Norm. Sup\'er.},  {\bf 50} (2015),  1447--1498.

\bibitem {Burq-Raugel-Schlag-2018}
R. Burq, G. Raugel, and W. Schlag,
Long time dynamics for weakly damped nonlinear Klein-Gordon equations,
 arXiv: 1801.06735v1.


\bibitem{C-D-S-W} G.-Q. Chen, C. Dafermos, M. Slemrod, and D. Wang,
On two-dimensional sonic-subsonic flow,
{\it Commun. Math. Phys.}, {\bf 271} (2007),  635--647.

\bibitem{Chen-Pan} G. Chen, R.  Pan, and S.  Zhu,
Singularity formation for the compressible Euler equations,
{\it SIAM J. Math. Anal.}, {\bf 49} (2017),  2591--2614.

\bibitem{Chen-Li-Li-Mei-Zhang}  S. G. Chen, H. Li, J. Li, M. Mei, and K. Zhang,
Global and blow-up solutions to compressible Euler equations with time-dependent damping,
{\it J. Differential Equations}, {\bf 268} (2020), 5035--5077.

\bibitem{Courant-F} R. Courant and O.K. Friedrichs,
Supersonic Flow and Shock Waves, Springer-Verlag, New York, 1948.

\bibitem{Cui-Yin-Zhang-Zhu} H.-B. Cui, H.-Y. Yin, J.-S. Zhang, and C.-J. Zhu,
Convergence to nonlinear diffusion waves for solutions of Euler equations
with time-depending damping,
{\it J. Differential Equations}, {\bf 264}  (2018), 4564--4602.


\bibitem{Dafermos} C. Dafermos,
Hyperbolic Conservation Laws in Continuum Physics, 3rd ed., Springer-Verlag, New York, 2010.

\bibitem{Geng-Huang} S. Geng and F. Huang,
$L^1$-convergence rates to the Barenblatt solution for the damped compressible Euler equations,
{\it J. Differential Equations}, {\bf 266} (2019), 7890--7908.

\bibitem{Geng-Lin-Mei} S. Geng, Y. Lin, and M. Mei,
Asymptotic behavior of solutions to Euler equations with time-dependent damping in critical case,
{\it SIAM J. Math. Anal.},  {\bf 52}  (2020),  1463--1488.

\bibitem{Guo}
Y. Guo and B. Pausader,
Global smooth ion dynamics in the Euler-Poisson system,
{\it Commun. Math. Phys.}, {\bf 303} (2011), 89--125.

\bibitem{Hou-Witt-Yin} F. Hou, I. Witt, and H.C. Yin,
Global existence and blowup of smooth solutions of 3-D potential equations with time-dependent damping, {\it Pacific J. Math.}, {\bf 292} (2018), 389--426.

\bibitem{Hou-Yin} F. Hou and H.C. Yin,
On the global existence and blowup of smooth solutions to the multi-dimensional
compressible Euler equations with time-depending damping,
{\it Nonlinearity}, {\bf 30} (2017), 2485--2517.


\bibitem{Hsiao-Liu} L. Hsiao and T.-P. Liu,
Convergence to diffusion waves for solutions of a system of hyperbolic
conservation laws with damping,
{\it Commun. Math. Phys.}, {\bf 143} (1992), 599--605.

\bibitem{Ji-Mei-1} S. Ji and M. Mei, Optimal decay rates of the compressible Euler equations
with time-dependent damping in $\mathbb R^n$: (I) under-damping case, preprint, 2020.

\bibitem{Huang-Marcati-Pan} F.M. Huang, P. Marcati, and R.H. Pan,
Convergence to the Barenblatt solution for the
compressible Euler equations with damping and vacuum,
{\it Arch. Ration. Mech. Anal.}, {\bf 176} (2005), 1--24.

\bibitem{Huang-Pan} F.M. Huang and R. H. Pan,
Convergence rate for compressible Euler equations with damping and vacuum,
{\it Arch. Ration. Mech. Anal.}, {\bf 166} (2003), 359--376.

\bibitem{Huang-Pan-Wang} F.M. Huang, R. Pan and Z. Wang,
$L^1$ convergence to the Barenblatt solution for compressible Euler equations with damping,
{\it Arch. Ration. Mech. Anal.}, {\bf 200} (2011),  665--689.

\bibitem{Lax} P.D. Lax,
Development of singularities of solutions of nonlinear hyperbolic partial differential equations, {\it J. Math. Phys.}, {\bf 5} (1964), 611--614.

\bibitem{Li-Li-Mei-Zhang} H. Li, J. Li, M. Mei, and K. Zhang,
Convergence to nonlinear diffusion waves for solutions of p-system with time-dependent damping,
{\it J. Math. Anal. Appl.},   {\bf 456} (2017),  849--871.

\bibitem{Hailiang-Li} H.-L. Li and X. Wang,
Formation of singularities of spherically symmetric solutions to the 3D compressible Euler equations and Euler-Poisson equations,
{\it Nonlinear Differential Equations Appl.}, {\bf  25} (2018), 1--15.

\bibitem{Luo-Zeng} T. Luo and H.H. Zeng,
Global existence of smooth solutions and convergence to
Barenblatt solutions for the physical vacuum free boundary problem of
compressible Euler equations with damping,
{\it Comm. Pure Appl. Math.}, {\bf 69} (2016), 1354--1396.

\bibitem{Marcati-Milani} P. Marcati and A. Milani,
The one-dimensional Darcy's law as the limit of a compressible Euler flow,
{\it J. Differential Equations}, {\bf 84} (1990), 129--147.

\bibitem{Mei} M. Mei,
Best asymptotic profile for hyperbolic $p$-system with damping,
{\it SIAM J. Math. Anal.}, {\bf 42} (2010), 1--23.

\bibitem{Nishihara} K. Nishihara,
Convergence rates to nonlinear diffusion waves for solutions of system of hyperbolic
conservation laws with damping,
{\it J. Differential Equations}, {\bf 131} (1996), 171--188.

\bibitem{Nishihara-Wang-Yang} K. Nishihara, W. K. Wang, and T. Yang,
$L_p$-convergence rates to nonlinear diffusion waves for $p$-system with damping,
{\it  J. Differential Equations}, {\bf  161} (2000), 191--218.

\bibitem{Pa1} X. Pan,
Blow up of solutions to $1$-d Euler equations with time-dependent damping,
{\it  J. Math. Anal. Appl.},  {\bf 442} (2016),   435--445.

\bibitem{Pa2} X. Pan,
Global existence of solutions to $1$-d Euler equations with time-dependent damping,
{\it  Nonlinear Anal.}, {\bf   132} (2016), 327--336.

\bibitem{PanXH-AA}
X. Pan,
Global existence and asymptotic behavior of solutions to the Euler equations with time-dependent damping,
{\it Applicable Analysis}, (2020), 1--30.

\bibitem{Sideris-Thomas-Wang} T. Sideris, B. Thomases, and D. Wang,
Long time behavior of solutions to the 3D compressible Euler equations with damping,
{\it  Comm. Partial Differential Equations}, {\bf  28} (2003),  795--816.

\bibitem{Smoller} J. Smoller,
Shock Waves and Reaction-Diffusion Equations, Springer-Verlag, New York, 1982.

\bibitem{Sugiyama1} Y. Sugiyama,
Singularity formation for the 1D compressible Euler equations with variable damping coefficient, {\it  Nonlinear Anal.},  {\bf 170} (2018), 70--87.

\bibitem{Sugiyama2} Y. Sugiyama,
Remark on the global existence for the 1D
compressible Euler equation with time-dependent damping, arXiv: 1909.05683.



\bibitem{TanZ-JDE13}
Z. Tan and Y. Wang,
Global solution and large-time behavior of the $3$D compressible Euler equations with damping,
{\it J. Differential Equations}, {\bf 254} (2013), 1686--1704.

\bibitem{TanZ-JDE12}
Z. Tan and G. Wu,
Large time behavior of solutions for compressible Euler equations with damping in $\mathbb R^3$,
{\it J. Differential Equations}, {\bf 252} (2012), 1546--1561.

\bibitem{Tao}
T. Tao,
Nonlinear dispersive equations, local and global analysis,
CBMS. Regional Conference Series in Mathematics,
{\bf 106}.
Published for the Conference Board of the Mathematical Science,
Washington, DC; Providence, RI: Amer. Math. Soc., 2006.


\bibitem{Todoraova}
G. Todorova and B. Yordanov,
Weighted $L^2$-estimates for dissipative wave equations with variable coefficients,
{\it J. Differential Equations}, {\bf 246} (2009), 4497--4518.

\bibitem{Wang-Chen} D. Wang and G.-Q. Chen,
Formation of singularities in compressible Euler-Poisson fluids with heat diffusion and damping relaxation, {\it  J. Differential Equations}, {\bf  144} (1998), 44--65.

\bibitem{Wirth-JDE06}
J. Wirth,
Wave equations with time-dependent dissipation I Non-effective dissipation,
{\it J. Differential Equations}, {\bf  222} (2006), 487--514.

\bibitem{Wirth-JDE07}
J. Wirth,
Wave equations with time-dependent dissipation II Effective dissipation,
{\it J. Differential Equations}, {\bf  232} (2007), 74--103.

\bibitem{Wirth-MMAS04}J. Wirth, Solution representations for a wave equation with weak dissipation,
{\it Math. Methods Appl. Sci.}, {\bf 27} (2004), 101--124.


\end{thebibliography}
\end{document}